\documentclass[10pt,a4paper]{article}
\usepackage{amssymb,amsmath,amsthm}
\usepackage{enumitem}
\usepackage[nosort]{cite}
\usepackage{pict2e,graphicx}
\usepackage{color}
\usepackage{comment}
\graphicspath{{./figs/}}
\numberwithin{equation}{section}
\newtheorem{thm}{Theorem}[section]
\newtheorem{df}[thm]{Definition}
\newtheorem{prop}[thm]{Proposition}
\newtheorem{lem}[thm]{Lemma}
\newtheorem{rem}[thm]{Remark}

\newtheorem{cor}[thm]{Corollary}
\newtheorem{ex}[thm]{Example}
\let\oldproofname=\proofname
\renewcommand{\proofname}{\rm\bf{\oldproofname}}
\newenvironment{Proof}[1][\unskip]%
 {\begin{trivlist} \item[]{\bf Proof #1. }}%
 {\hspace*{\fill}$\rule{.4\baselineskip}{.4\baselineskip}$\end{trivlist}}

\newcommand{\N}{\mathbb{N}}

\newcommand{\R}{\mathbb{R}}

\newcommand{\cN}{\mathcal{N}}
\newcommand{\cO}{\mathcal{O}}

\newcommand{\cS}{\mathcal{S}}
\newcommand{\cT}{\mathcal{T}}

\newcommand{\dd}{\,{\rm d}}
\newcommand{\D}{{\rm d}}

\newcommand{\supp}{\mathop{\mathrm{supp}}}
\newcommand{\dist}{\mathop{\mathrm{dist}}}
\newcommand{\erfc}{\mathop{\mathrm{erfc}}}

\newcommand{\loc}{\mathrm{loc}}

\newcommand{\QED}{\mbox{}\hfill$\Box$}
\renewcommand{\:}{\thinspace :}
\renewcommand{\epsilon}{\varepsilon}

\DeclareMathOperator*{\esssup}{ess\,sup}
\DeclareMathOperator*{\essinf}{ess\,inf}
\DeclareSymbolFont{bbold}{U}{bbold}{m}{n}
\DeclareSymbolFontAlphabet{\mathbbold}{bbold}
\newcommand{\id}{\mathbbold{1}}

\begin{document}

\title{Viscous shocks and long-time behavior of scalar conservation laws}

\author{Thierry Gallay and Arnd Scheel}


\maketitle


\begin{abstract}
We study the long-time behavior of scalar viscous conservation laws via the
structure of $\omega$-limit sets. We show that $\omega$-limit sets always
contain constants or shocks by establishing convergence to shocks for arbitrary
monotone initial data. In the particular case of Burgers' equation, we review
and refine results that parametrize entire solutions in terms of probability
measures, and we construct initial data for which the $\omega$-limit set is
not reduced to the translates of a single shock. Finally we propose several
open problems related to the description of long-time dynamics. 
\end{abstract}

\section{Introduction and main results}\label{s:1}

We are interested in the long-time dynamics of viscous scalar conservation laws,
\begin{equation}\label{SCL}
  \partial_t u(t,x) + f'(u(t,x))\partial_x u(t,x) \,=\, \partial_x^2 u(t,x)\,,
  \qquad t > 0\,, \quad x \in \R\,, 
\end{equation}
with smooth and strictly convex flux function $f : \R \to \R$, that is,
$f''(u) > 0$ for all $u \in \R$.  A typical example is Burgers' equation where
$f(u) = u^2/2$. The Cauchy problem for \eqref{SCL} is globally well-posed in the
space $L^\infty(\R)$, see e.g. \cite{Se2}.  More precisely, given initial data
$u_0 \in L^\infty(\R)$, equation \eqref{SCL} has a unique global solution
$u \in C^0((0,+\infty),L^\infty(\R))$ such that $u(t,\cdot)$ converges to $u_0$
in the weak-$*$ topology of $L^\infty(\R)$ as $t \to 0+$.  By parabolic
regularity, the function $u(t,x)$ is smooth for all positive times.  For any
$t > 0$, let $\cS_t : L^\infty(\R) \to L^\infty(\R)$ be the nonlinear map
defined by $u(t,\cdot) = \cS_t(u_0)$, where $u(t,x)$ is the solution of
\eqref{SCL} with initial data $u_0 \in L^\infty(\R)$. We also write
$\cS_0 = \id$, the identity map.

To a first approximation, the long-time behavior of $\cS_t(u_0)$ as $t\to\infty$ is
described by the collection of all limit points, usually referred
to as the $\omega$-limit set. The unboundedness of the spatial domain $\R$
implies a typical lack of compactness of the trajectory $\{\cS_t(u_0)\,|\,t>0\}$,
and the $\omega$-limit set may indeed be empty when convergence is measured in
the uniform topology defined by the norm in $L^\infty(\R)$. It is therefore
preferable to rely on the local topology induced by $L^\infty_\loc(\R)$, which is the
topology of uniform convergence on compact intervals $[-R,R] \subset \R$. The
$\omega$-limit set is commonly defined as follows:
\begin{equation}\label{e:omega0}
  \omega_0(u_0) \,:=\, \bigl\{v \in L^\infty(\R)\,\big|\, \exists\,
  t_k\to\infty \,\text{ s.th. } \cS_{t_k}(u_0) \to v
  \text{ in } L^\infty_\loc(\R)\bigr\}\,,
\end{equation}
but this definition assigns a particular role to the laboratory frame and is not invariant
under Galilean transformations. A somewhat richer description of asymptotic behavior is
obtained by considering the set of limit points modulo translations,
\begin{equation}\label{e:omega}
  \omega(u_0) \,:=\, \bigl\{v \in L^\infty(\R)\,\big|\, \exists\, t_k\to\infty
  \text{ and }x_k\in\R \,\text{ s.th. } \cT_{x_k}\cS_{t_k}(u_0) \to v \text{ in }
  L^\infty_\loc(\R)\bigr\}\,,
\end{equation}
where $(\cT_y u)(x)= u(x-y)$. Note that we use the zero subscript in the
definition of $\omega_0$ in \eqref{e:omega0} to emphasize the fixed origin in
the definition of locally uniform convergence in \eqref{e:omega0}.

Fairly standard results assert that both $\omega_0(u_0)$ and $\omega(u_0)$ are
non-empty, compact, connected, fully invariant, attractive, and chain recurrent
(up to translations) in the topology of $L^\infty_\loc(\R)$;
see Propositions~\ref{omega-prop0} and \ref{omega-prop}. Full invariance implies
in particular that for any $v_0$ in the $\omega$-limit set, there exists a solution
$v(t,x)$ of \eqref{SCL} that is defined for all $t\in\R$ and satisfies
$v(0,\cdot)=v_0$. We refer to such solutions, defined for all positive and
negative times, as {\em entire solutions}. Describing all possible long-term
dynamics can then be rephrased as describing all subsets of the family of entire
solutions that can occur as $\omega$-limit sets for bounded initial data. 

The results we present can be seen as small steps in this direction. Somewhat
trivial candidates for the $\omega$-limit sets are first spatially constant
states $v(x)\equiv m$ and then viscous shocks, found as traveling-wave solutions
$v(t,x)=\phi_{\beta,\alpha}(x-ct)$ with $\phi(-\infty)=\beta>\phi(+\infty)=\alpha$,
$\phi'(\xi)<0$ for all $\xi$, and $c$ given by the Rankine-Hugoniot formula
$c_{\beta,\alpha}=(f(\beta)-f(\alpha))/(\beta-\alpha)$. It is known since the
classical work of Il'in and Oleinik \cite{IO} that large sets of initial
data give rise to solutions of \eqref{SCL} that converge uniformly to
shocks as $t \to +\infty$, see also \cite{FS}. In Proposition~\ref{mono-prop}
below, we show that this is the case for all initial data that are
monotonically decreasing, without any assumption on the rate at which the
limits are approached as $x \to \pm\infty$. At this level of generality, we cannot prove
convergence to a fixed translate of the shock as $t \to +\infty$. In fact, as
discussed in Remark~\ref{mono-rem}, there exist monotone initial data
$u_0 \in L^\infty(\R)$ with $u_0(-\infty) = \beta>u_0(+\infty)=\alpha$ and
$c_{\beta,\alpha}=0$ such that, for instance, $\omega_0(u_0) = \beta$; in particular
one observes that $\cT_y\phi_{\beta,\alpha}\notin \omega_0(u_0)$ for all $y \in \R$.

Our first general result establishes a property reminiscent of the
Poincar\'e-Bendixson theorem, in the sense that it describes the long-time
behavior of solutions with initial data in $\omega$-limit sets. 

\begin{prop}\label{p:shock}
For every $u_0\in L^\infty(\R)$ and any nonconstant $v \in \omega(u_0)$, there exist
real numbers $\alpha < \beta$ such that 
\begin{equation}\label{omvshock}
  \omega(v) \,=\, \overline{\{\cT_y\phi_{\beta,\alpha}\,;\,y\in\R\}}^{\,L^\infty_\loc}
  \,\subset\, \omega(u_0)\,.
\end{equation}
In particular, the set $\omega(u_0)$ contains a shock unless it consists entirely of
constants.
\end{prop}

In other words, if $v \in \omega(u_0)$ is nonconstant, the $\omega$-limit set
$\omega(v)$ consists of all translates of a viscous shock $\phi_{\beta,\alpha}$,
together with the constant states $\alpha$ and $\beta$ that arise as limits of
the shock profile at $\pm\infty$. The proof relies on the simple observation
that any such $v$ is necessarily monotonically decreasing, as a consequence
of Oleinik's inequality \eqref{Oleinik}. We can thus invoke
Proposition~\ref{mono-prop} to determine the $\omega$-limit $\omega(v)$,
which is included in $\omega(u_0)$ since the latter set is invariant under
the dynamics of \eqref{SCL}. 

Our remaining results focus on the specific case of Burgers' equation, which
through the Cole-Hopf transformation allows for a somewhat explicit
representation of any solution in terms of its initial data. Interestingly, as
pointed out in \cite{UPK}, bounded entire solutions of Burgers' equation
can be represented in terms of probability measures $\mu$ on the real line,
\begin{equation}\label{UPK-rep0}
  u(t,x) \,=\, \genfrac{}{}{1pt}{0}{\int z\,e^{-zx/2 + z^2t/4} \dd\mu(z)}{
  \int e^{-zx/2 + z^2t/4} \dd\mu(z)}\,, \qquad t \in \R\,, \quad x \in \R\,.
\end{equation}
This remarkable formula gives, in particular, an explicit characterization
of candidates for elements in $\omega$-limit sets. In Section~\ref{s:5} we give a
short proof of the representation \eqref{UPK-rep0}, showing that the
measure $\mu$ is unique and supported in the closure of the range of the entire
solution $u$. We also relate the measure $\mu$ to backward-in-time asymptotics
of the entire solution $u$.  A striking result in this direction is:

\begin{prop}\label{p:suppmu}
Assume that $u$ is given by \eqref{UPK-rep0} for some probability measure
$\mu$ on $\R$. A real number $c \in \R$ belongs to $\supp(\mu)$ if and only
if $u(t,\cdot+ct)$ converges to $c$ in $L^\infty_\loc(\R)$ as $t \to -\infty$. 
\end{prop}

It is also possible to determine the asymptotic behavior of $u(t,x)$
as $t \to -\infty$ in Galilean frames with speeds $c \notin \supp(\mu)$.
In that case we define
\[
  m_-(c) \,=\, \sup\{z<c\,|\,z\in \supp\,\mu\}\,, \quad \text{ and }\quad
  m_+(c) \,=\, \inf\{z>c\,|\,z\in \supp\,\mu\}\,.
\]

\begin{prop}\label{p:self-similar}
If $c\notin\supp(\mu)$ and $c\neq (m_+(c)+m_-(c))/2$, the solution $u$ defined
by \eqref{UPK-rep0} satisfies
\[
  \lim_{t\to -\infty} u(t,\cdot + ct) \,=\,
  \left\{\begin{array}{ll}
          m_-(c) & \text{if } ~c<(m_+(c)+m_-(c))/2\,,\\[1mm] 
          m_+(c) & \text{if } ~c>(m_+(c)+m_-(c))/2\,,
        \end{array}\right.
\]
where convergence is understood in $L^\infty_\loc(\R)$.
\end{prop}

We refer to Propositions~\ref{mu-prop1} and \ref{mu-prop2} below for more
general statements, which also cover the somewhat delicate situation where
$c =(m_+(c)+m_-(c))/2$. Other properties of the measure $\mu$, such as
the presence of atoms, can also be detected in the ancient behavior of
the corresponding entire solution $u$. 

Lastly, we show that out of this plethora of entire solutions, the
$\omega$-limit set may contain elements that are not simply shocks or constants.

\begin{prop}\label{p:merger}
There exist initial data $u_0\in L^\infty(\R)$ for Burgers' equation such that
$\omega_0(u_0)$ contains a solution $v(t,x)$ that is neither a constant nor a shock.
In fact $v$ describes the merging of a pair of shocks into a single shock.  
\end{prop}

The construction is carried out in a somewhat explicit fashion in Section~\ref{s:6}.
In the terminology of dynamical systems, the $\omega$-limit set contains a heteroclinic
trajectory connecting the zero solution to a steady shock $\phi$, as well as a
continuous family of steady shocks interpolating between $\phi$ and $0$. This can
be compared to a famous example of coarsening dynamics due to Eckmann and Rougemont
\cite{ER}, also rigorously studied by Pol\'a\v{c}ik \cite{Po1,Po2}, where the
$\omega$-limit set is a heteroclinic loop. 

\medskip\noindent
\textbf{Outline.} We recall basic properties of conservation laws and
shocks in Section~\ref{s:2}. We then formulate and establish properties of
both $\omega$-limit sets $\omega_0(u_0)$ and $\omega(u_0)$ in Section~\ref{s:3}.
Our first main result, the convergence to shocks for monotone initial data,
is proved in Section~\ref{s:4}. Section~\ref{s:5} derives the representation
of entire solutions in terms of probability measures, displays some key
examples, and relates measures to ancient limits. Lastly, Section~\ref{s:6}
is devoted to the proof of Proposition~\ref{p:merger}. We conclude with a brief
discussion.

\medskip\noindent
\textbf{Acknowledgements.} This project started from discussions held in the
stimulating atmosphere of the Mathematisches Forschungsinstitut Oberwolfach,
in August 2021. ThG would like to thank Denis Serre for his expert advice on
several points addressed in this work. The authors were partially
supported by the grants ISDEEC ANR-16-CE40-0013 (ThG) and NSF DMS-1907391,
DMS-2205663 (AS).

\section{Properties of scalar conservation laws and shocks solutions}\label{s:2}

We first recall some basic properties of scalar conservation laws of the form \eqref{SCL}.

\vspace{-7pt}
\paragraph{A priori bounds and monotonicity.} The evolution semigroup
$(\cS_t)_{t \ge 0}$ defined by \eqref{SCL} in $L^\infty(\R)$ has the
following properties\: 

\begin{enumerate}[leftmargin=20pt,itemsep=2pt,topsep=2pt,label=\alph*)]

\item {\em Monotonicity\:} if $u_0,u_1 \in L^\infty(\R)$ and $u_0 \le u_1$
almost everywhere, then $\cS_t(u_0) \le \cS_t(u_1)$ everywhere when $t > 0$;

\item {\em Contraction in $L^1$\:} if $u_0,u_1 \in L^\infty(\R)$ satisfy
$u_0 - u_1 \in L^1(\R)$, then $\cS_t(u_0) - \cS_t(u_1) \in L^1(\R)$ and 
$\|\cS_t(u_0) - \cS_t(u_1)\|_{L^1} \le \|u_0 - u_1\|_{L^1}$ for all $t > 0$;

\item {\em Conservation of mass\:} under the assumptions of b), we also
have
\[
  \int_\R \bigl(\cS_t(u_0) - \cS_t(u_1)\bigr)(x)\dd x \,=\,
  \int_\R \bigl(u_0 - u_1\bigr)(x)\dd x\,, \qquad t > 0\,. 
\]  
\end{enumerate}

Assertions a), b), c) are readily established using the parabolic
maximum principle \cite{PW} and the fact that \eqref{SCL} is a
conservation law, see e.g. \cite{Se1,Se2}. 

Another remarkable property
of the solutions of \eqref{SCL} is a universal upper bound for the
derivative $\partial_x u$, which is known as Oleinik's
inequality. Given $u_0 \in L^\infty(\R)$, we define
\begin{equation}\label{alphabet}
  \alpha \,:=\, \essinf_{x \in \R} u_0(x)\,, \qquad \beta \,:=\,
  \esssup_{x \in \R} u_0(x) \,. 
\end{equation}
Since constants are steady states of \eqref{SCL}, monotonicity implies that
the solution $u(t) = \cS_t(u_0)$ satisfies $\alpha \le u(t,x) \le \beta$ for
all $t > 0$ and all $x \in \R$ (in fact, due to the strong maximum principle,  
both inequalities are strict as soon as $\alpha < \beta$). Oleinik's
inequality asserts that, for all $t > 0$ and all $x \in \R$, 
\begin{equation}\label{Oleinik}
  \partial_x u(t,x) \,<\, \frac{1}{kt}\,, \qquad\hbox{where}\quad
   k \,:=\, \min\bigl\{f''(u)\,;\, u \in [\alpha,\beta]\bigr\} \,>\, 0\,.
\end{equation}
For convenience, we include a short proof of \eqref{Oleinik} in
Section~\ref{ssecA1}.

\vspace{-7pt}
\paragraph{Viscous shocks.}
Given $\alpha,\beta \in \R$ with $\alpha < \beta$, equation \eqref{SCL} has
a unique traveling wave solution of the form $u(t,x) = \phi_{\beta,\alpha}(x-ct)$,
such that $\phi (-\infty) = \beta$ and $\phi(+\infty) = \alpha$, and $\phi(0) =
(\alpha+\beta)/2$. The profile $\phi$ is strictly decreasing and solves  
\begin{equation}\label{phiODE}
  \phi'(y) \,=\, f(\phi(y)) - c\phi(y) - d\,, \qquad y \in \R\,,
\end{equation}
where 
\begin{equation}\label{RH}
  c \,:=\, \frac{f(\beta) - f(\alpha)}{\beta-\alpha}\,, \qquad
  d \,:=\, f(\beta) - c\beta \,\equiv\, f(\alpha) - c\alpha\,,
\end{equation}
Strict convexity of $f$ gives the Lax condition $f'(\beta) > c > f'(\alpha)$, and the
ODE \eqref{phiODE} then implies that $\phi_{\beta,\alpha}$ converges exponentially to
its limits at $\pm\infty$.

Stability of viscous shocks has been known since the classical work of Il'in and
Oleinik \cite{IO}. For localized perturbations, that is, for initial data
$u_0 \in L^\infty(\R)$ with $u_0 - \phi_{\beta,\alpha} \in L^1(\R)$ for some
$\alpha < \beta$, the solution $u(t,x)$ of \eqref{SCL} converges uniformly to
$\phi_{\beta,\alpha}(x-ct-x_0)$ as $t \to +\infty$, with $c$ as in \eqref{RH} and
\begin{equation}\label{x0shift}
  x_0 \,:=\, \frac{1}{\beta-\alpha}\int_\R \bigl(u_0(x) -
  \phi_{\beta,\alpha}(x)\bigr)\dd x\,.
\end{equation}
Extensions towards viscous conservation laws with more general flux function,
allowing for degenerate shocks, can be found in the references
\cite{MN,FS,Se2,NZ,HST,Hen}. Rates of convergence can be obtained under
stronger localization of the perturbations.  However, the hypothesis that
$u_0 - \phi_{\beta,\alpha} \in L^1(\R)$, which allows one to determine the
asymptotic shift \eqref{x0shift}, seems to play an important role in all
existing results. Our analysis in Section \ref{s:4} removes this restriction
for monotone solutions. 

\section{Properties of $\omega$-limit sets}\label{s:3}

In this section we establish the properties of the $\omega$-limit sets
\eqref{e:omega0}, \eqref{e:omega} that were announced in the introduction.
Here and in what follows, we denote by $d$ the distance on $L^\infty(\R)$
defined by
\begin{equation}\label{distdef}
  d(u,v) \,=\, \|u-v\|_{\exp}\,, \qquad \hbox{where}\quad \|u\|_{\exp} \,=\,
  \esssup_{x \in \R} \bigl(e^{-|x|} |u(x)|\bigr)\,.
\end{equation}
As is easily verified, on any bounded set $\Sigma \subset L^\infty(\R)$, the topology
defined by the distance \eqref{distdef} coincides with the topology of $L^\infty_\loc(\R)$,
namely the topology of uniform convergence on compact subsets of $\R$.

In view of the properties recalled in Section~\ref{s:2}, for any initial
data $u_0 \in L^\infty(\R)$ the solution $u(t,\cdot) = \cS_t(u_0)$ of
\eqref{SCL} belongs for all times to the ball 
\begin{equation}\label{Sigmau}
  \Sigma(u_0) \,:=\, \bigl\{u \in L^\infty(\R)\,;\, \|u\|_{L^\infty} \le
  \|u_0\|_{L^\infty} \bigr\} \,\subset\, L^\infty(\R)\,.
\end{equation}
The following standard result plays a fundamental role: 

\begin{lem}\label{contlem}
When equipped with the topology of $L^\infty_\loc(\R)$, the  ball $\Sigma(u_0)$ defined
by \eqref{Sigmau} is closed and the solution map $\cS_t:\Sigma(u_0)\to \Sigma(u_0)$
is continuous for any $t\ge 0$. 
\end{lem}

\begin{Proof}
It is easy to check that $\Sigma(u_0)$ is closed in $L^\infty_\loc(\R)$, and the
properties recalled in Section~\ref{s:2} imply that the semiflow $\cS_t$ maps
the ball $\Sigma(u_0)$ into itself. The key point is the continuous
dependence of the solution $\cS_t(u)$ upon the initial data $u \in \Sigma(u_0)$,
in the topology of $L^\infty_\loc(\R)$. This a rather standard result for
parabolic PDEs on unbounded domains, see e.g. \cite{MS}.  For the reader's
convenience, the argument showing continuity is reproduced in Section~\ref{ssecA2}
below.
\end{Proof}

We are now in position to establish the main properties of the $\omega$-limit
set \eqref{e:omega0}. 

\begin{prop}\label{omega-prop0}
For any $u_0 \in L^\infty(\R)$, the $\omega$-limit set $\omega_0(u_0)$
defined by \eqref{e:omega0} is bounded in $L^\infty(\R)$ and, when equipped
with the topology of $L^\infty_\loc(\R)$, has the following properties\:

\begin{enumerate}[leftmargin=20pt,itemsep=2pt,topsep=0pt,label=\alph*)]

\item $\omega_0(u_0)$ is non-empty, compact, connected, and
\begin{equation}\label{omega0exp}
  \omega_0(u_0) \,=\, \bigcap_{T>0} \,\overline{\bigl\{
  \cS_t(u_0)\,;\, t \ge T\bigr\}}^{\,L^\infty_\loc}\,;
\end{equation}

\vspace{-10pt}
\item $\omega_0(u_0)$ is fully invariant, attractive, and chain recurrent, namely\:

\begin{itemize}[leftmargin=10pt,itemsep=3pt,topsep=2pt]
  \item $\cS_t(\omega_0(u_0)) = \omega_0(u_0)$ for all $t \ge 0$;
  \item for any neighborhood $\cN$ of $\omega_0(u_0)$, there exists $T>0$ such that
  $\cS_t(u_0)\in \cN$ for all $t \ge T$; 
  \item for each $v_0\in \omega_0(u_0)$ and any $T,\epsilon>0$, there exists a closed
    $(\epsilon,T)$-pseudo-orbit in $\omega_0(u_0)$ starting at $v_0$, that is, there exist
    finite sequences $v_j\in \omega(u_0)$ and $t_j\ge T$ for $0\le j\le N-1$, such that
    $v_N=v_0$ and $d(v_{j+1},\cS_{t_j}(v_j))<\epsilon$ for all $j \in \{0,\dots,N{-}1\}$.  
\end{itemize}
In particular, if $v_0\in \omega_0(u_0)$, there exists an entire solution $v \in C^0(\R,L^\infty(\R))$
of \eqref{SCL} such that $v(t,\cdot) \in \omega_0(u_0)$ for all $t \in \R$ and $v(0,\cdot) = v_0$;
moreover $\omega_0(v_0)\subset \omega_0(u_0)$.
\item $\omega_0(u_0)$ is a bounded subset of $C^k_b(\R)$ for all $k \in \N$, and any $v
\in \omega_0(u_0)$ satisfies $v'(x) \le 0~\forall x \in \R$. 
\end{enumerate}
\end{prop}

\begin{Proof}
Smoothing properties of the parabolic equation \eqref{SCL} and a priori bounds
for the solutions and their derivatives guarantee that, for any $u_0 \in L^\infty(\R)$
and any $k \in \N$, the solution $\cS_t(u_0)$ is uniformly bounded in $C^k_b(\R)$
for $t\ge 1$. This does not imply that the forward trajectory $\gamma_+(u_0) :=
\{\cS_t(u_0) \,;\, t\ge 0\}$ is compact in $L^\infty_\loc(\R)$, because in
general the map $t \mapsto \cS_t(u_0)$ is not continuous at $t = 0$ in that
topology. However, for any $T > 0$, the trajectory $\gamma_+(\cS_T(u_0)) =
\{\cS_t(u_0) \,;\, t\ge T\}$ is relatively compact and connected in 
$L^\infty_\loc(\R)$. Topological properties (a) and dynamic properties (b) of 
the $\omega$-limit set $\omega_0(u_0)$ follow in a standard fashion. We include 
some details here for later reference.

\smallskip\noindent 1) {\em Compactness and attractivity.} It is easy to 
verify that the relation \eqref{omega0exp} is equivalent to the definition 
\eqref{e:omega0}. Now \eqref{omega0exp} shows that $\omega_0(u_0)$ is 
the intersection of a decreasing family of non-empty compact sets, so that 
$\omega_0(u_0)$ is itself compact and non-empty. By the same argument, if
$\cN$ is any neighborhood of $\omega_0(u_0)$ in $L^\infty_\loc(\R)$, there
exists $T > 0$ such that $\gamma_+(\cS_T(u_0)) \subset \cN$, which
proves attractivity. 

\smallskip\noindent 2) {\em Connectedness.} We argue by contradiction: if
$\omega_0(u_0) = A_1 \cup A_2$ where $A_1, A_2$ are non-empty disjoint closed
sets, then $A_1, A_2$ are in fact compact and are therefore separated by a
distance $\epsilon > 0$. If $\cN_1,\cN_2$ are $\epsilon/3$-neighborhoods of
$A_1,A_2$, respectively, then $\cN_1,\cN_2$ are non-empty disjoint open sets, and
the attractivity property shows that, for $T > 0$ sufficiently large, the
connected forward orbit $\gamma_+(\cS_T(u_0))$ is contained in the neighborhood
$\cN := \cN_1\cup \cN_2$, without being included in either $\cN_1$ or $\cN_2$,
which is clearly impossible.

\smallskip\noindent 3) {\em Full invariance.} If $v_0\in\omega_0(u_0)$ there
exists a sequence $t_k \to +\infty$ such that $\cS_{t_k}(u_0) \to v_0$ in
$L^\infty_\loc$. By Lemma~\ref{contlem}, for any $t > 0$, we thus have
\[
  \cS_t(v_0) \,=\, \cS_t\bigl(\lim_{k\to\infty}\cS_{t_k}(u_0)\bigr) \,=\, 
  \lim_{k\to\infty}\cS_{t_k+t}(u_0) \,\in\, \omega_0(u_0)\,,
\]
which proves that $\cS_t(\omega_0(u_0)) \subset \omega_0(u_0)$. Similarly, we can
extract a subsequence (still denoted by $t_k$) such that $\cS_{t_k-t}(u_0) \to v_{-t}
\in\omega(u_0)$, where $v_{-t}$ satisfies $\cS_t(v_{-t}) = v_0$. Altogether,
this shows that $\cS_t(\omega_0(u_0)) = \omega_0(u_0)$ for all $t \ge 0$.
It is easy to deduce that, given any $v_0\in\omega_0(u_0)$, there exists
an entire solution $v \in C^0(\R,L^\infty_\loc(\R))$ of \eqref{SCL} such that
$v(0) = v_0$.

\smallskip\noindent 4) {\em Chain recurrence.} This is a consequence of
continuity and attractivity, which can be established as follows. 
Fix $\epsilon,T > 0$ and take $v_0\in\omega_0(u_0)$. By continuity,
there exists $\delta \in (0,\epsilon/2)$ such that, for all $u_1,u_2
\in \Sigma(u_0)$ such that $d(u_1,u_2) < \delta$, one has $d(\cS_t(u_1),
\cS_t(u_2)) < \epsilon/2$ for all $t \in [T,2T]$. By attractivity, we
can then choose $t_* > 0$ such that $\dist(\cS_t(u_0),\omega_0(u_0)) < \delta$
for all $t\ge t_*$. Now, we take $T_0 \ge t_*$ such that $d(\cS_{T_0}(u_0),v_0)
<\delta$, and also $T_* \ge t_* + T$ such that $d(\cS_{T_*}(u_0),v_0) < \delta$.
For some $N \in \N^*$, we define intermediate times $T_1,\dots,T_N$ such
that $T_N = T_*$ and 
\[
  t_j \,:=\, T_{j+1} - T_j \,\in\, [T,2T]\,, \quad
  \hbox{for all } \,j \in \{0,\dots,N-1\}\,.
\]
Finally we denote $\tilde u_j = \cS_{T_j}(u_0)$ for $j = 0,\dots,N$, and we take
$v_j \in \omega_0(u_0)$ such that $d(v_j,\tilde u_j) < \delta$. Note that
$v_0$ is given from the beginning, and we can take $v_N = v_0$. We claim that
the sequence $v_j$ for $j = 0,\dots,N$ is the desired pseudo-orbit. Indeed
for $j = 0,\dots,N-1$ we have $\tilde u_{j+1} = \cS_{t_j}(\tilde u_j)$, hence
\[
  d\bigl(v_{j+1},\cS_{t_j}(v_j)\bigr) \,\le\, d\bigl(v_{j+1},\tilde u_{j+1}\bigr)
  + d\bigl(\cS_{t_j}(\tilde u_j),\cS_{t_j}(v_j)\bigr) \,<\, \delta + \epsilon/2
  \,<\, \epsilon\,,
\]
where we used the uniform continuity of $\cS_{t_j}$ and the fact
that $t_j \in [T,2T]$. 

\smallskip\noindent 5) Assertion (c) is an easy consequence of parabolic
smoothing and Oleinik's inequality \eqref{Oleinik}.  
\end{Proof}

We next consider the larger $\omega$-limit set \eqref{e:omega}, where
limit points are considered up to translations in space. The analogue
of Proposition~\ref{omega-prop0} is: 

\begin{prop}\label{omega-prop}
For any $u_0 \in L^\infty(\R)$, the $\omega$-limit set $\omega(u_0)$
defined by \eqref{e:omega} is bounded in $L^\infty(\R)$ and, when equipped
with the topology of $L^\infty_\loc(\R)$, has the following properties\:

\begin{enumerate}[leftmargin=20pt,itemsep=2pt,topsep=0pt,label=\alph*)]

\item $\omega(u_0)$ is non-empty, compact, connected, and
\begin{equation}\label{omegaexp}
  \omega(u_0) \,=\, \bigcap_{T>0} \,\overline{\bigl\{
  \cT_y \cS_t(u_0)\,;\, t \ge T,\,y \in \R\bigr\}}^{\,L^\infty_\loc}\,;
\end{equation}

\vspace{-10pt}
\item $\omega(u_0)$ is fully invariant in time, translation invariant
in space, uniformly attractive, and chain recurrent up to translations;

\item $\omega(u_0)$ is a bounded subset of $C^k_b(\R)$ for all $k \in \N$,
and any $v \in \omega(u_0)$ satisfies $v'(x) \le 0~\forall x \in \R$. 
\end{enumerate}
\end{prop}

\begin{Proof}
The proof is completely parallel to that of Proposition~\ref{omega-prop0},
and we just indicate here the main differences. The starting point is the formula
\eqref{omegaexp}, which is easily derived from the definition \eqref{e:omega}.
Since the space-time trajectory $\bigl\{\cT_y \cS_t(u_0)\,;\, t \ge T,\,y \in
\R\bigr\}$ is relatively compact in $L^\infty_\loc(\R)$ for any $T > 0$,
we see that $\omega(u_0)$ is non-empty and compact as the decreasing
intersection of non-empty compact sets. Moreover, if $\cN$ is any neighborhood
of $\omega(u_0)$ in $L^\infty_\loc(\R)$, we have $\bigl\{\cT_y \cS_t(u_0)\,;\,
t \ge T,\,y \in \R\bigr\} \subset \cN$ for any sufficiently large $T > 0$, which
means that $\omega(u_0)$ attracts the trajectory $\cT_y\cS_t(u_0)$ uniformly
in $y \in \R$ as $t \to +\infty$. As a consequence, since the space-time trajectory 
is connected for all $T > 0$, the same argument as in Proposition~\ref{omega-prop0} 
shows that $\omega(u_0)$ is a connected set. There is no difference either 
in the reasoning showing that $\cS_t(\omega(u_0)) = \omega(u_0)$ for all 
$t \ge 0$. Finally, the definition \eqref{e:omega} immediately implies 
that $\cT_y(\omega(u_0)) = \omega(u_0)$ for all $y \in \R$, and the boundedness
properties (c) are established exactly as before. 

The main difference we would like to point out is that $\omega(u_0)$ is 
not chain recurrent in the sense of Proposition~\ref{omega-prop0}, 
but only in a weaker sense that can be called ``chain recurrence up 
to translations''. The precise definition is as follows: for each $v_0\in
\omega(u_0)$ and any $T,\epsilon>0$, there exist finite sequences
$v_j\in \omega(u_0)$, $t_j\ge T$, and $y_j \in \R$ for $0\le j\le N-1$,
such that $v_N=v_0$ and $d(v_{j+1},\cT_{y_j}\cS_{t_j}(v_j))<\epsilon$
for all $j \in \{0,\dots,N{-}1\}$. In other words, the definition of
the $(\epsilon,T)$-pseudo-orbit involves spatial shifts $y_j$ in
addition to the time shifts $t_j$, which is natural in view of
\eqref{e:omega}. The existence of such a pseudo-orbit for all
$v_0 \in \omega(u_0)$, all $\epsilon > 0$, and all $T > 0$ is
established by the same argument as in Proposition~\ref{omega-prop0}. 
\end{Proof}

Although the $\omega$-limit sets \eqref{e:omega0}, \eqref{e:omega} are
relatively easy to define and enjoy the nice properties listed in
Propositions~\ref{omega-prop0} and \ref{omega-prop}, it is notoriously difficult
to compute them for arbitrary initial data. In the case of equation~\eqref{SCL},
general results in this direction are only available under monotonicity
assumptions. If $u_0$ is increasing, the solution $\cS_t(u_0)$ remains
increasing for all $t > 0$ by the maximum principle, and
$\|\partial_x \cS_t(u_0)\|_{L^\infty} \to 0$ as $t \to +\infty$ by Oleinik's
inequality \eqref{Oleinik}. This implies that $\omega(u_0)$ consists of the
constant states $u \equiv \gamma$ for all $\gamma \in [\alpha,\beta]$, where
$\alpha,\beta$ are as in \eqref{alphabet}. On the other hand, if $u_0$ is
decreasing, Proposition~\ref{mono-prop} below implies that $\omega(u_0)$ is the
set of all translates of the viscous shock $\phi_{\beta,\alpha}$, supplemented
with the constant states $u \equiv \alpha$ and $u \equiv \beta$. A similar
conclusion is reached if $u_0$ satisfies the assumptions of Il'in and Oleinik's
result \cite{IO}. Incidentally, we observe in this example that $\omega(u_0)$ is
a heteroclinic orbit in the terminology of dynamical systems, so that $\omega(u_0)$
is not chain recurrent. 

More generally, if $u_0 \in L^\infty(\R)$ and if there is a
$v \in \omega(u_0)$ that is not a constant, then $v$ is decreasing by
Proposition~\ref{omega-prop}, and since $\omega(v) \subset \omega(u_0)$ we
deduce that $\omega(u_0)$ contains the translates of a viscous shock, as
asserted in \eqref{omvshock}. So we see that Proposition~\ref{p:shock} is a
direct consequence of Propositions~\ref{omega-prop} and \ref{mono-prop}.  In
addition we have

\begin{cor}\label{omega-cor}
For any $u_0 \in L^\infty(\R)$, the $\omega$-limit set $\omega(u_0)$
contains a constant state or a viscous shock. 
\end{cor}

This statement can be compared with a result by S. Slijep\v{c}evi\'c and the
first author \cite{GS} which shows that, for a general class of dissipative
systems including reaction-diffusion equations on the real line, the
$\omega$-limit set of a bounded trajectory always contains an equilibrium. For the
viscous conservation law \eqref{SCL}, where all Galilean frames are equivalent, 
the role of equilibria is played by the constant states and the viscous shocks. 

From a different perspective, one may wonder which collections of entire
solutions $v\in C^0(\R,L^\infty_\loc(\R))$ of \eqref{SCL} may occur as $\omega$-limit
sets of bounded initial data. From the results presented thus far, only
non-empty, compact, connected, invariant, and chain-recurrent sets are
candidates. In general, compact, connected sets with a chain-recurrent flow are
precisely the possible $\omega$-limit sets of flows, in the sense that any such
set is topologically conjugated to an $\omega$-limit set of some flow
\cite{FrSe}.  It is however not clear at all if any compact, connected,
invariant, and chain-recurrent set within the family of entire solutions is
realized as the $\omega$-limit set of the particular flow generated by the
conservation law \eqref{SCL}. 

\begin{rem}\label{other-def}
In the definition \eqref{e:omega} of the $\omega$-limit set $\omega(u_0)$,
we allow for arbitrary spatial shifts $x_k \in \R$ while the temporal shifts
$t_k$ must go to infinity. This is clearly not the only possibility. One the
one hand, we could restrict the class of spatial shifts by imposing,
for instance, a Galilean constraint of the form $|x_k| \le ct_k$ with,
typically, $c > \|u_0\|_\infty$. Actually, we could even require that $x_k/t_k$
converges to some limit in $[-c,c]$ as $k \to +\infty$. In a different
direction, we could consider all spatio-temporal shifts such that $|x_k| + t_k
\to \infty$ as $k \to \infty$, which potentially gives an $\omega$-limit set 
even larger than \eqref{e:omega}. However, in the examples we are
aware of, these alternative possibilities do not seem to change the nature
of the $\omega$-limit set in a profound way, so in what follows we stick to
the original definition \eqref{e:omega}. 
\end{rem}

\section{Convergence to shocks for monotone initial data}\label{s:4}

The main result of this section is: 

\begin{prop}\label{mono-prop}
Assume that $u_0 \in L^\infty(\R)$ is nonincreasing and satisfies
$\alpha < \beta$, where
\begin{equation}\label{alphabet2}
  \alpha \,:=\, \lim_{x \to +\infty}u_0(x)\,, \qquad
  \beta \,:=\, \lim_{x \to -\infty}u_0(x)\,.
\end{equation}
Then there exists a smooth function $s : (0,+\infty) \to \R$
such that the solution $u$ of \eqref{SCL} with initial
data $u_0$ satisfies
\begin{equation}\label{mono-conv}
  \sup_{x \in \R} \,\bigl|u(t,x) - \phi_{\beta,\alpha}(x-s(t))\bigr|
  \,\xrightarrow[t \to +\infty]{}\, 0\,.
\end{equation}
Moreover $s(t)/t \to c$ as $t \to +\infty$, where $c$ is
given by the Rankine-Hugoniot formula \eqref{RH}. 
\end{prop}

\begin{rem}\label{mono-rem}
It is not difficult to find examples for which the shift function
$s$ in \eqref{mono-conv} is not asymptotically linear, namely
$s(t) - ct$ has no limit as $t \to +\infty$. For instance, assume
that $\beta - u_0\in L^1(\R_-)$ but $u_0 - \alpha \notin L^1(\R_+)$.
Given any $\gamma > 0$, we define $\tilde u_0(x) = u_0(x)$ for $x \le
\gamma$ and $\tilde u_0(x) = \alpha$ for $x > \gamma$. As 
$\tilde u_0 \le u_0$, monotonicity implies that $\cS_t \tilde u_0
\le \cS_t u_0$ for all $t > 0$. On the other hand, since $\tilde
u_0 - \phi_{\beta,\alpha} \in L^1(\R)$, we can apply the result of
\cite{IO} to deduce that $(\cS_t \tilde u_0)(x)$ converges uniformly
to $\phi_{\beta,\alpha}(x-ct-\tilde x_0)$ as $t \to +\infty$, where
\[
  \tilde x_0 \,:=\, \frac{1}{\beta-\alpha}\left(\int_{-\infty}^\gamma
  \bigl(u_0(x) - \phi_{\beta,\alpha}(x)\bigr)\dd x + \int_\gamma^{+\infty}
  \bigl(\alpha - \phi_{\beta,\alpha}(x)\bigr)\dd x\right)\,.
\]
In particular one has $\liminf_{t \to +\infty} (s(t) - ct) \ge \tilde x_0$
by monotonicity. Now, taking taking $\gamma \to +\infty$, we see that
$\tilde x_0 \to +\infty$ by assumption on $u_0$, and we conclude that
$s(t) - ct \to +\infty$ as $t \to +\infty$. A more explicit example of
such a ``sublinear shift'' will be given in Section~\ref{s:52} below. 
\end{rem}

The remainder of this section is devoted to the proof of Proposition~\ref{mono-prop}. 
Assume that the initial data $u_0 \in L^\infty(\R)$ are nonincreasing
and satisfy \eqref{alphabet2} for some $\alpha < \beta$. The solution
$u(t) = \cS_t u_0$ of \eqref{SCL} is smooth for positive times,
and the strong maximum principle implies that $\partial_x u(t,x) < 0$
for all $t > 0$ and all $x \in \R$. On the other hand, using for instance
Lemma~\ref{contlem}, it is not difficult to verify that the limits of
$u(t,x)$ as $x \to \pm\infty$ are independent of time. As a consequence,
for each $t > 0$, there exists a unique point $s(t) \in \R$ such that
\begin{equation}\label{shift-def}
  u(t,s(t)) \,=\, \frac{\alpha+\beta}{2}\,.
\end{equation}
Moreover $s(t)$ is a smooth function of time thanks to the implicit
function theorem.

\begin{lem}\label{shift-lem}
The shift function $s : (0,+\infty) \to \R$ defined by \eqref{shift-def}
satisfies
\begin{equation}\label{shift-asym}
  \lim_{t \to +\infty} \frac{s(t)}{t} \,=\, c \,:=\, \frac{f(\beta)-
  f(\alpha)}{\beta-\alpha}\,.
\end{equation}
\end{lem}

\begin{Proof}
We use the monotonicity of the evolution map $\cS_t$ to compare
the solution $u(t) = \cS_t u_0$ with suitably translated viscous shocks.
Take $\epsilon > 0$ small enough so that
\begin{equation}\label{eps-cond}
  0 \,<\, \epsilon \,<\, \frac{\beta-\alpha}{2}\,, \qquad \hbox{hence}
  \quad \alpha + \epsilon \,<\, \frac{\alpha+\beta}{2} \,<\, \beta - \epsilon\,.
\end{equation}
Since $u_0$ is nonincreasing and satisfies \eqref{alphabet2}, there exist
$x_+(\epsilon) \in \R$ and $x_-(\epsilon) \in \R$ such that
\begin{equation}\label{u0-comp}
  \phi_{\beta-\epsilon,\alpha-\epsilon}(x - x_-(\epsilon)) \,\le\, u_0(x)
  \,\le\, \phi_{\beta+\epsilon,\alpha+\epsilon}(x - x_+(\epsilon))\,,
  \qquad \forall\,x \in \R\,,
\end{equation}
where $\phi_{\beta\pm\epsilon,\alpha\pm\epsilon}$ denotes the viscous shock
connecting $\beta\pm\epsilon$ with $\alpha\pm\epsilon$. In fact, it
is straightforward to verify that \eqref{u0-comp} holds as soon as
$x_+(\epsilon) \gg 1$ and $x_-(\epsilon) \ll -1$ are sufficiently large,
depending on $\epsilon$. By monotonicity, we deduce from \eqref{u0-comp} that
\begin{equation}\label{u-comp}
  \phi_{\beta-\epsilon,\alpha-\epsilon}\bigl(x - c_-(\epsilon)t - x_-(\epsilon)\bigr)
  \,\le\, u(t,x) \,\le\, \phi_{\beta+\epsilon,\alpha+\epsilon}\bigl(x - c_+(\epsilon)t
  - x_+(\epsilon)\bigr)\,,  
\end{equation}
for all $t \ge 0$ and all $x \in \R$, where the speeds $c_\pm(\epsilon)$ are
given by the Rankine-Hugoniot formulas
\begin{equation}\label{cpm-def}
  c_+(\epsilon) \,:=\, \frac{f(\beta{+}\epsilon) - f(\alpha{+}\epsilon)}{\beta-\alpha}\,,
  \qquad
  c_-(\epsilon) \,:=\, \frac{f(\beta{-}\epsilon) - f(\alpha{-}\epsilon)}{\beta-\alpha}\,.
\end{equation}
On the other hand, due to the second inequality in \eqref{eps-cond}, there
exist $s_+(\epsilon) \in \R$ and $s_-(\epsilon) \in \R$ such that
\begin{equation}\label{shiftpm}
  \phi_{\beta+\epsilon,\alpha+\epsilon}\bigl(s_+(\epsilon)\bigr) \,=\,
  \phi_{\beta-\epsilon,\alpha-\epsilon}\bigl(s_-(\epsilon)\bigr) \,=\,
  \frac{\alpha+\beta}{2}\,.
\end{equation}
In view of \eqref{shiftpm}, we deduce from \eqref{u-comp} that the
shift function defined by \eqref{shift-def} satisfies
\begin{equation}\label{shift-comp}
  c_-(\epsilon)t + x_-(\epsilon) + s_-(\epsilon) \,\le\, s(t) \,\le\,
  c_+(\epsilon)t + x_+(\epsilon) + s_+(\epsilon)\,, \qquad
  \forall\,t > 0\,.
\end{equation}
In particular we infer from \eqref{shift-comp} that
\begin{equation}\label{speed-comp}
  c_-(\epsilon) \,\le\, \liminf_{t \to +\infty} \frac{s(t)}{t} \,\le\,
  \limsup_{t \to +\infty} \frac{s(t)}{t} \,\le\, c_+(\epsilon)\,. 
\end{equation}
Finally it is clear from \eqref{cpm-def} that $c_\pm(\epsilon) \to c$
as $\epsilon \to 0$, which concludes the proof of \eqref{shift-asym}.
\end{Proof}

In a second step, we consider the auxiliary function $v$ defined by
\begin{equation}\label{v-def}
  v(t,x) \,=\, \partial_x u(t,x) - f(u(t,x)) + c u(t,x) + d\,,
  \qquad t > 0\,,\quad x \in \R\,,
\end{equation}  
where the constants $c,d$ are defined in \eqref{RH}. This function
is smooth for positive times and a direct calculation shows that it
satisfies the evolution equation
\begin{equation}\label{v-eq}
  \partial_t v(t,x) + f'\bigl(u(t,x)\bigr)\partial_x v(t,x) \,=\,
  \partial_x^2 v(t,x)\,,   \qquad t > 0\,,\quad x \in \R\,.
\end{equation}
The key step in the proof of Proposition~\ref{mono-prop} is\:

\begin{lem}\label{v-lem}
The function $v$ defined by \eqref{v-def} converges uniformly to zero
as $t \to +\infty$\:
\begin{equation}\label{v-conv}
  \sup_{x \in \R}\,|v(t,x)|   \,\xrightarrow[t \to +\infty]{}\, 0\,.
\end{equation}
\end{lem}

\begin{Proof}
Shifting the initial time if needed, we can assume without loss of
generality that the functions  $u(t,x)$ and $v(t,x)$ are smooth for all
$t \ge 0$. We consider the linear advection-diffusion equation
\begin{equation}\label{w-eq}
  \partial_t w(t,x) + f'\bigl(u(t,x)\bigr)\partial_x w(t,x) \,=\,
  \partial_x^2 w(t,x)\,,   \qquad t > 0\,,\quad x \in \R\,,
\end{equation}
where the function $u(t,x)$ is considered as given. The following
$L^p$--$L^q$ estimates are known for the solution of \eqref{w-eq} with initial
data $w_0$\:
\begin{equation}\label{w-est}
  \sup_{t \ge 0} \|w(t,\cdot)\|_{L^\infty(\R)} \,\le\, \|w_0\|_{L^\infty(\R)}\,,
  \qquad
  \sup_{t > 0} t^{1/2}\|w(t,\cdot)\|_{L^\infty(\R)} \,\le\, C\|w_0\|_{L^1(\R)}\,,
\end{equation}
where $C > 0$ is a universal constant. While the first bound in \eqref{w-est}
is a direct consequence of the parabolic maximum principle, the second one
takes into account the convexity of the flux function $f$ as well as
the monotonicity of the solution $u$. For the reader's convenience, we
provide a proof of the second estimate \eqref{w-est} in Section~\ref{ssecA3}. 

For the time being, to conclude the proof of Lemma~\ref{v-lem}, we consider the
solution $v$ of \eqref{v-eq} with initial data $v_0 := \partial_x u_0 - f(u_0)
+ cu_0 + d$. We know that $\partial_x u_0 \in L^1(\R)$ because $u_0$ is
decreasing and bounded, and that $f(u_0) - cu_0 - d$ converges to zero as
$x \to \pm\infty$ because of \eqref{RH}. As a consequence, given any $\epsilon > 0$,
we can decompose $v_0 = w_1 + w_2$, where $w_1 \in L^1(\R)$ and
$\|w_2\|_{L^\infty} \le \epsilon$. For $j = 1,2$ we denote by $w_j(t)$ the
solution of \eqref{w-eq} with initial data $w_j$, so that $v(t) = w_1(t) + w_2(t)$
by linearity. Using \eqref{w-est}, we infer that
\[
  \|v(t)\|_{L^\infty} \,\le\, \|w_1(t)\|_{L^\infty} + \|w_2(t)\|_{L^\infty}
  \,\le\, C t^{-1/2} \|w_1\|_{L^1} + \|w_2\|_{L^\infty}\,, \qquad
  \forall\,t > 0\,,
\]
so that $\limsup_{t \to +\infty} \|v(t)\|_{L^\infty} \le \epsilon$. Since $\epsilon > 0$
was arbitrary, this gives \eqref{v-conv}.   
\end{Proof}

Equipped with Lemmas~\ref{shift-lem} and \ref{v-lem}, it is now straightforward
to conclude the proof of Proposition~\ref{mono-prop}. Let $g : \R \to \R$ be
the convex function defined by $g(u) = f(u) - cu - d$, where $c,d$ are given
by \eqref{RH}, and let
\[
  L \,:=\, \max\bigl\{|g'(u)|\,;\, u \in [\alpha,\beta]\bigr\} \,=\,
  \max\bigl\{c - f'(\alpha)\,,\,f'(\beta)-c\bigr\}\,.
\]
Fix any $t > 0$. In view of \eqref{v-def}, the function $\psi : \R \to \R$ defined
by $\psi(y) = u\bigl(t,y+s(t)\bigr)$ satisfies the ODE
\[
  \psi'(y) \,=\, g\bigl(\psi(y)\bigr) + v\bigl(t,y+s(t)\bigr)\,, \quad
  \forall\, y\in \R\,, \qquad \hbox{with}~\,\psi(0) \,=\, \frac{\alpha+\beta}{2}\,.
\]
This is to be compared with the ODE \eqref{phiODE} satisfied by the viscous
shock $\phi(y) := \phi_{\beta,\alpha}(y)$, namely $\phi'(y) = g(\phi(y))$ and
$\phi(0) = (\alpha+\beta)/2$. If $w = \psi - \phi$, we infer that $|w'(y)| \le
L |w(y)| + \|v(t)\|_{L^\infty}$. Integrating this differential inequality and
recalling that $w(0) = 0$, we obtain
\begin{equation}\label{uloc-conv}
  |w(y)| \,\equiv \,\bigl|u(t,y+s(t)) - \phi_{\beta,\alpha}(y)\bigr| \,\le\,
  |y|\,e^{L|y|}\,\|v(t)\|_{L^\infty}\,, \qquad \forall\,y \in \R\,.
\end{equation}
Since $\|v(t)\|_{L^\infty} \to 0$ by Lemma~\ref{v-conv}, it follows from
\eqref{uloc-conv} that $u(t,y +s(t)) - \phi_{\beta,\alpha}(y)$ converges to zero
as $t \to +\infty$, uniformly for $y$ in any compact interval. Taking into
account the fact that both functions $u(t,\cdot)$ and $\phi_{\beta,\alpha}$
are decreasing and have the same limits as $y \to \pm\infty$, we deduce that
the convergence is in fact uniform for all $y \in \R$. This proves \eqref{mono-conv},
and we already established in Lemma~\ref{shift-lem} that $s(t)/t$ has a limit as
$t \to +\infty$. \QED

\begin{rem}\label{nonstrict-rem}
Neither the strict convexity of the flux nor the Lax condition for the
shock is used in the proof of Proposition~\ref{mono-prop}, which therefore
remains valid if we only assume that $f''(u) \ge 0$ for all $u \in [\alpha,\beta]$. 
\end{rem}

\section{Representation of entire solutions via probability measures}\label{s:5}

From now on we restrict our attention to the special case of Burgers' equation
\begin{equation}\label{Burgers}
  \partial_t u(t,x) + u(t,x)\partial_x u(t,x) \,=\, \partial_x^2 u(t,x)\,,
  \qquad t > 0\,, \quad x \in \R\,,
\end{equation}
which corresponds to taking $f(u) = u^2/2$ in \eqref{SCL}. The Cauchy problem for
\eqref{Burgers} can be solved in explicit form through the celebrated Cole-Hopf
transformation \cite{Co,Ho}. As is easily verified, if $U(t,x)$ is any positive
solution of the heat equation $\partial_t U = \partial_x^2 U$, a corresponding
solution $u(t,x)$ of \eqref{Burgers} is obtained by setting
\begin{equation}\label{CH}
  u(t,x) \,=\, -2 \,\frac{\partial_x U(t,x)}{U(t,x)}\,,
  \qquad t > 0\,, \quad x \in \R\,. 
\end{equation}
It is tempting to conclude that the dynamics of \eqref{Burgers} is trivial,
but one should keep in mind that bounded solutions of \eqref{Burgers} are
associated via \eqref{CH} to functions $U(t,x)$ that may grow exponentially
as $|x| \to \infty$, and this seriously complicates the process of computing
the long-time asymptotics, even if $U$ solves a simple equation. 

A beautiful application of the Cole-Hopf transformation is the derivation of
the representation formula \eqref{UPK-rep0} for bounded entire solutions of
\eqref{Burgers}. By Oleinik's inequality \eqref{Oleinik}, any entire solution
$u$ of \eqref{Burgers} necessarily satisfies $\partial_x u(t,x) \le 0$ for all
$(t,x) \in \R \times \R$. If we assume in addition that $u$ is bounded, the
limits $\alpha,\beta \in \R$ defined by
\begin{equation}\label{alphabet3}
  \alpha \,=\, \lim_{x \to +\infty}u(t,x)\,, \qquad
  \beta \,=\, \lim_{x \to -\infty}u(t,x)\,,
\end{equation}
are therefore finite, and independent of time (the last property follows
from Lemma~\ref{contlem} and translation invariance.) We then have the
following result: 

\begin{prop}\label{UPK-prop} \null\hspace{-2pt}{\bf\cite{UPK}}
If $u : \R \times \R \to \R$ is a bounded entire solution of Burgers'
equation \eqref{Burgers}, there exists a unique probability measure
$\mu$ supported on $[\alpha,\beta]$, where $\alpha,\beta$ are the limits
defined in \eqref{alphabet3}, such that
\begin{equation}\label{UPK-rep}
  u(t,x) \,=\, \genfrac{}{}{1pt}{0}{\int z\,e^{-zx/2 + z^2t/4} \dd\mu(z)}{
  \int e^{-zx/2 + z^2t/4} \dd\mu(z)}\,, \qquad t \in \R\,, \quad x \in \R\,.
\end{equation}
\end{prop}

The proof of Proposition~\ref{UPK-prop} uses a representation of positive, ancient
solutions of the heat equation of the form $U(t,x) = \int e^{-zx/2 + z^2t/4}\dd\mu(z)$,
which substituted into \eqref{CH} immediately gives \eqref{UPK-rep}. The formula
\eqref{UPK-rep} appears implicitly in the work of Kenig and Merle \cite{KM}, and
explicitly in the PhD thesis of U.~P.~Karunathilake \cite{UPK}. Since the
latter reference is not widely available, we reproduce the proof here
and establish the uniqueness of the measure $\mu$, which is not asserted
in \cite{UPK}. 

\begin{rem}\label{ancientrem}
When restricted to the time interval $(-\infty,T)$, for some $T \in \R$, 
the representation formula \eqref{UPK-rep} remains valid for ancient solutions
$u : (-\infty,T)\times\R \to \R$ that are not necessarily bounded. In that
situation $\mu$ is just a positive measure, supported on the closure of the
range of $u$, which may have finite or infinite mass (in the latter case it
cannot be normalized, and no uniqueness is claimed). In what follows we focus on
bounded entire solutions, due to their connection with $\omega$-limit sets, but we
allow ourselves occasional comments on the general case. 
\end{rem}

We now turn to the proof of Proposition \ref{UPK-prop}, then study some examples
of measures in Section \ref{s:52} and conclude with an analysis of
\eqref{UPK-rep} for $t\to-\infty$ in Section \ref{s:53}.

\subsection{Representation of ancient solutions}\label{s:51}

Assume that $u : \Omega_- \to \R$ is a smooth solution of Burgers' equation
\eqref{Burgers} on the space-time domain $\Omega_- := \bigl\{(t,x) \in \R^2 \,;\, t < 0
\bigr\}$. Our goal is to obtain a representation formula for $u$ in terms
of a positive measure on the real line. We proceed in three steps. 

\medskip\noindent{\bf Step 1\:} {\em Cole-Hopf transformation} \cite{Co,Ho}.
We first define
\begin{equation}\label{U-CH}
  U(t,x) \,=\, \exp\Bigl(-\frac12\int_0^x u(t,y)\dd y + a(t)\Bigr)\,,
  \qquad \forall\, (t,x) \in \Omega_-\,,
\end{equation}
where
\[
  a(t) \,=\, \int_{t_0}^t \Bigl(\frac14\,u(s,0)^2 - \frac12 \,\partial_x
  u(s,0)\Bigr)\dd s\,, \qquad \forall\, t < 0\,.
\]
Here $t_0 < 0$ is some arbitrary reference time. A direct calculation shows that
\[
  \partial_t U(t,x) \,=\, U(t,x)\Bigl(-\frac12\,\partial_x u(t,x) + \frac14\,
  u(t,x)^2\Bigr) \,=\, \partial_x^2 U(t,x)\,, \qquad \forall\, (t,x) \in \Omega_-\,.
\]
Our solution $u$ of Burgers' equation can therefore be expressed as
\begin{equation}\label{u-CH}
  u(t,x) \,=\, \frac{-2\partial_x U(t,x)}{U(t,x)}\,,
  \qquad \forall\, (t,x) \in \Omega_-\,,
\end{equation}
where $U :\Omega _- \to (0,+\infty)$ is a positive solution of the heat
equation $\partial_t U = \partial_x^2 U$. 

\medskip\noindent{\bf Step 2\:} {\em Appell transformation} \cite{Ap,Wi2}.
We next transform the ancient solution $U(t,x)$ of the heat equation into
a solution $V(t,x)$ of the same equation which is defined for positive
times, namely on the space-time domain $\Omega_+ := \bigl\{(t,x) \in \R^2
\,;\, t > 0\bigr\}$. This remarkable transformation, first discovered by
P.~Appell, takes the form
\begin{equation}\label{V-Ap}
  V(t,x) \,=\, K(t,x)\,U\Bigl(\frac{-1}{t}\,,\,\frac{-x}{t}\Bigr)\,,
  \qquad \forall\, (t,x) \in \Omega_+\,,
\end{equation}
where $K(t,x)$ is the fundamental solution of the one-dimensional heat
equation\:
\begin{equation}\label{K-def}
  K(t,x) \,=\, \frac{e^{-x^2/(4t)}}{\sqrt{4\pi t}}\,,
  \qquad \forall\, (t,x) \in \Omega_+\,.
\end{equation}  
A simple calculation shows that $\partial_t V(t,x) = \partial_x^2 V(t,x)$ for
all $(t,x) \in \Omega_+$, and by construction $V(t,x)$ is strictly positive
on the domain $\Omega_+$. 

\medskip\noindent{\bf Step 3\:} {\em Poisson representation} \cite{Wi1}.
A classical result due to Widder \cite[Theorem~6]{Wi1} asserts that,
if $V(t,x)$ is a {\em nonnegative} solution of the heat equation in
$\Omega_+$, there exists a (unique) positive Borel measure $\mu$ on
$\R$ such that
\begin{equation}\label{Poisson}
  V(t,x) \,=\, \int_\R K(t,x-z)\dd\mu(z)\,,   \qquad \forall\,
  (t,x) \in \Omega_+\,.
\end{equation}
It should be emphasized at this point that the convergence of the
integral in \eqref{Poisson} is part of the conclusion of Widder's
theorem. In particular, the measure $\mu(I)$ of any compact interval
$I \subset \R$ should be finite, which implies that $\mu$ is a regular
measure \cite[Theorem~2.18]{Ru}. In addition $\mu$ should have a
``moderate growth'' at infinity so that the integral in \eqref{Poisson}
is finite even when $t > 0$ is large. For instance, if
$\D\mu = e^{cz^2/4}\dd z$ for some $c > 0$, the right-hand side of
\eqref{Poisson} is infinite when $t \ge 1/c$, which contradicts the
assumption that $V$ is defined on the whole domain $\Omega_+$.

\begin{rem}\label{positive-rem}
The assumption that $V$ is nonnegative is crucial in Widder's theorem.
For instance the function $V(t,x) = (x/t)K(t,x)$ is a (sign-changing)
solution of the heat equation in $\Omega_+$ which converges to zero as
$t\to 0+$ for any $x \in \R$. As is easily verified, in that case one
cannot find any measure $\mu$ on $\R$ such that \eqref{Poisson} holds. 
\end{rem}

We now return to the ancient solution of Burgers' equation.
Combining \eqref{V-Ap} and \eqref{Poisson} we first obtain
\[
  U(t,x) \,=\, \frac{V(-1/t,x/t)}{K(-1/t,x/t)} \,=\,
  \frac{1}{K(-1/t,x/t)} \int_\R K\Bigl(\frac{-1}{t}\,,\,\frac{x}{t} -
  z\Bigr)\dd\mu(z)\,,
\]
for all $(t,x) \in \Omega_-$. The right-hand side can be simplified using
the explicit expression \eqref{K-def} of the heat kernel, leading to the
following representation formula for ancient positive solutions of the
heat equation\:
\begin{equation}\label{U-rep}
  U(t,x) \,=\, \int_\R e^{-zx/2 + z^2t/4}\dd\mu(z)\,,   \qquad
 \forall\, (t,x) \in \Omega_-\,,
\end{equation}
see also \cite[Theorem~8.1]{Wi2}. Finally, we deduce from \eqref{u-CH}
the desired representation of ancient solutions to Burgers' equation\:
\begin{equation}\label{UPK-rep3}
  u(t,x) \,=\, \genfrac{}{}{1pt}{0}{\int z\,e^{-zx/2 + z^2t/4} \dd\mu(z)}{
  \int e^{-zx/2 + z^2t/4} \dd\mu(z)}\,,  \qquad \forall\, (t,x) \in \Omega_-\,.
\end{equation}

Conversely, if $\mu$ is a positive measure on $\R$ that is finite on compact
intervals and has moderate growth at infinity, it is straightforward to verify
that the function $u$ defined by \eqref{UPK-rep3} is an ancient solution
of Burgers' equation \eqref{Burgers}. For any $(t,x) \in \Omega_-$, the quantity
$u(t,x)$ can be interpreted as the average of a random variable $z \in \R$ with
respect to the (non-normalized) measure $e^{-zx/2 + z^2t/4} \dd\mu(z)$. Introducing
the obvious notation $u(t,x) = \langle z\rangle_{t,x}$, we find by direct calculation
\[
  \partial_x u(t,x) \,=\, -\frac12 \Bigl(\langle z^2\rangle_{t,x} -
  \langle z\rangle_{t,x}^2\Bigr) \,=\, -\frac12 \Big\langle
  \bigl(z - \langle z\rangle_{t,x}\bigr)^2\Big\rangle_{t,x}
  \,\le\, 0\,.
\]
This shows that all solutions of the form \eqref{UPK-rep3} are non-increasing
in $x$, which is also a direct consequence of Oleinik's inequality \eqref{Oleinik}. 
Actually, we have $\partial_x u(t,x) < 0$ for all $x \in \R$ unless the measure
$\mu = \delta_\alpha$ is a single Dirac mass, in which case $u(t,x) = \alpha$ for
all $(t,x) \in \Omega_-$. 

\begin{ex}\label{Lebex}
The simple example where $\mu$ is just the Lebesgue measure on $\R$ is
already quite instructive. In that case, it is clear from \eqref{Poisson} that
$V \equiv 1$ on $\Omega_+$, and we deduce from \eqref{V-Ap}, \eqref{u-CH} that
\[
  U(t,x) \,=\, \frac{1}{K(-1/t,x/t)} \,=\, \frac{\sqrt{4\pi}}{\sqrt{-t}}
  \,e^{-x^2/(4t)}\,, \qquad u(t,x) \,=\, \frac{x}{t}\,, \qquad \forall\,
  (t,x) \in \Omega_-\,.
\]
We observe that the ancient solution $U(t,x)$ of the heat equation blows
up as $t \to 0-$ at any point $x \in \R$, and that the ancient solution $u(t,x)$
of Burgers' equation does so for any $x \neq 0$. 
\end{ex}

The blow-up phenomenon observed in Example~\ref{Lebex} only occurs for
unbounded solutions. Indeed, by the maximum principle, bounded ancient solutions
of either the heat equation or the Burgers equation remain uniformly bounded
at later times, and can therefore be extended to (bounded) entire solutions.
In what follows, we concentrate on bounded ancient solutions, which are
candidates for trajectories in $\omega$-limit sets of bounded initial
data. We have the following characterization: 

\begin{prop}\label{UPK-prop2}
The ancient solution $u : \Omega_- \to \R$ of Burgers' equation given by
\eqref{UPK-rep3} is bounded if and only if the measure $\mu$ has bounded
support. In that case $u$ satisfies \eqref{alphabet3} for all $t < 0$, with
\begin{equation}\label{alphabet4} 
  \alpha \,=\, \inf\,\bigl(\supp(\mu)\bigr)\,, \qquad
  \beta \,=\, \sup\,\bigl(\supp(\mu)\bigr)\,.
\end{equation}
\end{prop}

\begin{Proof}
Let $\mu$ be a (nontrivial) positive measure on $\R$ that is finite on
compact intervals and has moderate growth at infinity. We define
$\alpha,\beta$ by \eqref{alphabet4}, so that $\alpha \in [-\infty,+\infty)$
and $\beta \in (-\infty,+\infty]$. We shall show that, for any fixed $t < 0$,
the quantity $u(t,x)$ defined by \eqref{UPK-rep3} converges to $\beta$ as
$x \to -\infty$, and to $\alpha$ as $x \to +\infty$. We concentrate on the limit
at $-\infty$, the other case being similar. Setting $x = -2y$, where $y > 0$, 
we have the representation
\[
  u(t,-2y) \,=\, \genfrac{}{}{1pt}{0}{\int z\,e^{zy} \dd\nu_t(z)}{
  \int e^{zy} \dd\nu_t(z)}\,,  \qquad \hbox{where}
  \quad \D\nu_t(z) \,=\, e^{z^2 t/4}\dd\mu(z)\,.
\]
Given real numbers $a,b$ such that $a < b < \beta$, we decompose
\[
  \int_\R e^{zy} \dd\nu_t(z) \,=\, \int_{\{z < a\}} e^{zy} \dd\nu_t(z)
  \,+\, \int_{\{z \ge a\}} e^{zy} \dd\nu_t(z) \,=:\, I_0(y) + J_0(y)\,,
\]
and we observe that
\[
  I_0(y) \,\le\, e^{ay}\int_{\{z < a\}}\dd\nu_t(z)\,, \qquad
  J_0(y) \,\ge\,  \int_{\{z \ge b\}} e^{zy} \dd\nu_t(z) \,\ge\,
  e^{by}\int_{\{z \ge b\}}\dd\nu_t(z)\,.
\]
Since $b < \beta$, the last integral is strictly positive, and we deduce
that $I_0(y)/J_0(y)\to 0$ as $y \to +\infty$. A similar argument gives
\[
  \int_\R z\,e^{zy} \dd\nu_t(z) \,=\, \int_{\{z < a\}} z\,e^{zy} \dd\nu_t(z)
  \,+\, \int_{\{z \ge a\}} z\,e^{zy} \dd\nu_t(z) \,=:\, I_1(y) + J_1(y)\,,
\]
where $|I_1(y)|/J_0(y)\to 0$ as $y \to +\infty$. It follows that
\[
  \lim_{y \to +\infty} u(t,-2y) \,=\, \lim_{y \to +\infty}\,\frac{I_1(y) +
  J_1(y)}{I_0(y) + J_0(y)} \,=\, \lim_{y \to +\infty}
  \genfrac{}{}{1pt}{0}{\int_{\{z \ge a\}} z\,e^{zy} \dd\nu_t(z)}{
  \int_{\{z \ge a\}} e^{zy} \dd\nu_t(z)} \,\ge\, a\,.
\]
Since this is true for any $a < \beta$, we deduce that $\ell_-(t) :=
\lim_{x\to-\infty}u(t,x) \ge \beta$. This means that the solution $u(t,\cdot)$
is unbounded from above if $\beta = + \infty$. In the converse case, we
must have $\ell_-(t) = \beta$, because it easily follows from \eqref{UPK-rep3}
that $u(t,x) \le \beta$ for all $(t,x) \in \Omega_-$. A symmetric argument
shows that $u(t,\cdot)$ is bounded from below if and only if $\alpha > -\infty$,
in which case $\ell_+(t) := \lim_{x\to+\infty}u(t,x) = \alpha$ for all $t < 0$. 
\end{Proof}

It is now straightforward to conclude the proof of Proposition~\ref{UPK-prop}.
If $u$ is a bounded entire solution of \eqref{Burgers}, then $u$ is a fortiori a
bounded ancient solution on $\Omega_-$, and can therefore be represented as in
\eqref{UPK-rep3} for some positive Borel measure $\mu$ that is finite on compact
intervals. By Proposition~\ref{UPK-prop2} we know that $\supp(\mu) \subset
[\alpha,\beta]$, where $\alpha,\beta \in \R$ are the spatial limits defined
in \eqref{alphabet3}. In particular $\mu$ is a finite measure, which can be
normalized into a probability measure without affecting the representation
\eqref{UPK-rep3}. We conclude that \eqref{UPK-rep} holds for all $t < 0$, hence
for all $t \in \R$ because both members are bounded solutions of Burgers'
equation which coincide on the space-time domain $\Omega_-$.

It remains to verify that the representation \eqref{UPK-rep} is unique. Assume
that $\mu_1, \mu_2$ are two probability measures on $[\alpha,\beta]$ such that
\eqref{UPK-rep} holds.  Defining
\begin{equation}\label{U1U2}
  U_1(t,x) \,=\, \int_\R e^{-zx/2 + z^2t/4} \dd\mu_1(z)\,,
  \qquad
  U_2(t,x) \,=\, \int_\R e^{-zx/2 + z^2t/4} \dd\mu_2(z)\,,
\end{equation}
we see that $U_1, U_2$ are positive solutions of the heat equation such that
$(\partial_x U_1)/U_1 = (\partial_x U_2)/U_2$ for all $(t,x) \in \R \times \R$.
This means that the ratio $r(t) := U_1(t,x)/U_2(t,x)$ does not depend on the space
variable $x$. Setting $t = 0$ in \eqref{U1U2} we deduce that
\[
  \int_\R e^{-zx/2} \dd\mu_1(z)  \,=\, r(0) \int_\R e^{-zx/2} \dd\mu_2(z)\,,
  \qquad \hbox{for all } x \in \R\,,
\]
which implies that $\mu_1 = r(0)\mu_2$ since the Laplace transform is
one-to-one. Finally, as $\mu_1, \mu_2$ are both probability measures, we
conclude that $\mu_1 = \mu_2$. \QED

\subsection{Examples: shocks, mergers, and continuous shock superposition}\label{s:52}

In this section we examine some examples of bounded entire solutions corresponding
to simple choices  for the measure $\mu$ in \eqref{UPK-rep}. As a preliminary
remark, we recall that Burgers' equation \eqref{Burgers} is invariant under
several continuous symmetries\: translations in space and time, Galilean
transformations, and parabolic scaling. It is instructive to observe, in the
case of bounded entire solutions, how the symmetry group acts on the
(not necessarily normalized) measure $\mu$. From the representation formula
\eqref{UPK-rep} we easily obtain the following group actions, where
$x_0, t_0, c \in \R$ and $\lambda > 0$\:

\begin{enumerate}[leftmargin=20pt,itemsep=2pt,topsep=2pt,label=\alph*)]

\item {\em Translation in space\:} $u(t,x) \mapsto u(t,x-x_0)$,
$\dd\mu(z) \mapsto e^{zx_0/2}\dd\mu(z)$;
  
\item {\em Translation in time\:} $u(t,x) \mapsto u(t+t_0,x)$,
$\dd\mu(z) \mapsto e^{z^2t_0/4}\dd\mu(z)$;
  
\item {\em Galilean transformation\:} $u(t,x) \mapsto u(t,x-ct) +c$,
$\dd\mu(z) \mapsto \dd\mu(z-c)$;
  
\item {\em Parabolic scaling\:} $u(t,x) \mapsto \lambda u(\lambda^2 t,\lambda x)$,
$\dd\mu(z) \mapsto \dd\mu(z/\lambda)$.  
  
\end{enumerate}

\begin{figure}[ht!]
    \includegraphics[width=0.3\textwidth]{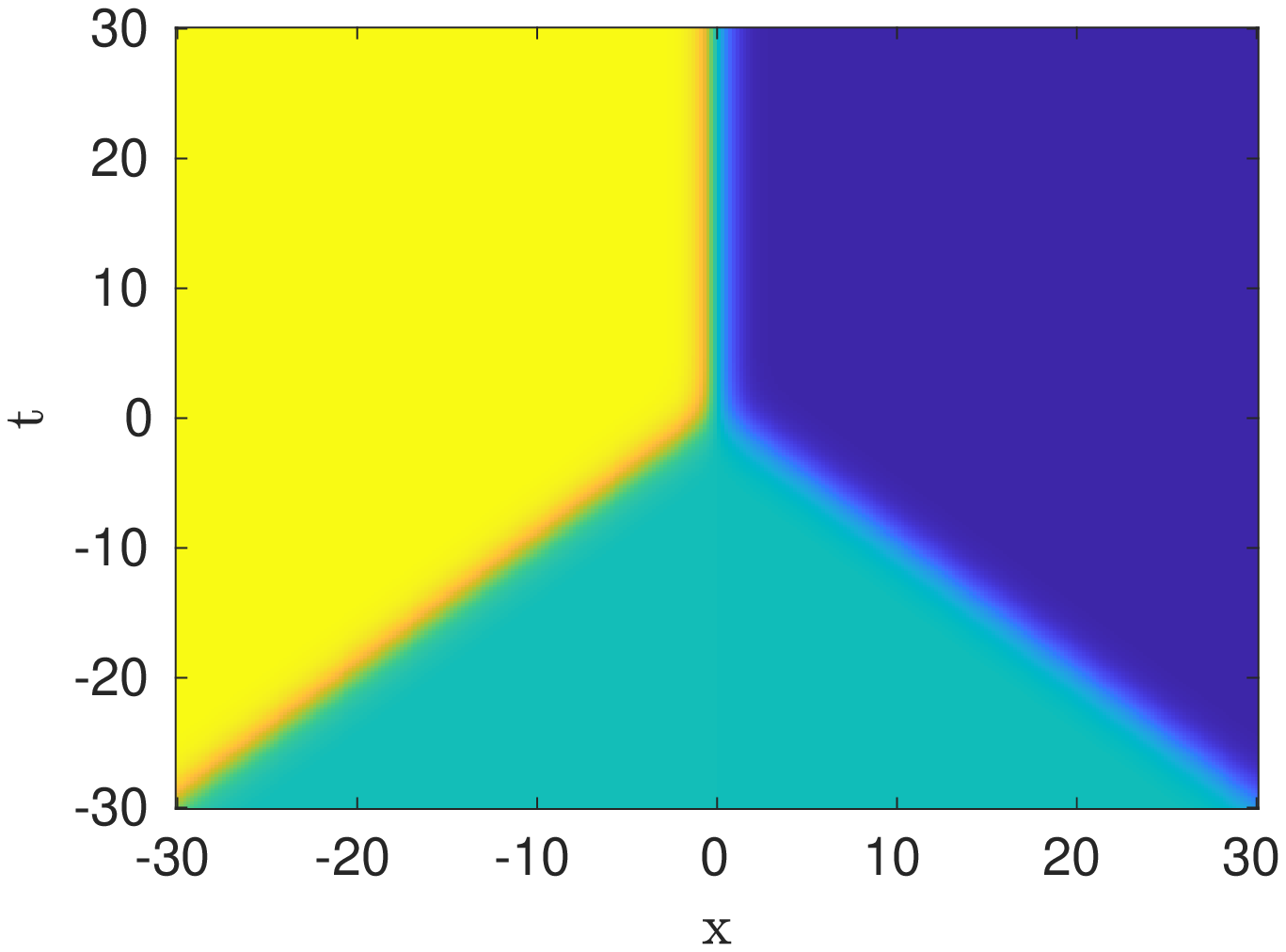}\hfill
    \includegraphics[width=0.3\textwidth]{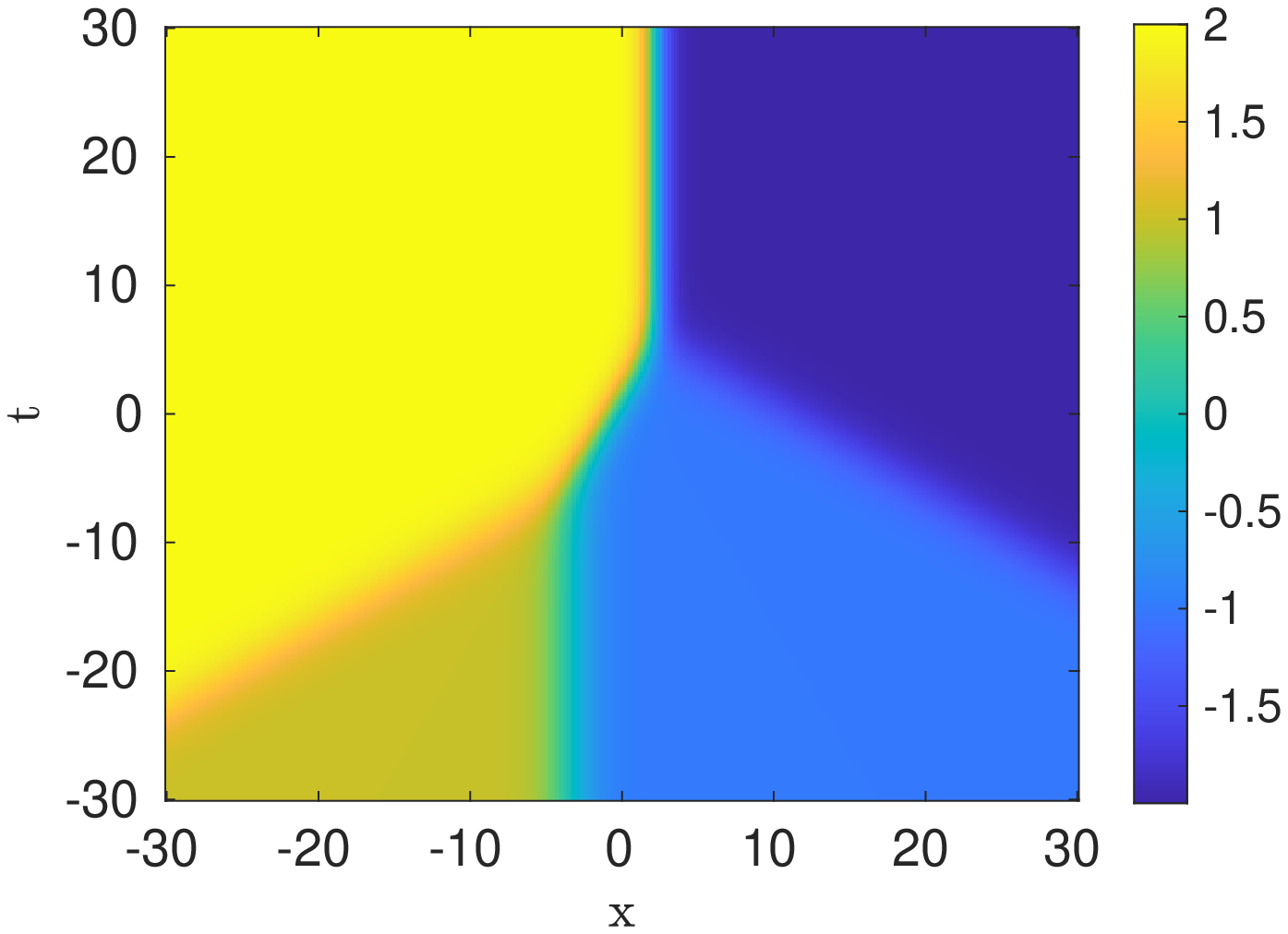}\hfill
    \includegraphics[width=0.3\textwidth]{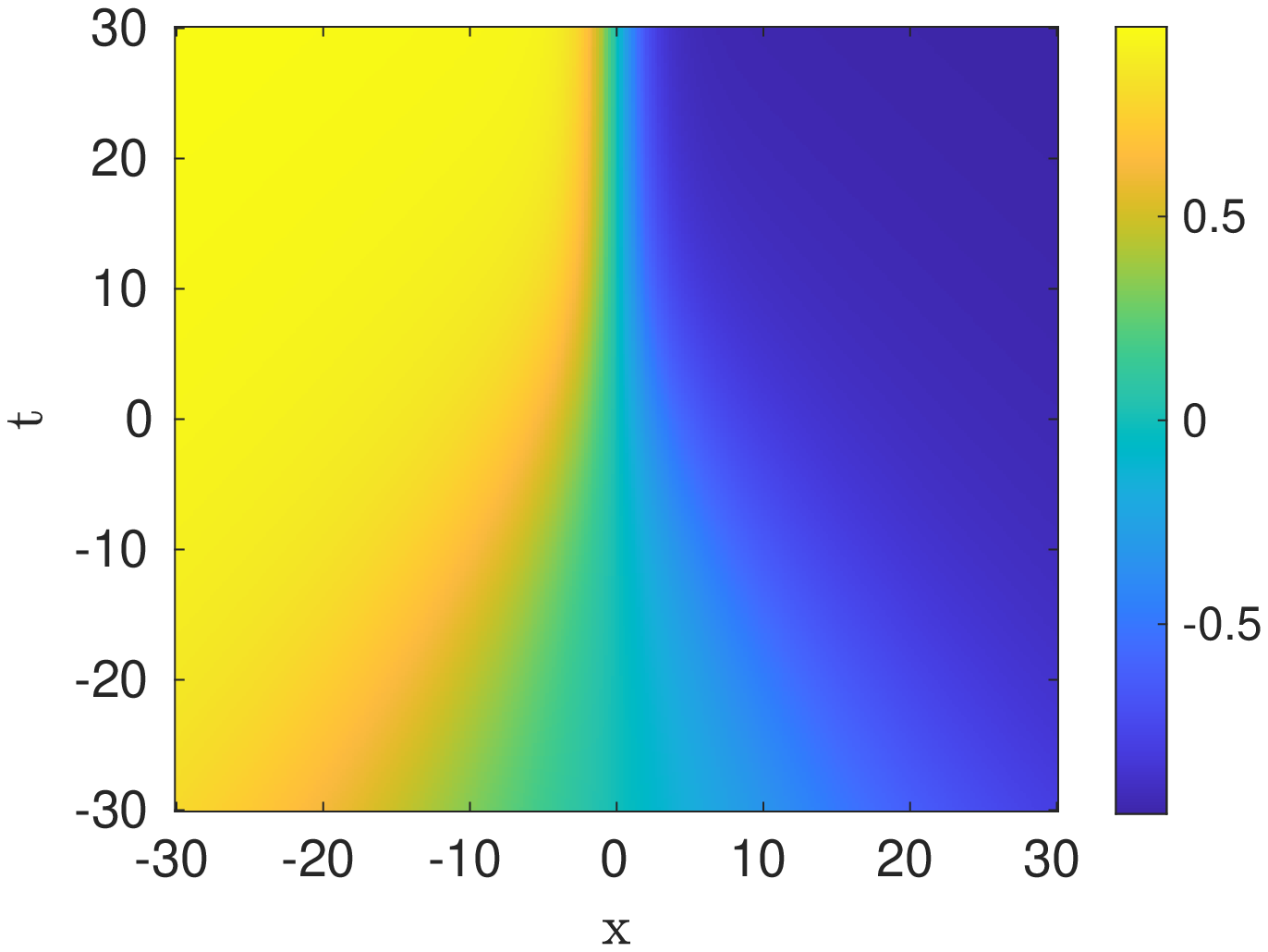}\\    
    \includegraphics[width=0.3\textwidth]{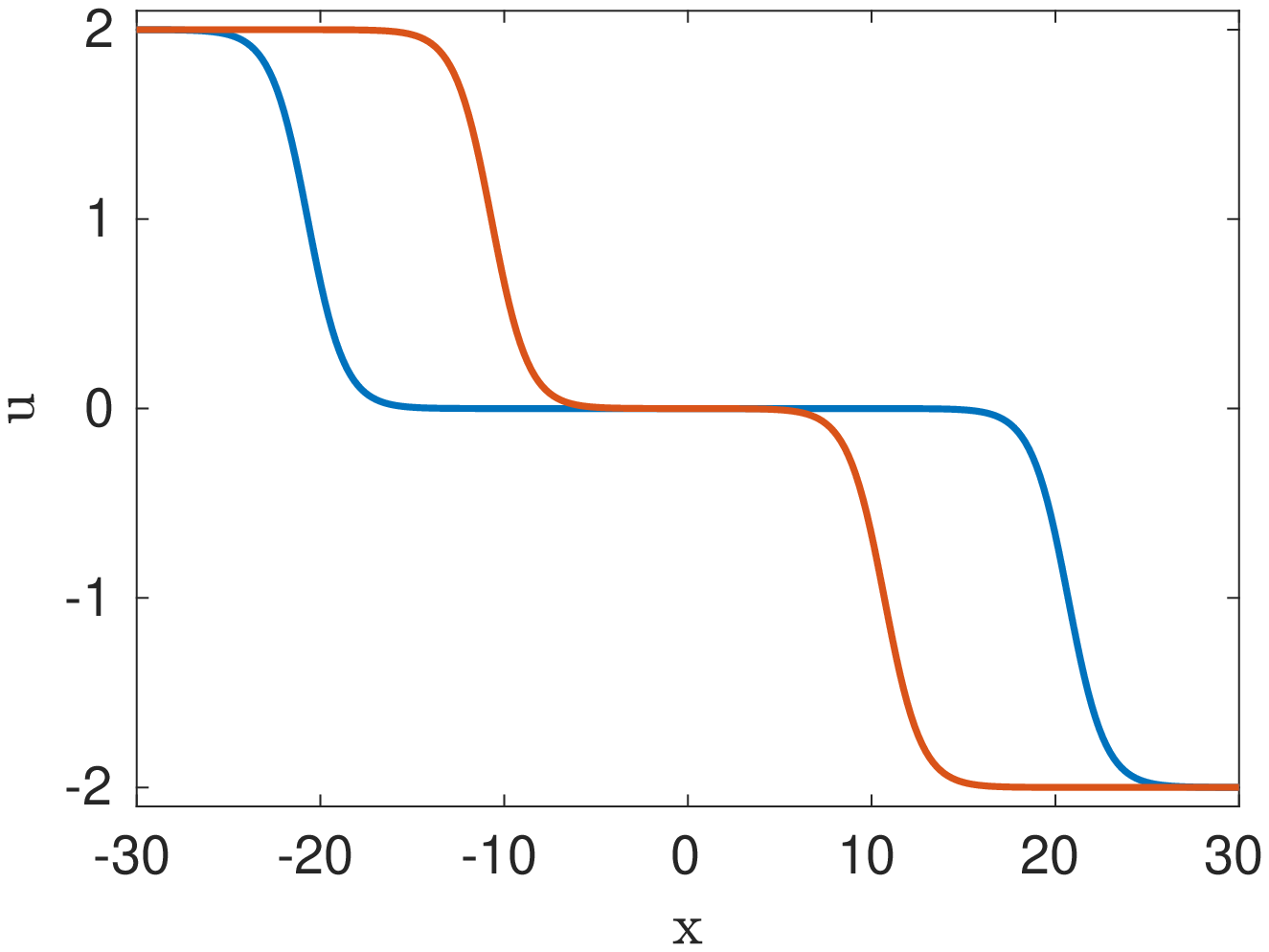}\hfill
    \includegraphics[width=0.3\textwidth]{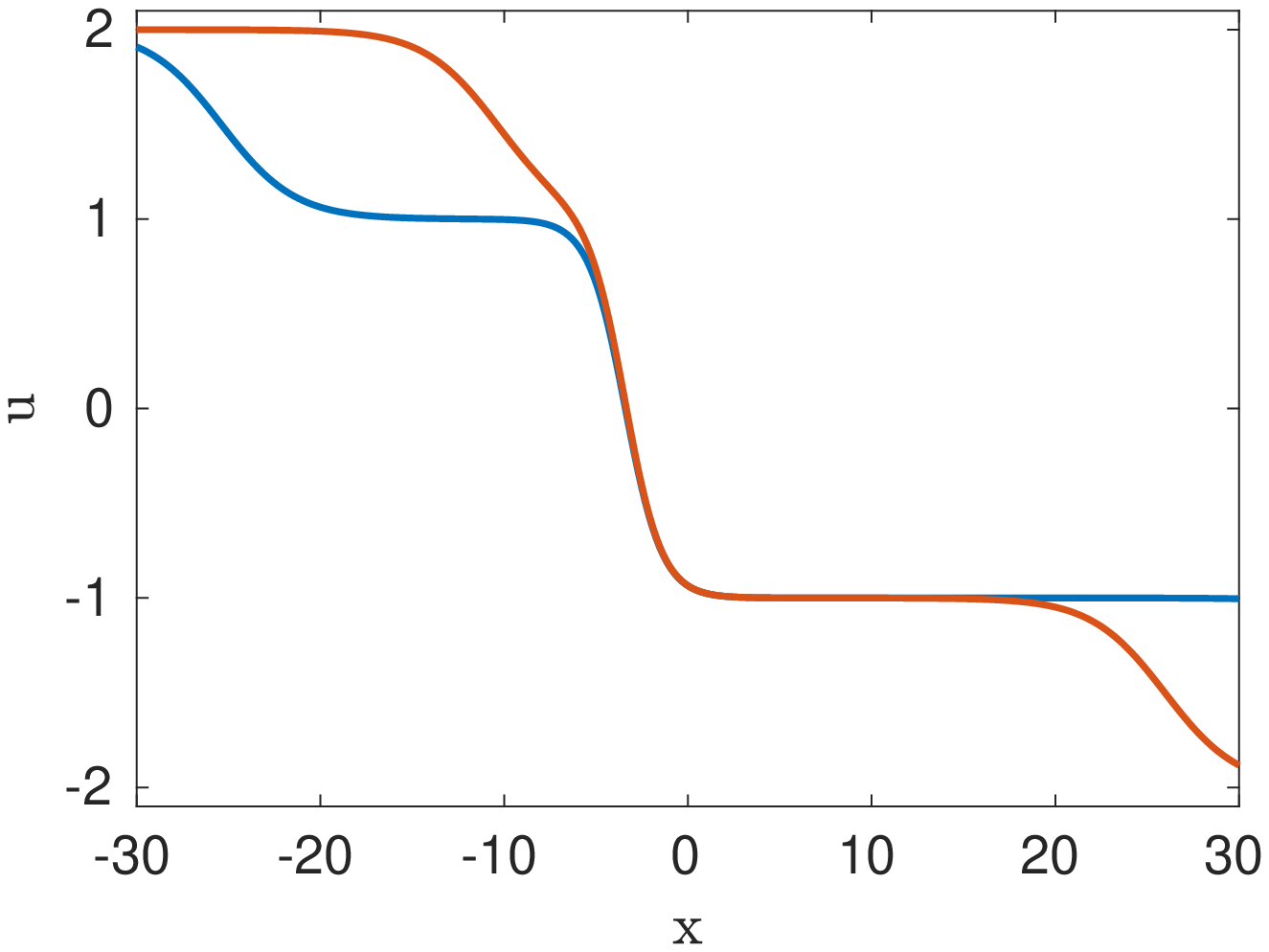}\hfill
    \includegraphics[width=0.3\textwidth]{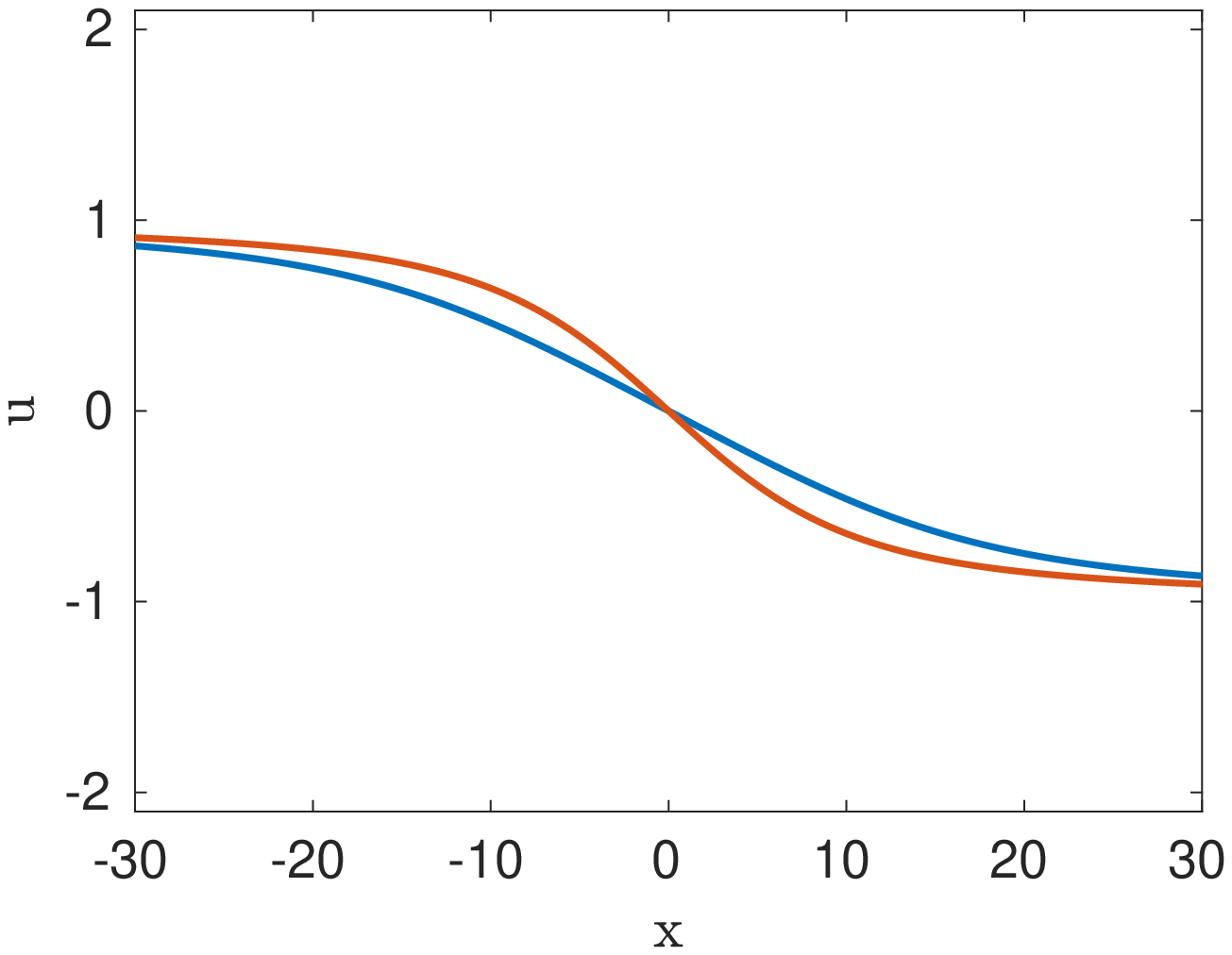}\\  
    \makebox[0.3\textwidth]{$\mu=\frac{1}{4}\delta_{-2}+\frac{1}{2}\delta_0+\frac{1}{4}\delta_{2}$}\hfill 
    \makebox[0.3\textwidth]{$\mu=\frac{1}{16}\delta_{-2}+15 \delta_{-1} +\frac{1}{2}\delta_{1}+ 5\delta_2$}
    \hfill
    \makebox[0.3\textwidth]{$\mu=\mu_L|_{[-1,1]}$}\\[0.5in]
    \includegraphics[width=0.3\textwidth]{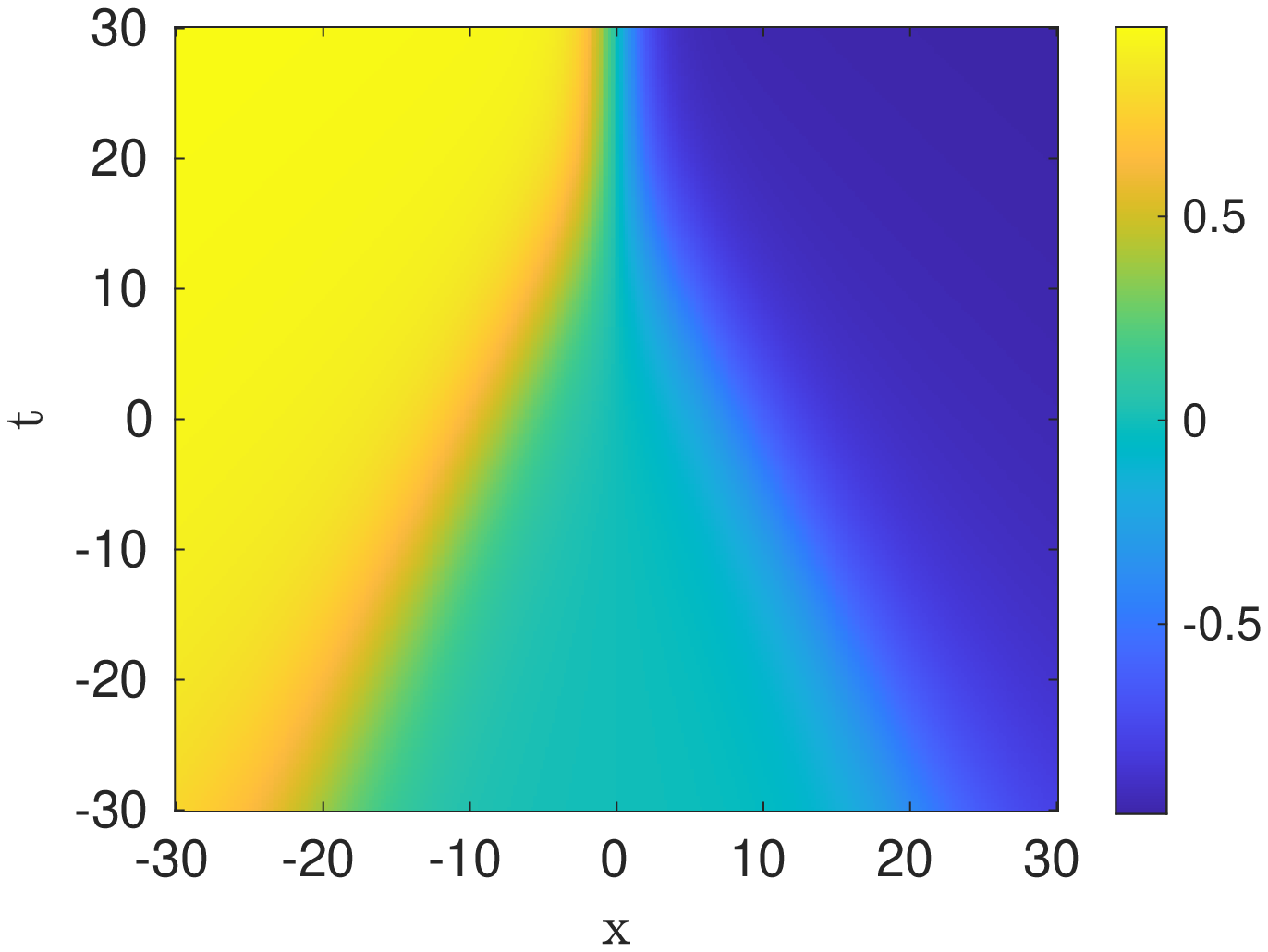}\hfill
    \includegraphics[width=0.3\textwidth]{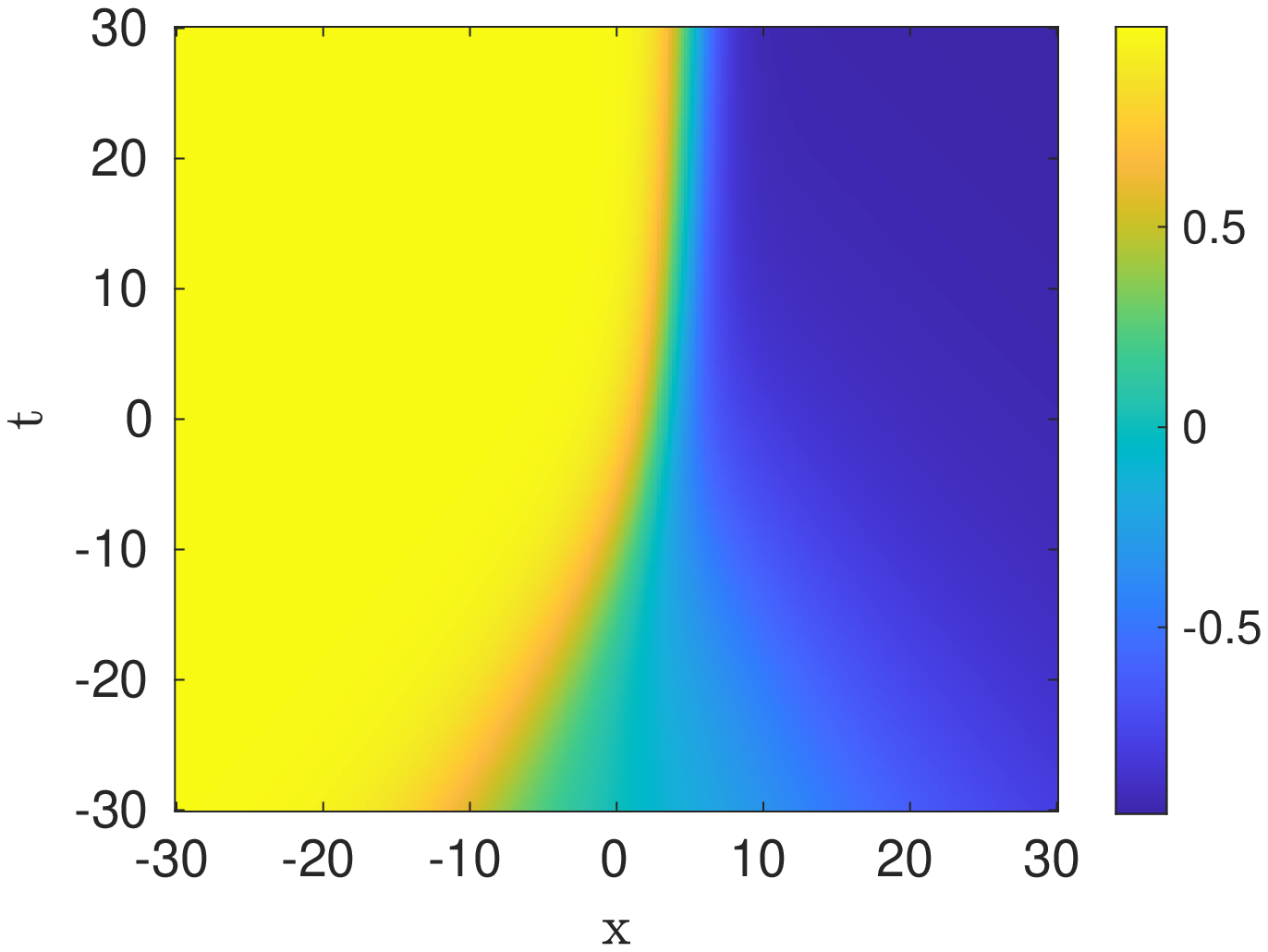}\hfill
    \includegraphics[width=0.3\textwidth]{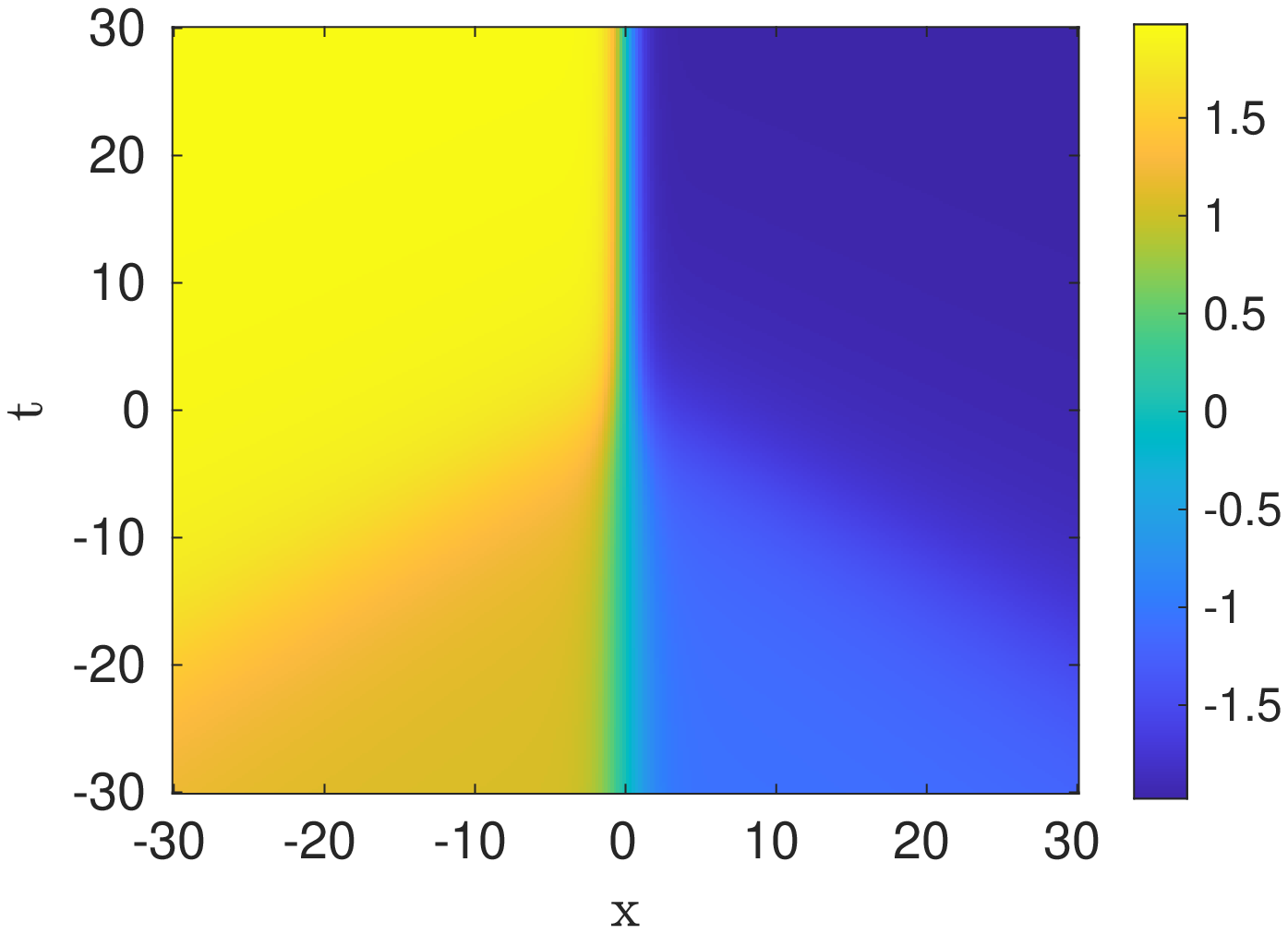}\\    
    \includegraphics[width=0.3\textwidth]{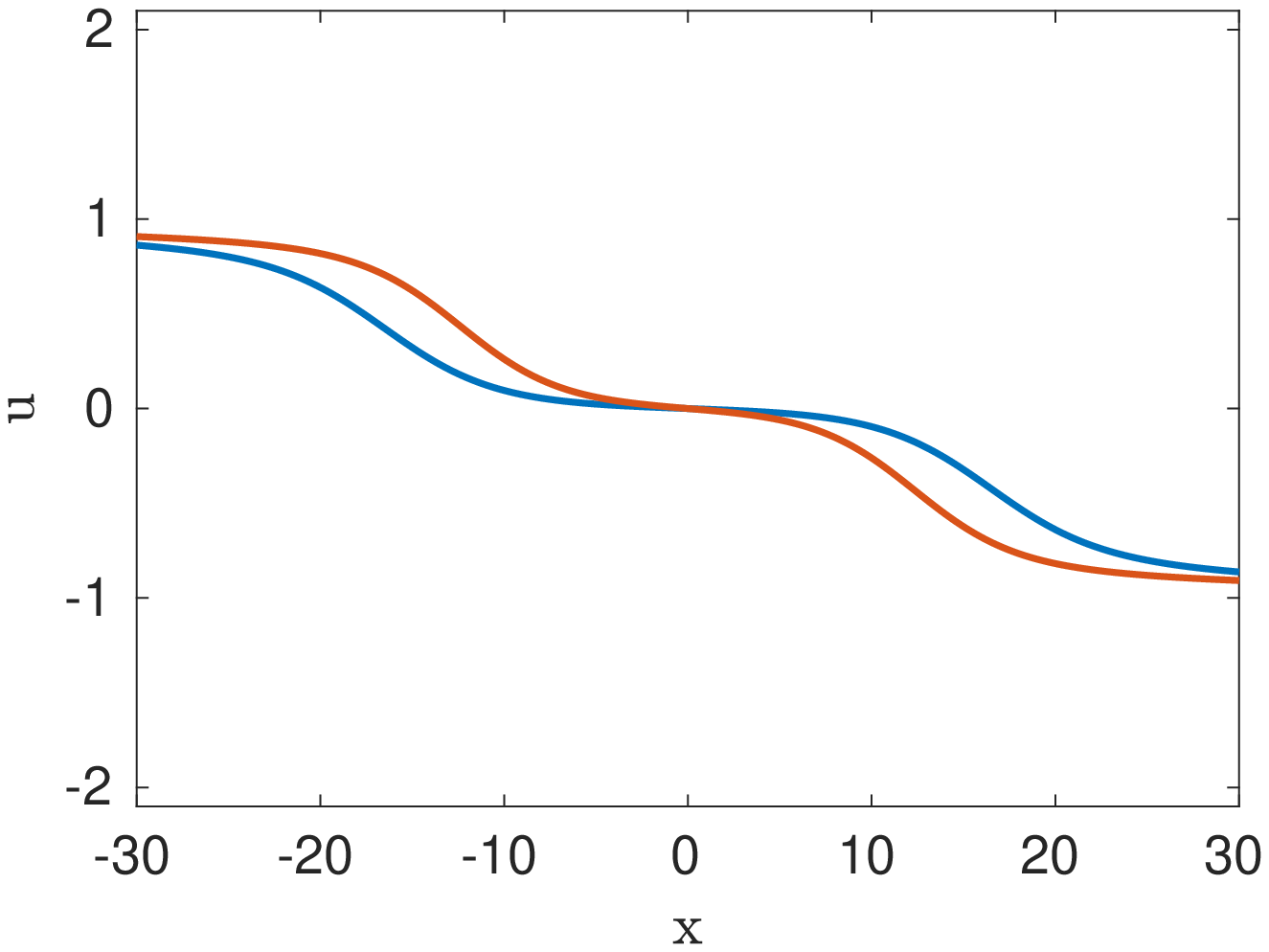}\hfill
    \includegraphics[width=0.3\textwidth]{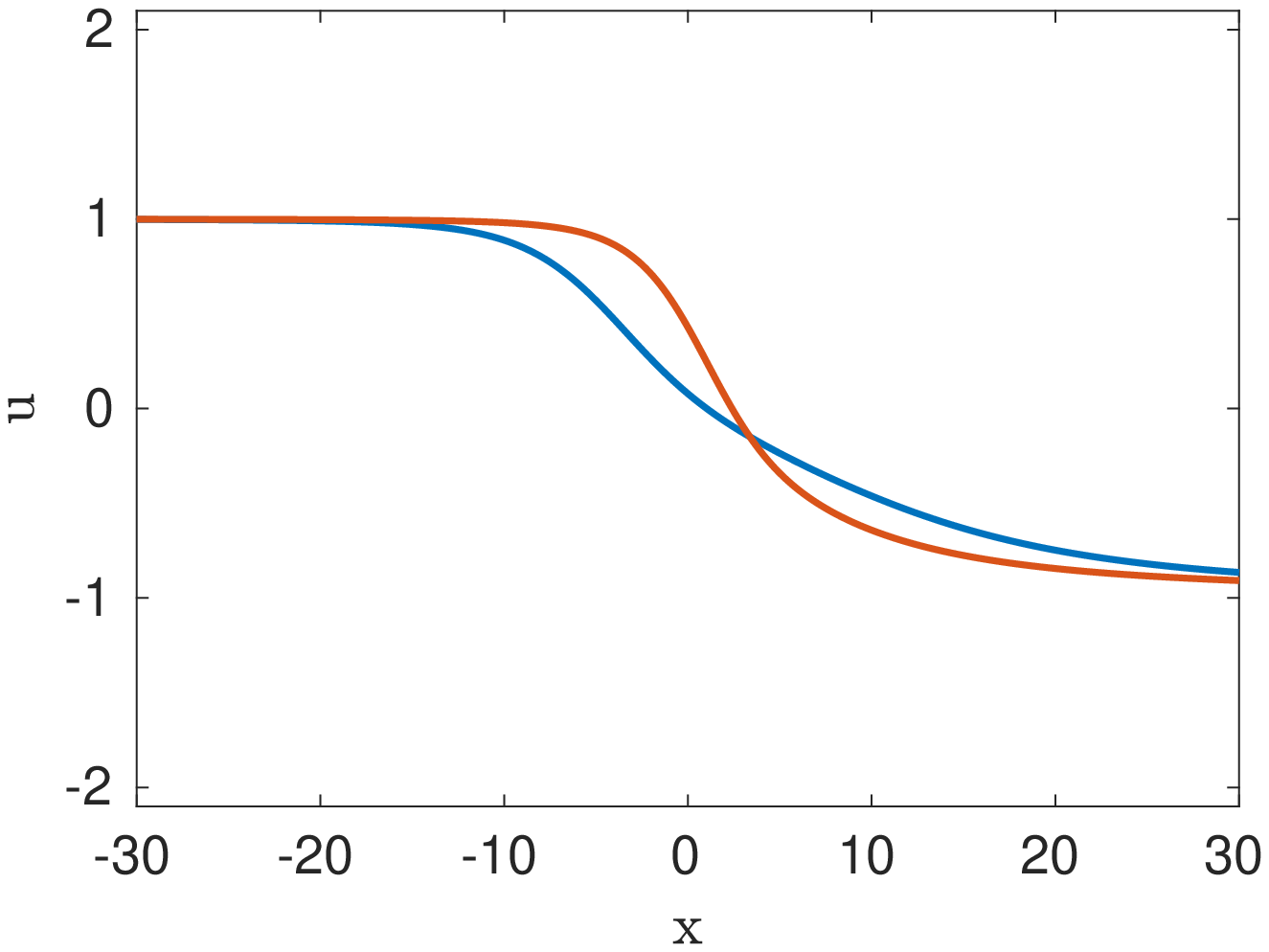}\hfill
    \includegraphics[width=0.3\textwidth]{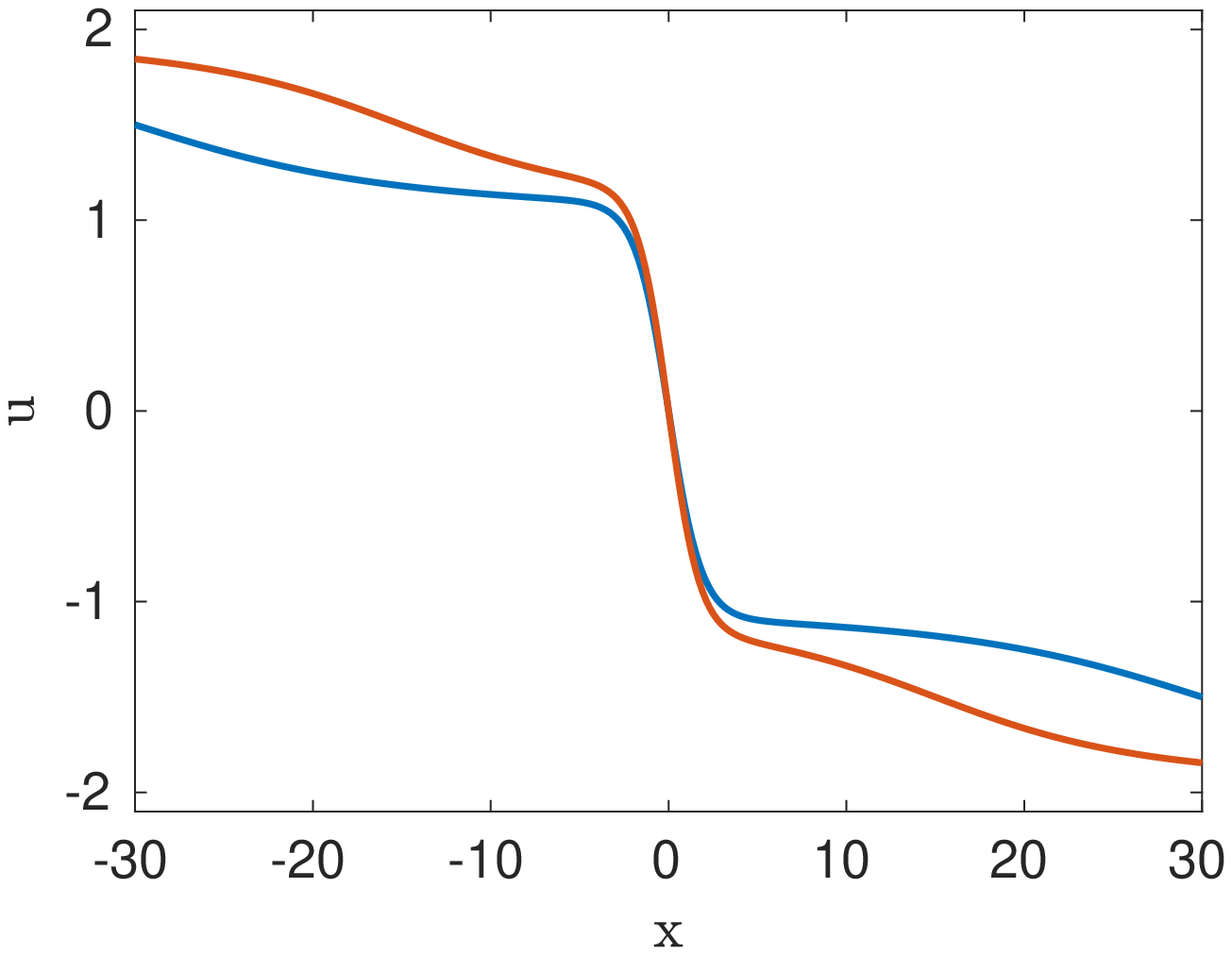} \\
    \makebox[0.3\textwidth]{$\mu=10\delta_0+\mu_L|_{[-1,1]}$}\hfill 
    \makebox[0.3\textwidth]{$\mu=10\delta_1+\mu_L|_{[-1,1]}$}
    \hfill
    \makebox[0.3\textwidth]{$\mu=\mu_L|_{[-2,-1]}+\mu_L|_{[1,2]}$}\\[0.0in]
    \caption{{\small Space-time plots of entire solutions and sample plots below at times $t=-20$ (blue) and
      $t=-10$ (red) for measures as indicated. Lebesgue's measure is denoted by $\mu_L$.}}\label{f:fig1}
\end{figure}

\medskip
In the following, we analyze several special cases of measures $\mu$. Illustrations of the corresponding entire solutions can be found in Figure~\ref{f:fig1}. 
The simplest possible example of a bounded entire solution corresponds to
$\mu = \delta_\alpha$ being a Dirac mass located at some point $\alpha \in \R$.
In that case we clearly have $u(t,x) = \alpha$ for all $(t,x) \in \R \times \R$.
A more interesting situation is obtained when $\mu = \frac12 \delta_\alpha +
\frac12 \delta_\beta$ for some $\alpha < \beta$. A direct calculation then 
shows that $u(t,x) = \phi_{\beta,\alpha}(x-ct)$, where $\phi_{\beta,\alpha}$ 
is the viscous shock profile given by
\begin{equation}\label{Burgshock}
  \phi_{\beta,\alpha}(y) \,=\, c - \delta \tanh\Bigl(\frac{\delta y}{2}\Bigr)\,,
  \qquad c \,=\, \frac{\alpha+\beta}{2}\,, \qquad \delta \,=\,
  \frac{\beta-\alpha}{2}\,.
\end{equation}
As soon as $\mu$ contains more than two Dirac masses, the solution
$u(t,x)$ given by \eqref{UPK-rep} describes the merger of several
viscous shocks into a single one. A typical example is $\mu =
\frac14 \delta_{-2} + \frac12 \delta_0 + \frac14 \delta_2$ for which
\begin{equation}\label{Burgmerger}
  u(t,x) \,=\, \frac{-2\sinh(x)}{e^{-t} + \cosh(x)}\,,
  \qquad t \in \R\,, \quad x \in \R\,.
\end{equation}
It is clear that $u(t,x) \approx \phi_{2,-2}(x)$ as $t \to +\infty$,
whereas for large negative times a direct calculation shows that
$u(t,x) \approx \phi_{2,0}(x-t+x_0) + \phi_{0,-2}(x+t-x_0)$ with
$x_0 := \log(2)$.  Thus the solution \eqref{Burgmerger} realizes the
merger of a pair of traveling viscous shocks into a single steady
shock. Mergers of more than two shocks can be described in a similar fashion. 

We next consider examples where the measure $\mu$ has an absolutely continuous
component. In analogy to finitely many Dirac masses describing discrete
superpositions of shocks and subsequent mergers, a continuous measure 
$\mu$ can be thought of as representing a continuous superposition of shocks with
continuous merger events. Such an interpretation is reminiscent of Hamel and
Nadirashvili's characterization of entire solutions to the Fisher-KPP equation in
$\R^N$, see \cite{HN}.

The prototypical example of an absolutely continuous measure $\mu$ is the
Lebesgue measure on the interval $[-1,1]$, which was considered in
\cite{UPK,KM}. In that case, the entire solution \eqref{U-rep} of the heat
equation takes the simple form
\begin{equation}\label{Uex}
  U(t,x) \,=\, \int_{-1}^1 e^{-zx/2 + z^2t/4}\dd z\,, \qquad (t,x) \in \R^2\,,
\end{equation}
so that $U(0,x) = (4/x)\sinh(x/2)$. When $t \neq 0$, we obtain after
integrating by parts
\begin{align*}
  -2\partial_x U(t,x) \,&=\, \frac{x}{t}\int_{-1}^1 e^{-zx/2 + z^2t/4}\dd z
  + \int_{-1}^1 \bigl(z-\frac{x}{t}\bigr)\,e^{-zx/2 + z^2t/4}\dd z\\
  \,&=\, \frac{x}{t}\,U(t,x) + \frac{2}{t}\,e^{t/4}\bigl(e^{-x/2} -
  e^{x/2}\bigr)\,,
\end{align*}
so that the entire solution \eqref{u-CH} of Burgers' equation
has the following expression\:
\begin{equation}\label{uex}
  u(t,x) \,=\, \frac{-2\partial_x U(t,x)}{U(t,x)} \,=\, 
  \frac{x}{t} \,-\, \frac{4}{t}\,\frac{e^{t/4}\sinh(x/2)}{
  U(t,x)}\,, \qquad t \neq 0\,, \quad x \in \R\,.
\end{equation}
It remains to obtain a more explicit formula for $U(t,x)$. When $t < 0$
a direct calculation shows that
\begin{equation}\label{Uex2}
  U(t,x) \,=\, \frac{2}{\sqrt{|t|}}\,e^{-x^2/(4t)}\biggl\{E\biggl(\frac{
  \sqrt{|t|}}{2} + \frac{x}{2\sqrt{|t|}}\biggr) + E\biggl(\frac{\sqrt{|t|}}{2}
  - \frac{x}{2\sqrt{|t|}}\biggr)\biggr\}\,, \quad (t,x) \in \Omega_-\,,
\end{equation}
where $E(x) = \int_0^x e^{-y^2}\dd y$ is the (non-normalized) error function.
Using \eqref{uex}, \eqref{Uex2}, it is not difficult to verify that
\[
  \lim_{t \to -\infty} u(t,x+ct) \,=\, \left\{\begin{array}{lcr} 1 & \hbox{if} & c > 1\,,
  \\ c & \hbox{if} & |c| \le 1\,,\\ -1 & \hbox{if} & c < -1\,, \end{array}
  \right.
\]
the convergence being uniform for $|x| \le L(t)$ provided $L(t)/|t| \to 0$
as $t \to -\infty$. This is of course in full agreement with
Propositions~\ref{p:suppmu} and \ref{p:self-similar}. For positive times,
the analogue of \eqref{Uex2} is
\begin{equation}\label{Uex3}
  U(t,x) \,=\, \frac{2}{\sqrt{t}}\,e^{t/4}\biggl\{e^{x/2} D\biggl(
  \frac{\sqrt{t}}{2} + \frac{x}{2\sqrt{t}}\biggr) + e^{-x/2}D\biggl(
  \frac{\sqrt{t}}{2}- \frac{x}{2\sqrt{t}}\biggr)\biggr\}\,,
  \quad (t,x) \in \Omega_+\,,
\end{equation}
where $D : \R \to \R$ is the Dawson function
\[
  D(x) \,=\, e^{-x^2}\int_0^x e^{y^2}\dd y \,=\, \begin{cases}
  x - \frac{2x^3}{3} + \cO(|x|^5) & \hbox{as} \quad x \to 0\,, \\
  \frac{1}{2x} + \frac{1}{4x^3} + \cO(|x|^{-5}) & \hbox{as} \quad |x| \to \infty\,.
  \end{cases}
  \]
It easily follows from \eqref{uex} and \eqref{Uex3} that $u(t,x) \to -\tanh(x/2)$
as $t \to +\infty$, in agreement with Proposition~\ref{mono-prop}. 

We next investigate how the solution is modified if the measure $\mu$ contains
in addition a Dirac mass. Assume for instance that $\mu_0 = \mu + \delta_0$,
where $\mu$ is again the Lebesgue measure on $[-1,1]$, and let $u_0(t,x)$
be the entire solution of \eqref{Burgers} associated with the measure
$\mu_0$. The same calculations as before show that 
\begin{equation}\label{u0ex}
  u_0(t,x) \,=\, \frac{x}{t}\,\frac{U(t,x)}{1+U(t,x)} \,-\, \frac{4}{t}
  \,\frac{e^{t/4}\sinh(x/2)}{1+U(t,x)}\,, \qquad (t,x) \in \R^2\,,
\end{equation}
where $U(t,x)$ is given by \eqref{Uex}. As is easily verified, the
asymptotic behavior as $t \to +\infty$ is unchanged. However, the
presence of a Dirac mass at the origin can be detected by looking
at the solution for large negative times. Indeed, a direct calculation
reveals that, for all $x \in \R$, 
\begin{equation}\label{u0u1}
  \lim_{t \to -\infty}\sqrt{|t|}\,u_0\bigl(t,x\sqrt{|t|}\bigr) \,=\, 0\,,
  \qquad \hbox{whereas}\quad
  \lim_{t \to -\infty}\sqrt{|t|}\,u\bigl(t,x\sqrt{|t|}\bigr) \,=\, -x\,.
\end{equation}

Finally, we consider the measure $\mu_1 = \mu + \delta_1$, which
includes a Dirac mass at $z=1$, and we investigate the asymptotic 
behavior of the corresponding solution $u_1(t,x)$ as $t \to +\infty$.
It is clear that
\[
  u_1(t,x) \,=\, \frac{-2\partial_x U_1(t,x)}{U_1(t,x)}\,, \qquad
  \hbox{where} \qquad U_1(t,x) \,=\, U(t,x) + e^{-x/2+t/4}\,.
\]
Using the expression \eqref{Uex3} and the asymptotic behavior of
the Dawson function at infinity, we find
\begin{align*}
  e^{-t/4}\,U_1(t,x) \,&=\, e^{-x/2}\biggl\{1 + \frac{2}{\sqrt{t}}\, D\biggl(
  \frac{\sqrt{t}}{2}- \frac{x}{2\sqrt{t}}\biggr)\biggr\} + e^{x/2}
  \frac{2}{\sqrt{t}}\, D\biggl(\frac{\sqrt{t}}{2} + \frac{x}{2\sqrt{t}}\biggr) \\
  \,&=\, e^{-x/2}\Bigl(1 + \frac{2}{t} + \cO(t^{-2})\Bigr) +
  e^{x/2}\Bigl(\frac{2}{t} + \cO(t^{-2})\Bigr)\,, \quad t \to +\infty\,.
\end{align*}
Defining $\bar x(t) = \log(1+t/2)$, we see that $\sqrt{t}\,e^{-t/4}
\,U_1(t,x + \bar x(t)) \to 2\sqrt{2}\cosh(x/2)$ as $t \to +\infty$,
and we conclude that $u_1(t,x+\bar x(t)) \to -\tanh(x/2)$ as
$t \to +\infty$. In other words, the presence of a Dirac mass
at $z=1$ is responsible for a logarithmic shift in the position
of the viscous shock, as discussed in Remark~\ref{mono-rem}.

\subsection{Asymptotic analysis of entire solutions as $t\to-\infty$}\label{s:53}

The properties of the measure $\mu$ in \eqref{UPK-rep} are reflected in the
asymptotic behavior of the entire solution $u(t,x)$ in the ancient limit
$t \to -\infty$. Some results in this direction were already stated in
Propositions~\ref{p:suppmu} and \ref{p:self-similar}, and illustrated by
the examples considered in the previous section. Our goal here is to 
perform a more systematic study of the ancient limit for entire solutions of
\eqref{Burgers}. Our main results are Propositions~\ref{mu-prop1} and
\ref{mu-prop2} below, which immediately imply the statements given in
the introduction, and also extend the results obtained in \cite[Section~7]{UPK}.

To gain a first intuitive understanding, we consider the entire solution
\eqref{UPK-rep} in a Galilean frame moving with speed $c \in \R$. In the spirit
of Appel's transformation \eqref{V-Ap}, we also introduce the inverse time
$\tau = -1/t$, so that the ancient limit $t \to -\infty$ becomes a standard
short time limit in the new variable $\tau$. A simple calculation shows that
\begin{equation}\label{uheat}
  u(t,x+ct) \,=\, \frac{\int z e^{-zx/2+(z-c)^2t/4}\dd\mu(z)}{\int e^{-zx/2+(z-c)^2t/4}
  \dd\mu(z)} \,=\, \frac{\int K(\tau,c-z)\,z\,e^{-zx/2}\dd\mu(z)}{\int
  K(\tau,c-z)\,e^{-zx/2}\dd\mu(z)}\,,
\end{equation}
where $K(t,x)$ is the heat kernel \eqref{K-def}. We investigate the behavior
of \eqref{uheat} in the limit $\tau \to 0+$, for a fixed $x \in \R$ (for
simplicity, we assume here that $x = 0$). The denominator in the right-hand
side of \eqref{uheat} is exactly the solution at time $\tau > 0$ of the heat
equation with initial data $\mu$, evaluated at point $c \in \R$. When
$\tau$ is small, this quantity is an average of the measure $\mu$ in a small
neighborhood of the point $c$. The numerator has a similar interpretation,
except that the initial measure is now $z\dd\mu(z)$.

These observations strongly suggest that $u(t,ct)$ should converge to
$c$ as $t \to -\infty$, whenever $c$ belongs to the support of the measure
$\mu$. If $c \notin \supp(\mu)$, we expect that the ancient limit of $u(t,ct)$
will  depend on the behavior of the measure near the point in $\supp(\mu)$
that is closest to $c$. The results established below show that these heuristic
considerations are indeed correct.

In what follows, we always assume that $\mu$ is probability measure
supported in a bounded interval of $\R$, and we denote by $u(t,x)$ the
bounded entire solution of \eqref{Burgers} given by \eqref{UPK-rep}.  
We first consider the case where $c \in \supp(\mu)$. 

\begin{prop}\label{mu-prop1}
If $c \in \supp(\mu)$, the entire solution of \eqref{Burgers} defined by
\eqref{UPK-rep} satisfies
\begin{equation}\label{unif1}
  \sup_{|x| \le L(t)} \bigl|u(t,x+ct) - c\bigr|
  \,\xrightarrow[t \to -\infty]{}\, 0\,,
\end{equation}
where $L : \R_- \to \R_+$ is any function such that $L(t)/|t| \to 0$ as
$t \to -\infty$. 
\end{prop}

\begin{Proof}
By Galilean invariance, it is sufficient to prove \eqref{unif1} in the
particular case where $c = 0$. We proceed as in the proof of
Proposition~\ref{UPK-prop2}. Given $\epsilon > 0$ we observe that
\begin{equation}\label{urepIJ}
  u(t,x) \,=\, \frac{I_1(t,x) + J_1(t,x)}{I_0(t,x) + J_0(t,x)}\,,
  \qquad t \in \R\,, \quad x \in \R\,, 
\end{equation}
where for $k \in \{0,1\}$ we denote
\[
  I_k(t,x) \,=\, \int_{\{|z| > \epsilon\}} z^k\,e^{-zx/2+z^2t/4}\dd\mu(z)\,, \qquad
  J_k(t,x) \,=\, \int_{\{|z| \le \epsilon\}}z^k\,e^{-zx/2+z^2t/4}\dd\mu(z)\,.
\]
Assuming that $|x| \le L$ and recalling that $t < 0$, we observe that
$-zx/2+z^2t/4 \le -L^2/t + z^2t/6$, hence
\begin{equation}\label{Ibd1}
  I_k(t,x) \,\le\, e^{-L^2/t + \epsilon^2 t/6} \int_{\R}|z|^k\dd\mu(z)\,.
\end{equation}
On the other hand, we obviously have
\begin{equation}\label{Ibd2}
  J_0(t,x) \,\ge\, \int_{\{|z| \le \epsilon/2\}}e^{-zx/2+z^2t/4}\dd\mu(z)
  \,\ge\, e^{-\epsilon L/4 + \epsilon^2 t/16}\,\mu\bigl([-\epsilon/2,\epsilon/2]\bigr)\,,
\end{equation}
where $\mu\bigl([-\epsilon/2,\epsilon/2]\bigr) > 0$ since $0 \in \supp(\mu)$. Taking
$L = L(t)$ with $L(t) = o(|t|)$, we deduce from \eqref{Ibd1}, \eqref{Ibd2} that
$I_k(t,x)/J_0(t,x)$ converges to zero as $t \to -\infty$, uniformly for $|x| \le L(t)$.
If we now return to \eqref{urepIJ}, we conclude that
\begin{equation}\label{limsup}
  \limsup_{t \to -\infty}\sup_{|x| \le L(t)}\bigl|u(t,x)\bigr|  \,\le\,
  \limsup_{t \to -\infty} \sup_{|x| \le L(t)}\,\frac{|J_1(t,x)|}{J_0(t,x)}
  \,\le\, \epsilon\,.
\end{equation}
Since $\epsilon > 0$ was arbitrary, the left-hand side in \eqref{limsup}
actually vanishes, which gives \eqref{unif1}. 
\end{Proof}

\begin{rem}\label{lin-rem}
The assumption that $L(t)/|t| \to 0$ as $t\to -\infty$ is optimal in general. As can be
seen from Example~\ref{Lebex}, where $\mu$ is the Lebesgue measure on $\R$, the
conclusion \eqref{unif1} fails for any $c \in \R$ if $L(t)/|t|$ does not converge to zero.
However, if the measure $\mu$ has an atom at $c \in \R$ which is an isolated point
in $\supp(\mu)$, it is easy to verify that \eqref{unif1} holds with $L(t) = \epsilon |t|$
for some $\epsilon > 0$. 
\end{rem}
  
The case where $c \notin \supp(\mu)$ is more difficult to treat, see
\cite[Lemma~7.2]{UPK} for an attempt in this direction. For simplicity, we
formulate our result in the case where $c = 0$, but as already mentioned
this does not restrict the generality. 

\begin{prop}\label{mu-prop2}
Assume that $0 \notin \supp(\mu)$, and define $a \in [-\infty,0)$ and $b \in
(0,+\infty]$ by
\begin{equation}\label{abdef}
  a \,:=\, \sup\bigl\{c < 0\,;\, c \in \supp(\mu)\bigr\}\,, \qquad
  b \,:=\, \inf\bigl\{c > 0\,;\, c \in \supp(\mu)\bigr\}\,. 
\end{equation}
\begin{enumerate}[leftmargin=20pt,itemsep=2pt,topsep=2pt,label=\roman*)]

\item If $a+b \neq 0$ the ancient solution of \eqref{Burgers} given by
\eqref{UPK-rep} satisfies
\begin{equation}\label{unif2}
  \sup_{|x| \le L(t)} \bigl|u(t,x) - d\bigr| \,\xrightarrow[t \to -\infty]{}\, 0\,,
\end{equation}
where $L(t)$ is as in \eqref{unif1} and $d = b$ if $a+b < 0$, $d = a$ if $a+b > 0$. 

\item If $a+b = 0$ there exists a shift function $s : \R \to \R$ such that
$s(t)/t \to 0$ as $t \to -\infty$ and
\begin{equation}\label{unif3}
  \sup_{|x| \le L(t)} \bigl|u(t,x) - \phi_{b,a}(x-s(t))\bigr| \,\xrightarrow[t \to
  -\infty]{}\, 0\,,
\end{equation}
where $\phi_{b,a}$ denotes the viscous shock connecting the constant states $b$
and $a = -b$, see \eqref{Burgshock}. 
\end{enumerate}
\end{prop}

\begin{rem}
Since the measure $\mu$ in \eqref{UPK-rep} is nontrivial, even when the solution $u(t,x)$
vanishes identically, the quantities $a,b$ defined in \eqref{abdef} cannot be infinite
simultaneously. It follows that the sum $a + b \in [-\infty,\infty]$ is well defined,
and so is $d \in \R$ in case (i). In the other case, both $a$ and $b$ are finite.
\end{rem}

\begin{Proof}
The proof of case (i) uses exactly the same arguments as in Proposition~\ref{mu-prop1}.
Assume for instance that $a + b < 0$, so that $b$ is the point in $\supp(\mu)$ that
is closest to the origin. If $t < 0$ is large, the leading contributions in the
representation formula \eqref{UPK-rep} correspond to the restriction of the measure
$\mu$ to a small neighborhood of $b$. More precisely, given any $\epsilon > 0$, we find
\[
  \limsup_{t \to -\infty}\sup_{|x| \le L(t)}\bigl|u(t,x) - b\bigr|  \,\le\, 
  \limsup_{t \to -\infty} \sup_{|x| \le L(t)} \genfrac{}{}{1pt}{0}{\int_{N_\epsilon}
   |z-b|\,e^{-zx/2 + z^2t/4} \dd\mu(z)}{\int_{N_\epsilon} e^{-zx/2 + z^2t/4} \dd\mu(z)}
  \,\le\, \epsilon\,,
\]
where $N_\epsilon = [b,b+\epsilon]$ and $L(t)$ is as in \eqref{unif1}. This gives
\eqref{unif2} when $a+b < 0$, and the other case is treated similarly.

We now concentrate on the case (ii) where $a = -b$, which requires a more
careful analysis because both intervals $N_\epsilon^+ := [b,b{+}\epsilon]$ and
$N_\epsilon^- := [-b{-}\epsilon,-b]$ equally contribute to the representation
formula \eqref{UPK-rep} when $t < 0$ is large. Denoting
\[
  v_\epsilon(t,x) \,=\, \frac{I_1^+(t,x) + I_1^-(t,x)}{I_0^+(t,x) + I_0^-(t,x)}\,,
  \qquad \hbox{where} \quad
  I_k^\pm(t,x) \,=\, \int_{N_\epsilon^\pm} z^k\,e^{-zx/2+z^2t/4}\dd\mu(z)\,,
\]
we easily find that $|u(t,x) - v_\epsilon(t,x)|$ converges to zero as $t \to -\infty$
uniformly for $|x| \le L(t)$, provided $L(t)/|t| \to 0$ as $t \to -\infty$. So it
remains to determine the behavior of $v_\epsilon(t,x)$ for large negative times. 

For this purpose we first observe that
\[
  \bigl|I_1^+(t,x) - b I_0^+(t,x)\bigr| \,\le\, \epsilon I_0^+(t,x)\,, \qquad
  \bigl|I_1^-(t,x) + b I_0^-(t,x)\bigr| \,\le\, \epsilon I_0^-(t,x)\,,
\]
so that $|v_\epsilon(t,x) - w_\epsilon(t,x)| \le \epsilon$ where
\begin{equation}\label{waux}
  w_\epsilon(t,x) \,:=\, b\,\frac{I_0^+(t,x) - I_0^-(t,x)}{I_0^+(t,x) + I_0^-(t,x)}\,,
  \qquad t \in \R\,, \quad x \in \R\,.
\end{equation}
Using the change of variables $z = \pm(b+y)$ we can write
\begin{equation}\label{IJaux}
  \begin{split}
    I_0^+(t,x) \,&=\, e^{b^2t/4}\,e^{-bx/2} J^+_\epsilon(t,x)\,, \qquad J^+_\epsilon(t,x)
    \,:=\, \int_{[0,\epsilon]} e^{-xy/2}\,e^{t(by/2 + y^2/4)}\dd\nu_+(y)\,, \\
    I_0^-(t,x) \,&=\, e^{b^2t/4}\,e^{bx/2} J^-_\epsilon(t,x)\,, \hspace{28pt} J^-_\epsilon(t,x)
    \,:=\, \int_{[0,\epsilon]} e^{xy/2}\,e^{t(by/2 + y^2/4)}\dd\nu_-(y)\,,
  \end{split}
\end{equation}
where $\nu_\pm$ are positive measures on $[0,\epsilon]$ with $0 \in \supp(\nu_\pm)$.
If we substitute \eqref{IJaux} in \eqref{waux} we obtain the nicer expression
\begin{equation}\label{waux2}
  w_\epsilon(t,x) \,=\, -b\tanh\Bigl(\frac{b}{2}\bigl(x - S_\epsilon(t,x)\bigr)\Bigr)\,,
  \qquad \hbox{where}\quad S_\epsilon(t,x) \,=\, \frac{1}{b}\,
  \log\frac{J^+_\epsilon(t,x)}{J^-_\epsilon(t,x)}\,.
\end{equation}
Summarizing the results obtained so far, we have shown that, given any $\epsilon > 0$,
\begin{equation}\label{uwconv}
  \limsup_{t \to -\infty}\sup_{|x| \le L(t)}\bigl|u(t,x) - \phi\bigl(x - S_\epsilon(t,x)
  \bigr)\bigr| \,\le\, \epsilon\,,
\end{equation}
where $\phi = \phi_{b,a}$ is the viscous shock \eqref{Burgshock} connecting the states
$b$ and $a = -b$. The bound \eqref{uwconv} holds provided $L(t)/|t| \to 0$ as
$t \to -\infty$. 

To go further, we need properties of the shift function $S_\epsilon(t,x)$ that are
established in Section~\ref{ssecA4}. 

\begin{lem}\label{Seps-lem}
The shift function $S_\epsilon$ defined in \eqref{waux2} satisfies the uniform
bounds
\begin{equation}\label{Sepsder}
  |\partial_x S_\epsilon(t,x)| \,\le\, \frac{\epsilon}{b}\,, \qquad |\partial_t
  S_\epsilon(t,x)| \,\le\, \Bigl(\epsilon + \frac{\epsilon^2}{2b}\Bigr)\,,
  \qquad t \in \R\,, \quad x \in \R\,.
\end{equation}
Moreover $S_\epsilon$ is independent of the parameter $\epsilon$ in the ancient
limit, in the sense that 
\begin{equation}\label{S2epsconv}
  \lim_{t \to -\infty}\sup_{|x| \le L(t)}\bigl|S_\epsilon(t,x) - S_{\epsilon'}(t,x)
  \bigr| \,=\, 0\,, \qquad \hbox{if}\quad 0 < \epsilon' < \epsilon\,,
\end{equation}
provided $L(t)/|t| \to 0$ as $t \to -\infty$. 
\end{lem}

In the rest of the proof, we assume without loss of generality that $0 < \epsilon < b$.
For any $t \in \R$, we denote by $s_\epsilon(t)$ the unique real number satisfying
\begin{equation}\label{sepsdef}
  s_\epsilon(t) \,=\, S_\epsilon(t,s_\epsilon(t))\,.
\end{equation}
Note that $|\partial_x S(t,x)| \le \epsilon/b < 1$ by \eqref{Sepsder}, so that
the equation $x = S_\epsilon(t,x)$ for $x \in \R$ has indeed a unique solution. The
following properties of $s_\epsilon$ will also be established in Section~\ref{ssecA4}. 

\begin{lem}\label{seps-lem}
If $0 < \epsilon < b$ the function $s_\epsilon : \R \to \R$ defined by
\eqref{sepsdef} satisfies
\begin{equation}\label{sepslim}
  \lim_{t \to -\infty} \frac{s_\epsilon(t)}{t} \,=\, 0\,, \qquad \hbox{and}\qquad
  \lim_{t \to -\infty}\,\bigl|s_\epsilon(t) - s_{\epsilon'}(t)\bigr| \,=\, 0
  \quad \hbox{for all } \epsilon' \in (0,\epsilon)\,.
\end{equation}
\end{lem}

Equipped with Lemmas~\ref{Seps-lem} and \ref{seps-lem}, we now conclude the
proof of Proposition~\ref{mu-prop2}. For $t < 0$ large enough and $|x| \le L(t)$,
we want to estimate the quantity
\begin{equation}\label{usdiff}
  \bigl|u(t,x) - \phi(x - s_\epsilon(t))\bigr| \,\le\, \bigl|u(t,x) -
  \phi\bigl(x - S_\epsilon(t,x)\bigr)\bigr| + \bigl|\phi\bigl(x -
  S_\epsilon(t,x)\bigr) - \phi\bigl(x - s_\epsilon(t)\bigr)\bigr|\,.
\end{equation}
The first term in the right-hand side is controlled using \eqref{uwconv}.
To bound the second one, we consider three different regions:

\begin{enumerate}[leftmargin=20pt,itemsep=2pt,topsep=4pt,label=\arabic*)]

\item When $|x - s_\epsilon(t)| \le \epsilon^{-1/2}$, we have
\begin{align*}
  \bigl|\phi\bigl(x - S_\epsilon(t,x)\bigr) - \phi\bigl(x - s_\epsilon(t)\bigr)\bigr|
  \,&\le\, \|\phi'\|_{L^\infty} \bigl|S_\epsilon(t,x) - S_\epsilon(t,s_\epsilon(t))
  \bigr|  \\
  \,&\le\, \|\phi'\|_{L^\infty} \,\frac{\epsilon}{b}\,|x - s_\epsilon(t)|
   \,\le\, \frac{b \epsilon^{1/2}}{2}\,,
\end{align*}
because $\|\phi'\|_{L^\infty} = b^2/2$ and $|\partial_x S_\epsilon| \le \epsilon/b$
by \eqref{Sepsder}.

\item When $x - s_\epsilon(t) \ge \epsilon^{-1/2}$, the triangle inequality implies
that
\[
  x - S_\epsilon(t,x) \,\ge\, x - s_\epsilon(t) - |S_\epsilon(t,x) - S_\epsilon(t,s_\epsilon(t))|
  \,\ge\, x - s_\epsilon(t) - \frac{\epsilon}{b} \bigl(x - s_\epsilon(t)\bigr)\,,
\]
so that $x - S_\epsilon(t,x) \ge \epsilon^{-1/2}(1-\epsilon/b)$. It follows
that
\[
  \bigl|\phi\bigl(x - S_\epsilon(t,x)\bigr) - \phi\bigl(x - s_\epsilon(t)\bigr)\bigr|
  \,\le\, \bigl|\phi\bigl(x - S_\epsilon(t,x)\bigr) + b\bigr| +
  \bigl|\phi\bigl(x - s_\epsilon(t)\bigr) + b\bigr| \,=\, \cO\bigl(e^{-b\epsilon^{-1/2}}
  \bigr)\,,
\]
because $\phi(y) + b = b\bigl(1-\tanh(by/2)\bigr) \sim 2b\,e^{-by}$ as $y \to +\infty$.

\item The same bound holds when $x - s_\epsilon(t) \le -\epsilon^{-1/2}$, and
is established by a similar argument. 
\end{enumerate}

\noindent
Summarizing, we have shown that there exists a constant $C > 0$ such that
\begin{equation}\label{phidiff}
  \sup_{x \in \R}\,\bigl|\phi\bigl(x - S_\epsilon(t,x)\bigr) - \phi\bigl(x - s_\epsilon(t)
  \bigr)\bigr| \,\le\, C \epsilon^{1/2}\,, \qquad \hbox{if} \quad 0 < \epsilon < b\,.
\end{equation}
If we now combine \eqref{uwconv}, \eqref{usdiff}, and \eqref{phidiff}, we arrive
at
\begin{equation}\label{uphiaux}
  \limsup_{t \to -\infty}\sup_{|x| \le L(t)}\bigl|u(t,x) - \phi\bigl(x - s_\epsilon(t)
  \bigr)\bigr| \,\le\, \epsilon + C\epsilon^{1/2}\,.
\end{equation}
In fact, since $|u(t,x) - \phi\bigl(x - s_\epsilon(t)\bigr)| \le |u(t,x) - \phi\bigl(x
- s_{\epsilon'}(t)\bigr)| + \|\phi'\|_{L^\infty} |s_\epsilon(t) - s_{\epsilon'}(t)|$, it
follows from \eqref{sepslim}, \eqref{uphiaux} that 
\begin{equation}\label{uphiaux2}
  \limsup_{t \to -\infty}\sup_{|x| \le L(t)}\bigl|u(t,x) - \phi\bigl(x - s_\epsilon(t)
  \bigr)\bigr| \,\le\, \epsilon' + C(\epsilon')^{1/2}\,, \qquad
  \hbox{for any }\, \epsilon' \in (0,\epsilon)\,,
\end{equation}
where the constant $C$ does not depend on $\epsilon'$. Thus, taking the limit
$\epsilon' \to 0$ in \eqref{uphiaux2}, we arrive at \eqref{unif3} with
$s(t) = s_\epsilon(t)$. 
\end{Proof}

In view of Galilean invariance, Proposition~\ref{p:self-similar} is simply a
reformulation of case (i) in Proposition~\ref{mu-prop2}. In case (ii), the
solution $u(t,x)$ converges to a translate of the viscous shock $\phi_{b,-b}$ if
the shift function $s(t)$ has a finite limit as $t \to -\infty$, or to the
constant state $\pm b$ if $s(t) \to \pm\infty$. A priori it is also possible
that $u(t,x)$ does not converge at all, if $s(t)$ has an oscillatory behavior,
but we do not have any explicit example. In any case $u(t,x)$ cannot
converge to zero in $L^\infty_\loc(\R)$, because this would contradict either
\eqref{unif2} of \eqref{unif3}.  So we see that Propositions~\ref{mu-prop1} and
\ref{mu-prop2} together imply Proposition~\ref{p:suppmu}.

\begin{rem}\label{rem:atom}
It is also possible to detect the presence of atoms in the measure
$\mu$ by using a different scaling in the ancient limit. Assume for instance
that $0 \in \supp(\mu)$. For any $x \in \R$, the representation formula
\eqref{UPK-rep} can be written in the form
\begin{equation}\label{atom1}
  |t|^{1/2} u(t,x|t|^{1/2}) \,=\, \frac{\int |t|^{1/2}z \,e^{-zx|t|^{1/2}/2+z^2t/4}
  \dd\mu(z)}{\int e^{-zx|t|^{1/2}/2+z^2t/4}\dd\mu(z)}\,, \qquad t < 0\,,
  \quad x \in \R\,.
\end{equation}
By Lebesgue's dominated convergence theorem, the numerator in the right-hand side
vanishes in the ancient limit $t \to -\infty$, whereas the denominator converges to 
$\mu(\{0\})$. So, if the measure $\mu$ has an atom at the origin, we deduce that
$|t|^{1/2}u(t,x|t|^{1/2})$ converges to zero in $L^\infty_\loc(\R)$ as $t \to -\infty$.
Now, assume on the contrary that $\dd\mu(z) = f(z)\dd z$ near the origin, where
the density $f$ is continuous and satisfies $f(0) > 0$. Using the change of
variable $z = y|t|^{-1/2}$, we can transform \eqref{atom1} into
\begin{equation}\label{atom2}
  |t|^{1/2} u(t,x|t|^{1/2}) \,=\, -x + \frac{\int (y+x)\,e^{-(y+x)^2/4} f(y|t|^{-1/2})
  \dd y}{\int \,e^{-(y+x)^2/4} f(y|t|^{-1/2})\dd y}\,, \qquad t < 0\,,
  \quad x \in \R\,.
\end{equation}
Applying Lebesgue's theorem again, we see that $|t|^{1/2}u(t,x|t|^{1/2})$ converges
to $-x$ in $L^\infty_\loc(\R)$ as $t \to -\infty$. This explains the observations
made in \eqref{u0u1}. 
\end{rem}

\section{Long-time the asymptotics beyond shocks}\label{s:6}

Equipped with the representation formula \eqref{UPK-rep}, we now return to the
discussion of $\omega$-limit sets, focusing our attention to the particular case
of Burgers' equation. If $u_0 \in L^\infty(\R)$, we know from
Proposition~\ref{omega-prop} that $\omega(u_0)$ is bounded and fully invariant
under the evolution semigroup $(\cS_t)_{t \ge 0}$ defined by
\eqref{Burgers}. This implies that any $\phi \in \omega(u_0)$ is the evaluation
at time $t = 0$ of some bounded entire solution of \eqref{Burgers}. Applying
Proposition~\ref{UPK-prop}, we thus find:

\begin{cor}\label{UPK-cor}
For any $\phi \in \omega(u_0)$, where $\omega(u_0)$ is the $\omega$-limit
set \eqref{e:omega} corresponding to Burgers' equation, there exists a
unique probability measure $\mu$ on $\R$ such that
\begin{equation}\label{UPK-rep2}
  \phi(x) \,=\, \genfrac{}{}{1pt}{0}{\int z\,e^{-zx/2} \dd\mu(z)}{
  \int e^{-zx/2} \dd\mu(z)}\,, \quad x \in \R\,.
\end{equation}
\end{cor}

This result means that, at least for Burgers' equation, the $\omega$-limit set
of any solution of \eqref{Burgers} with values in some interval $[\alpha,\beta]$
can be identified with a subset of all probability measures supported on that
interval. This does not imply, however, that any probability measure on
$[\alpha,\beta]$ can be realized in this way.  To make the discussion more
precise, let us denote
\begin{equation}\label{Sigmadef}
  \Sigma \,:=\, \bigcup_{\beta \ge \alpha}\,\bigcup_{y \in \R}\,\bigl\{\cT_y
  \,\phi_{\beta,\alpha}\bigr\}\,,
\end{equation}
where $\cT_y$ is the translation operator and $\phi_{\beta,\alpha}$ is given
by \eqref{Burgshock}. In other words $\Sigma$ is the collection of all translates
of all viscous shocks, including the constants. Any $\phi \in \Sigma$
corresponds, via \eqref{UPK-rep}, to a probability measure $\mu$ that is a
convex combination of at most two Dirac masses. 

We know from Proposition~\ref{mono-prop} that $\omega(u_0) \subset \Sigma$
whenever $u_0 \in L^\infty(\R)$ is monotonically decreasing. On the other
hand, Oleinik's inequality \eqref{Oleinik} indicates that all solutions
of \eqref{Burgers} are ``eventually decreasing'' when $t \to +\infty$.
Combining these observations, it is rather tempting to conjecture that
$\omega(u_0) \subset \Sigma$ for any $u_0 \in L^\infty(\R)$. Our last result
provides an example that contradicts this hasty conclusion. To make a
precise statement we introduce for any $\gamma > 0$ the function
$\Psi_\gamma : \R \to \R$ defined by
\begin{equation}\label{Psigamdef}
  \Psi_\gamma(x) \,=\, \frac{-2\sinh(x)}{\gamma + \cosh(x)}\,,
  \qquad x \in \R\,.
\end{equation}
Note that $\Psi_\gamma(x)$ is just the evaluation at time $t = \log(1/\gamma)$
of the two-shock solution \eqref{Burgmerger}. 

\begin{prop}\label{complex-prop}
There exist initial data $u_0 \in L^\infty(\R)$ for Burgers' equation such that
\begin{equation}\label{complex-om}
  \omega_0(u_0) \,\supset\, \bigl\{\Psi_\gamma \,;\, \gamma > 0\bigr\}
  \,\cup\, \bigl\{\phi_{\delta,-\delta}\,;\,\delta \in [0,2]\bigr\}\,,
\end{equation}
where $\omega_0(u_0)$ is the $\omega$-limit set \eqref{e:omega0}. 
In particular $\omega_0(u_0) \not\subset \Sigma$. 
\end{prop}

In other words, Proposition~\ref{complex-prop} gives an example of
bounded initial data $u_0$ such that even the ``small'' $\omega$-limit
set $\omega_0(u_0)$ contains the two-shock solution \eqref{Burgmerger},
in addition to a continuum of steady shocks. More generally we conjecture
that, for any probability measure $\mu$ on the interval $[\alpha,\beta]$,
there exist initial data $u_0 \in L^\infty(\R)$ satisfying \eqref{alphabet}
such that the $\omega$-limit set $\omega(u_0)$ contains the function $\phi$
defined by \eqref{UPK-rep2}. We believe that this (nontrivial) extension
of Proposition~\ref{complex-prop} can be obtained following the same
lines of thought as in Section~\ref{s:61} below. This question is
left for future work. 

Examples of $\omega$-limit sets with complicated structure were also
constructed for reaction-diffusion equations on the real line, see e.g.
\cite{ER,Po1,Po2}. In those examples, nonstationary solutions appear in the
$\omega$-limit set as a result of a {\em coarsening dynamics}. The same idea is
exploited here in our proof of Proposition~\ref{complex-prop}, but the result is
in some sense more surprising because it is not clear a priori if something like
a coarsening dynamics is compatible with the constraints imposed by Oleinik's
inequality \eqref{Oleinik}.

\begin{rem}\label{lastrem}
In the spirit of the work of Slijep\v{c}evi\'c and the first author,
one may ask if, for general initial data (including those
considered in Proposition~\ref{complex-prop}), the solution
$u(t) = \cS_t u_0$ approaches locally uniformly the set $\Sigma$
at least for ``almost all times'', in the precise sense
considered in \cite{GS}. We hope to address that interesting
question in a near future.   
\end{rem}

\subsection{Shock mergers in the $\omega$-limit set}\label{s:61}

In this section, we construct bounded initial data for 
Burgers' equation \eqref{Burgers} such that the corresponding solution
exhibits mergers of viscous shocks at the origin for an infinite
sequence of times. The construction is based on the Cole-Hopf
representation formula \eqref{u-CH}. 

\begin{df}\label{Vm-def}
For any $m \ge 0$, let $V_m(t,x)$ be the solution of the
linear heat equation $\partial_t V_m = \partial_x^2 V_m$ with
initial data 
\begin{equation}\label{Vm}
  V_m(0,x) \,=\, \begin{cases} \cosh(x) & \hbox{if }~ |x| \le m\,,\\
  \cosh(m) & \hbox{if }~ |x| > m\,.\end{cases}
\end{equation}
\end{df}

Since $e^t\cosh(x)$ is an exact solution of the heat equation, the
parabolic maximum principle implies that $1 \le V_m(t,x) \le \max\bigl(\cosh(m),
e^t\cosh(x)\bigr)$ for all $t \ge 0$ and all $x \in \R$. The
following two lemmas give more precise estimates on the function
$V_m$ and its derivative. 

\begin{lem}\label{long-lem}
For any $t > 0$ and any $x \in \R$ the following estimates hold\:
\begin{equation}\label{long-ineq}
  \cosh(m)\Bigl(1 - \frac{m}{\sqrt{\pi t}}\Bigr) \,\le\, V_m(t,x) \,\le\,
  \cosh(m)\,, \qquad \bigl|\partial_x V_m(t,x)\bigr| \,\le\, \frac{m\cosh(m)}{
  \sqrt{4\pi}\,t}\,.
\end{equation}
\end{lem}

\begin{Proof}
Let $W_m(t,x) = \cosh(m) - V_m(t,x)$. Then $W_m(t,x)$ satisfies the heat equation
on $\R$, so that
\begin{equation}\label{Wmdef}
  W_m(t,x) \,=\, \frac{1}{\sqrt{4\pi t}}\int_{-m}^m e^{-(x-y)^2/(4t)}
  \bigl(\cosh(m) - \cosh(y)\bigr)\dd y\,, \qquad (t,x) \in \Omega_+\,.
\end{equation}
It is clear from this representation that
\[
  0 \,\le\, W_m(t,x) \,\le\, \frac{1}{\sqrt{4\pi t}}\int_{-m}^m 
  \cosh(m)\dd y \,=\, \frac{m\cosh(m)}{\sqrt{\pi t}}\,,
\]
which gives the first two inequalities in \eqref{long-ineq}. Similarly,
differentiating \eqref{Wmdef}, we find
\[
  |\partial_x W_m(t,x)| \,\le\, \frac{1}{\sqrt{4\pi t}}\int_{-m}^m
  \frac{|x-y|}{2t}\,e^{-(x-y)^2/(4t)} \cosh(m) \dd y \,\le\,
  \frac{m\cosh(m)}{\sqrt{4\pi}\,t}\,,
\]
where we used the fact that $z\,e^{-z^2} \le 1/2$ for all $z \ge 0$. This
concludes the proof of \eqref{long-ineq}. 
\end{Proof}

Estimates \eqref{long-ineq} provide a good approximation of the
solution $V_m(t,x)$ for large times. The short time behavior near
the origin is described by the following result. 

\begin{lem}\label{short-lem}
Assume that $t > 0$ and $|x| + 2t \le m/2$. Then
\begin{equation}\label{short-ineq}
  \begin{split}
  e^t\cosh(x)\Bigl(1 - e^{-m^2/(16t)}\Bigr) \,&\le\, V_m(t,x) \,\le\,
  e^t\cosh(x)\,, \\
  \bigl|\partial_x V_m(t,x) - e^t\sinh(x)\bigr| \,&\le\, e^{-m^2/(16t)}
  \,e^t \cosh(x)\,.
  \end{split}
\end{equation}
\end{lem}

\begin{Proof}
Let $\tilde W_m(t,x) = e^t\cosh(x) - V_m(t,x)$. Then $\tilde W_m(t,x)$ is
again a solution of the heat equation, hence
\begin{equation}\label{tWmdef}
  \tilde W_m(t,x) \,=\, \frac{1}{\sqrt{4\pi t}}\int_{|y| \ge m} e^{-(x-y)^2/(4t)}
  \bigl(\cosh(y) - \cosh(m)\bigr)\dd y\,, \qquad (t,x) \in \Omega_+\,.
\end{equation}
Our main goal is to find an upper bound on $\tilde W_m(t,x)$ in the region region
where $|x| + 2t \le m/2$. We observe that $0 \le \tilde W_m(t,x) \le I_+(t,x)
+ I_-(t,x)$ where
\begin{equation}\label{Ipm}
  \begin{split}
  I_+(t,x) \,&=\, \frac{1}{\sqrt{4\pi t}}\int_m^\infty e^{-(x-y)^2/(4t)}
  \,e^y\dd y \,=\, \frac12\,e^{t+x} \erfc\Bigl(\frac{m{-}x{-}2t}{2\sqrt{t}}
  \Bigr)\,, \\
  I_-(t,x) \,&=\, \frac{1}{\sqrt{4\pi t}}\int_{-\infty}^{-m} e^{-(x-y)^2/(4t)}
  \,e^{-y}\dd y \,=\, \frac12\,e^{t-x} \erfc\Bigl(\frac{m{+}x{-}2t}{2\sqrt{t}}
  \Bigr)\,,
  \end{split}
\end{equation}
where $\erfc$ denotes the complementary error function. In the last equalities
in \eqref{Ipm}, we used the changes of variables $y = x \pm (2t + 2\sqrt{t}z)$
to reduce the integrals to an error function. By assumption, we have
$m \pm x -2t \ge m - |x| - 2t \ge m/2$, and it is known that $\erfc$
is a decreasing function on $\R_+$ which satisfies $\erfc(z) \le
e^{-z^2}$ for all $z \ge 0$. This leads to the upper bound
\[
  \tilde W_m(t,x) \,\le\, e^t \cosh(x)\,\erfc\Bigl(\frac{m}{4\sqrt{t}}\Bigr)
  \,\le\, e^t\cosh(x)\,e^{-m^2/(16t)}\,,
\]
which proves the first part of \eqref{short-ineq}. The second inequality
is established by a similar calculation based on the identity
\[
  \partial_x\tilde W_m(t,x) \,=\, \frac{1}{\sqrt{4\pi t}}\int_{|y| \ge m}
  e^{-(x-y)^2/(4t)} \sinh(y)\dd y\,, \qquad (t,x) \in \Omega_+\,.
\]
This concludes the proof of Lemma~\ref{short-lem}. 
\end{Proof}

We now explain our strategy to prove Proposition~\ref{complex-prop}.
If $m > 0$ is large enough, the solution of Burgers' equation given
by $u(t,x) = -2\partial_x V_m(t,x)/(1+V_m(t,x))$ satisfies,
by Lemma~\ref{short-lem}, 
\[
  u(t,x) \,\approx\, \frac{-2\sinh(x)}{e^{-t} + \cosh(x)}\,,
  \qquad \hbox{when}\quad |x| + 2t \le m/2\,,
\]
whereas $u(t,x) \approx 0$ when $t \ge m^2$ by Lemma~\ref{long-lem}. 
In other words $u(t,x)$ describes, for relatively small times, the merger of
a pair of viscous shocks at the origin, but in the long time regime
$u(t,x)$ actually converges to zero, uniformly in $x$ on compact intervals. 
The idea is to construct, by superposition, a solution of \eqref{Burgers}
which exhibits infinitely many such mergers, along an appropriate sequence
of times.

\medskip\noindent{\bf Proof of Proposition~\ref{complex-prop} (first part).}
Fix $N \ge 10$ and let $1 \le t_1 < t_2 < t_3 < \dots$ be a
sequence of times such that
\begin{equation}\label{tjdef}
  t_{j+1} \,\ge\, N^2 t_j^2 \quad \hbox{for all } j \ge 1\,, \qquad
  \hbox{and}\quad \frac{t_j^2}{t_{j+1}} \,\xrightarrow[j \to +\infty]{}\, 0\,. 
\end{equation}
We consider the function $U : [0,+\infty)\times\R \to (0,+\infty)$ defined by
\begin{equation}\label{Ucomp}
  U(t,x) \,=\, 1 + \sum_{j=1}^\infty e^{-t_j} \,V_{Nt_j}(t,x)\,, \qquad
  t \ge 0\,, \quad x \in \R\,,
\end{equation}
where $V_{Nt_j}$ is given by Definition~\ref{Vm-def} with $m = Nt_j$. 
Since $V_{Nt_j}(t,x) \le e^t\cosh(x)$, it is clear that the series in \eqref{Ucomp}
converge uniformly on compact sets in space-time, and that $U(t,x)$ is a
solution of the heat equation on $\R_+ \times \R$. We are interested in the
corresponding solution $u(t,x)$ of Burgers' equation, given by
\begin{equation}\label{u-ch2}
  u(t,x) \,=\, \frac{-2\partial_x U(t,x)}{U(t,x)}\,, \qquad
  t \ge 0\,,\quad x \in \R\,.
\end{equation}
As $|\partial_x V_m(0,x)| \le V_m(0,x)$ for all $x \in \R$, we have $|\partial_x
V_m(t,x)| \le V_m(t,x)$ for all $t \ge 0$ by the maximum principle, and it
follows that $|u(t,x)| \le 2$, so that $u(t,x)$ is a bounded solution
of Burgers' equation. We shall show that, for any sufficiently large $k \in \N$,
this solution exhibits a merger of viscous shocks at the origin on the time
interval $[t_k,2t_k]$. More precisely, we shall prove that
\begin{equation}\label{umerge}
  \sup_{|x| \le t_k}\,\Bigl|u(\tau_k,x) + \frac{2\sinh(x)}{1 + \cosh(x)}\Bigr|
  \,\xrightarrow[k \to +\infty]{}\, 0\,, \qquad \hbox{where}\quad
  \tau_k \,=\, t_k + (N{-}1)t_{k-1}\,.
\end{equation}

To establish \eqref{umerge}, we fix a large $k \in \N$, and we assume that
$t \in [t_k,2t_k$] and $|x| \le t_k$. If $j \ge k$, we know from
Lemma~\ref{short-lem} that
\begin{equation}\label{sum+}
  e^{t-t_j}\cosh(x)\Bigl(1 - e^{-N^2t_j^2/(16t)}\Bigr) \,\le\,
  e^{-t_j}\,V_{Nt_j}(t,x) \,\le\, e^{t-t_j}\cosh(x)\,, 
\end{equation}
because $|x| + 2t \le 5t_k \le Nt_j/2$ since $N \ge 10$. In view of
\eqref{tjdef}, estimate \eqref{sum+} shows that, in the space-time region
under consideration, all terms with $j \ge k+1$ can be neglected in the sum
\eqref{Ucomp} defining $U(t,x)$, namely
\[
  \sum_{j=k+1}^\infty e^{-t_j} \,V_{Nt_j}(t,x) \,\ll\, e^{-t_k} \,V_{Nt_k}(t,x)\,,
  \qquad \hbox{uniformly for }\,t \in [t_k,2t_k]\,, ~|x| \le t_k\,.
\]
If $j < k$, we apply Lemma~\ref{long-lem} and obtain the bound
\begin{equation}\label{sum-}
  e^{-t_j}\cosh(Nt_j)\Bigl(1 - \frac{Nt_j}{\sqrt{\pi t}}\Bigr) \,\le\,
  e^{-t_j}\,V_{Nt_j}(t,x) \,\le\, e^{-t_j}\cosh(Nt_j)\,,
\end{equation}
where the left-hand side is strictly positive since $t \ge t_k \ge N^2 t_j^2$.
Again it follows from \eqref{tjdef}, \eqref{sum-} that the terms with
$j \le k-2$ can be neglected in the sum \eqref{Ucomp}, in the sense that
\[
  1 + \sum_{j=1}^{k-2} e^{-t_j} \,V_{Nt_j}(t,x) \,\ll\, e^{-t_{k-1}} \,V_{Nt_{k-1}}(t,x)\,,
  \qquad \hbox{uniformly for }\,t \in [t_k,2t_k]\,, ~|x| \le t_k\,.
\]

Summarizing, we have shown that the function $U(t,x)$ defined by \eqref{Ucomp}
satisfies $U(t,x) \approx U_k(t,x)$ when $t \in [t_k,2t_k$] and $|x| \le t_k$,
where
\begin{equation}\label{Ukdef}
  U_k(t,x) \,:=\, e^{-t_{k-1}} \,V_{Nt_{k-1}}(t,x) + e^{-t_k} \,V_{Nt_k}(t,x)
  \,\approx\, e^{(N-1)t_{k-1}} + e^{t-t_k}\cosh(x)\,.
\end{equation}
More precisely we have
\begin{equation}\label{Uapprox}
  \sup_{t \in [t_k,2t_k]} \sup_{|x| \le t_k}\biggl(\frac{|U(t,x) - U_k(t,x)|}{
  U_k(t,x)} + \frac{|\partial_x U(t,x) - \partial_x U_k(t,x)|}{
  U_k(t,x)}\biggr) \,\xrightarrow[k \to +\infty]{}\, 0\,,
\end{equation}
where the estimate for the derivative is obtained in the same way 
using Lemmas~\ref{long-lem} and \ref{short-lem}.  We now take $t = \tau_k :=
t_k + (N{-}1)t_{k-1}$, so that both terms in the right-hand side
of \eqref{Ukdef} are of comparable size. With that choice, we have
$U_k(\tau_k,x) \approx e^{\tau_k - t_k}\bigl(1 + \cosh(x)\bigr)$, and using
estimates \eqref{Uapprox} we easily deduce that the function
$u(\tau_k,x) = -2\partial_x U(\tau_k,x)/U(\tau_k,x)$ indeed satisfies
\eqref{umerge}.

So far we have shown that the function $\Psi_\gamma$ introduced in
\eqref{Psigamdef} belongs to the $\omega$-limit set $\omega_0(u_0)$
when $\gamma = 1$. But the argument above also implies that, given any $t \in \R$, 
\[
  \lim_{k\to +\infty} \sup_{|x| \le t_k} \,\Bigl|u(t+\tau_k,x) +
  \frac{2\sinh(x)}{e^{-t} + \cosh(x)}\Bigr| \,=\, 0\,,
\]
which shows that the $\omega$-limit set $\omega_0(u_0)$ contains the entire
two-shock solution \eqref{Burgmerger} (as is clear from time invariance). 
So we conclude that $\omega_0(u_0) \supset \bigl\{\Psi_\gamma \,;\, \gamma > 0\bigr\}$,
as asserted in \eqref{complex-om}. \QED

\subsection{Repair along the family of steady shocks}\label{s:62}

In the topology of $L^\infty_\loc(\R)$, the two-shock solution
\eqref{Burgmerger} converges to zero as $t \to -\infty$ and to the steady shock
$\phi_{2,-2}$ as $t \to +\infty$. Such a heteroclinic connection is obviously not
chain recurrent in the sense of Proposition~\ref{omega-prop0}. As a consequence,
for the initial data $u_0 \in L^\infty(\R)$ constructed in the previous section,
the $\omega$-limit set must be larger than the heteroclinic orbit given by the
two-shock solution. In this section, we show that $\omega_0(u_0)$ contains in
addition a continuum of steady shocks, as stated in Proposition~\ref{complex-prop}.

To prove the claim, we need to control the function $V_m(t,x)$ introduced
in Definition~\ref{Vm-def} for some intermediate times that are not covered
by Lemmas~\ref{long-lem} and \ref{short-lem}.

\begin{lem}\label{inter-lem}
For any fixed $\delta \in (0,2)$ we have, as $m \to +\infty$,
\begin{equation}\label{inter-exp}
  V_m(m/\delta,x) \,=\, \frac{1}{\sqrt{\pi m \delta}}\,
  \frac{2+\delta}{2-\delta}\, e^{m(1-\delta/4)}\,\cosh(\delta x/2)
  \Bigl(1 + \cO(m^{-1})\Bigr)\,,
\end{equation}
where convergence is uniform in $x \in [-L,L]$ for any $L > 0$.
A similar asymptotic expansion also holds for the spatial derivative
$\partial_xV_m(m/\delta,x)$. 
\end{lem}

\begin{Proof}
To establish \eqref{inter-exp} it is convenient to use an explicit
expression for the function $V_m(t,x)$. Starting from the definition
\[
  V_m(t,x) \,=\, \frac{1}{\sqrt{4\pi t}}\int_{|y| \le m} e^{-|x-y|^2/(4t)}
  \cosh(y)\dd y + \frac{1}{\sqrt{4\pi t}}\int_{|y| > m} e^{-|x-y|^2/(4t)}
  \cosh(m)\dd y\,,
\]
and proceeding as in \eqref{Ipm}, we obtain the tractable formula
\begin{equation}\label{Vmdecomp}
  V_m(t,x) \,=\, \hat V_m(t,x) + \hat V_m(t,-x) + \hat W_m(t,x) +
  \hat W_m(t,-x)\,, \qquad t > 0\,, \quad x \in \R\,,
\end{equation}
where  
\begin{equation}\label{Vmexp}
  \hat V_m(t,x) \,=\, \frac{\cosh(m)}{2}\,\erfc\Bigl(\frac{m{+}x}{2\sqrt{t}}
  \Bigr)\,, \quad 
  \hat W_m(t,x) \,=\, \frac{e^{t+x}}{2}\biggl\{\erfc\Bigl(\frac{2t{-}m{+}x}{
  2\sqrt{t}}\Bigr) - \erfc\Bigl(\frac{2t{+}m{+}x}{2\sqrt{t}}\Bigr)\biggr\}\,.
\end{equation}
In what follows, we assume that $t = m/\delta$ for some fixed $\delta \in
(0,2)$, and we consider the limit $m \to +\infty$ for $x$ in some
fixed interval $[-L,L]$. Using the asymptotic expansion of the complementary
error function
\[
  \erfc(z) \,=\, \frac{e^{-z^2}}{z\sqrt{\pi}}\,\Bigl(1 + \cO(z^{-2})\Bigr)\,,
  \qquad \hbox{as } z \to +\infty\,,
\]
we easily find that
\begin{equation}\label{hatVm}
  \hat V_m(m/\delta,x) \,=\, \frac12\,\frac{1}{\sqrt{\pi m \delta}}\,
  e^{m(1-\delta/4)}\,e^{-\delta x/2} \Bigl(1 + \cO(m^{-1})\Bigr)\,,
  \qquad m \to +\infty\,.
\end{equation}
Next, observing that $2t+m \gg 2t-m \gg 1$ as $m \to +\infty$ because
$2/\delta > 1$, we see that the second term in the expression
\eqref{Vmexp} of $\hat W_m(t,x)$ is negligible compared to the
first one. Since
\[
  t+x - \biggl(\frac{2t{-}m{+}x}{2\sqrt{t}}\biggr)^2 \,=\, m -
  \frac{(m{-}x)^2}{4t} \,=\, m\bigl(1-\delta/4\bigr) + \delta x/2
  + \cO(m^{-1})\,,
\]
we thus find
\begin{equation}\label{hatWm}
  \hat W_m(m/\delta,x) \,=\, \frac{1}{\sqrt{\pi m \delta}}\,
  \frac{\delta}{2-\delta}\,e^{m(1-\delta/4)}\,e^{\delta x/2}
  \Bigl(1 + \cO(m^{-1})\Bigr)\,,   \qquad m \to +\infty\,. 
\end{equation}
Finally, replacing \eqref{hatVm} and \eqref{hatWm} into \eqref{Vmdecomp},
we arrive at \eqref{inter-exp}. A corresponding asymptotic expansion
for the derivative $\partial_x V_m(m/\delta,x)$ is obtained in the
same way. 
\end{Proof}

\medskip\noindent{\bf Proof of Proposition~\ref{complex-prop} (second part).}
We consider again the solution $U(t,x)$ of the heat equation given
by \eqref{Ucomp}, and we evaluate it along the sequence of times
$\hat\tau_k := Nt_k/\delta$, for some fixed $\delta \in (0,2)$. The
main contribution to the sum comes from the term where $j = k$, and using
Lemma~\ref{inter-lem} with $m = Nt_k$ we obtain
\begin{equation}\label{jeqk}
  e^{-t_k}\,V_{Nt_k}(Nt_k/\delta,x) \,=\, \frac{e^{-t_k}}{\sqrt{\pi N t_k
  \delta}}\, \frac{2+\delta}{2-\delta}\, e^{N t_k(1-\delta/4)}\,\cosh(\delta x/2)
  \Bigl(1 + \cO(t_k^{-1})\Bigr)\,,
\end{equation}
where convergence holds uniformly for $x \in [-L,L]$. The terms where
$j < k$ can be easily estimated using the trivial bound $V_m(t,x) \le
\cosh(m)$, leading to
\begin{equation}\label{jltk}
  1 + \sum_{j=1}^{k-1} e^{-t_j} \,V_{Nt_j}(Nt_k/\delta,x) \,\le\,
  1 + \sum_{j=1}^{k-1} e^{-t_j} \cosh(Nt_j) \,\le\, C\,e^{(N-1)t_{k-1}}\,, 
\end{equation}
where the right-hand side is much smaller than \eqref{jeqk} since
$t_k \ge N^2 t_{k-1}^2$. Finally, for the terms where $j > k$, we
use the simple bound $V_m(t,x) \le e^t \cosh(x)$ which gives
\begin{equation}\label{jgtk}
  \sum_{j=k+1}^\infty e^{-t_j} \,V_{Nt_j}(Nt_k/\delta,x) \,\le\,
  \sum_{j=k+1}^\infty e^{-t_j}\,e^{Nt_k/\delta}\,\cosh(x)\,. 
\end{equation}
Again, for $x \in [-L,L]$, the right-hand side is much smaller than
\eqref{jeqk} because the sequence $t_j$ grows fast enough as $j \to
\infty$ and $t_j \ge t_{k+1} \ge N^2 t_k^2$.

Summarizing, it follows from \eqref{jeqk}, \eqref{jltk}, \eqref{jgtk}
that the function $U(t,x)$ defined by \eqref{Ucomp} satisfies
\[
  U(Nt_k/\delta,x) \,=\, \frac{e^{-t_k}}{\sqrt{\pi N t_k\delta}}\,
  \frac{2+\delta}{2-\delta}\, e^{N t_k(1-\delta/4)}\,\cosh(\delta x/2)
  \Bigl(1 + \cO(t_k^{-1})\Bigr)\,, \quad k \to +\infty\,,
\]
uniformly for $x \in [-L,L]$, and a similar expansion also holds for
the first derivative $\partial_x U(t,x)$. So, we deduce that the
solution of Burgers' equation defined by \eqref{Ucomp}, \eqref{u-ch2} 
satisfies, for any $L > 0$, 
\[
  \sup_{|x| \le L}\,\bigl|u(Nt_k/\delta,x) + \delta\tanh(\delta x/2)\bigr| 
  \,\xrightarrow[k \to +\infty]{}\, 0\,.
\]
This implies that the $\omega$-limit set $\omega_0(u_0)$ contains
the viscous shock $\phi_{\delta,-\delta}$ for any value of $\delta \in (0,2)$, 
hence for any $\delta \in [0,2]$ since $\omega_0(u_0)$ is closed in 
$L^\infty_\loc(\R)$. The proof of \eqref{complex-om} is now complete.  
\QED

\begin{figure}[ht!]
  \begin{center}
  \begin{picture}(420,190)
  \put(30,0){\includegraphics[width=1.00\textwidth]{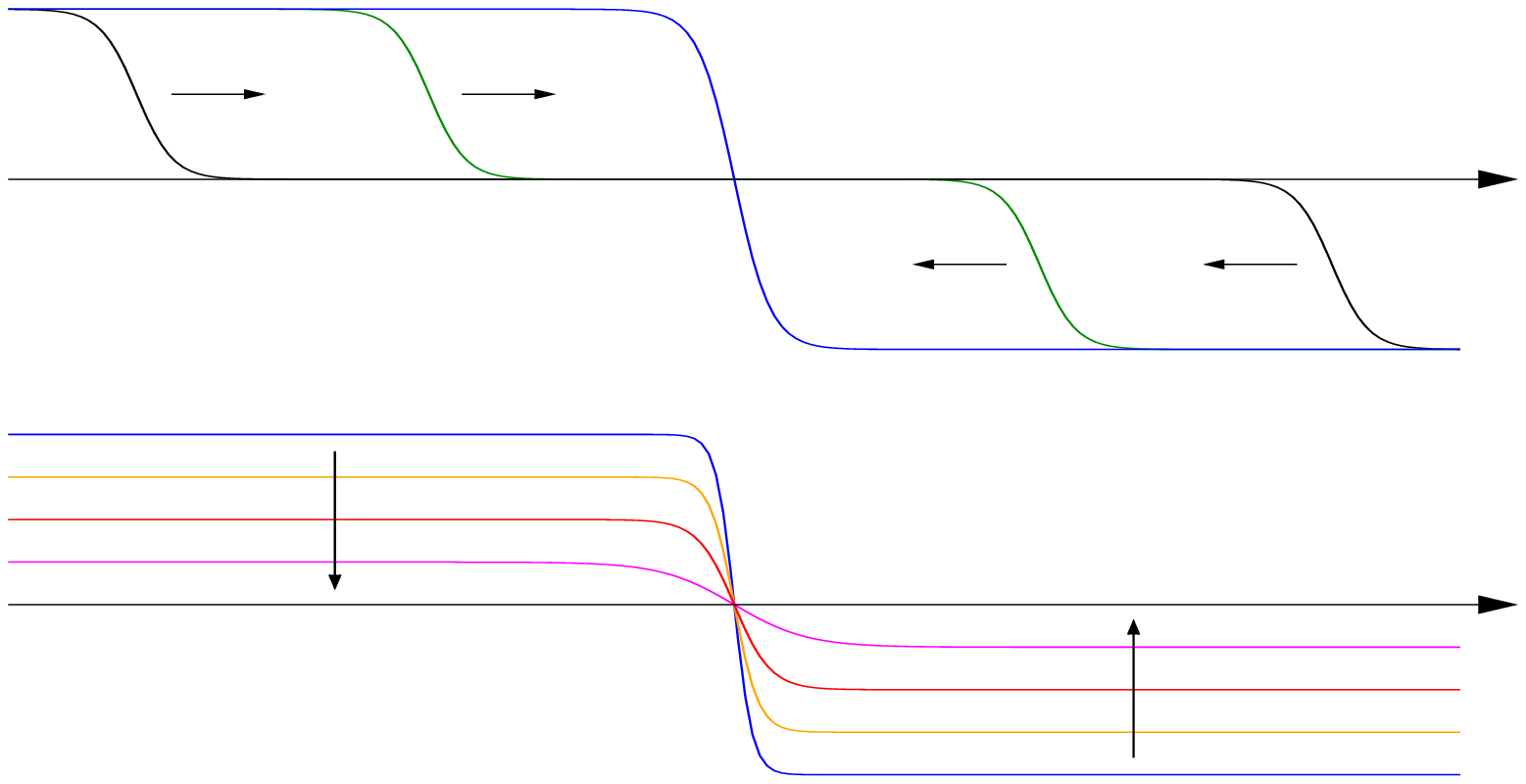}}
  \put(358,144){$x$}
  \put(358,51){$x$}
  \put(185,06){$-2$}
  \put(31,44){$0$}
  \put(31,82){$2$}
  \put(190,101){$-2$}
  \put(31,137){$0$}
  \put(31,173){$2$}
  \put(210,165){(a) Merger of two viscous shocks}
  \put(210,75){(b) Repair along the family of shocks}
  \end{picture}
  \caption{{\small Illustration of the long-time behavior of the solution 
  $u(t,x)$ of \eqref{Burgers} given by \eqref{Ucomp}, \eqref{u-ch2}. (a) Along 
  a sparse sequence of times $\tau_k \to +\infty$, the solution describes 
  the merger of a pair of viscous shocks near the origin, as in the explicit solution
  \eqref{Burgmerger}. (b) Between the times $\tau_k$ and $\tau_{k+1}$, the 
  solution slowly returns to zero along the family of steady shocks 
  $\phi_{\delta,-\delta}$, where $0 < \delta < 2$. Both processes recur 
  infinitely often, and are therefore reflected in the $\omega$-limit set
  of the solution $u(t,x)$, as asserted in Proposition~\ref{complex-prop}.}}
  \label{f:fig2}
  \end{center}
\end{figure}

\begin{rem}\label{conj-rem}
We conjecture that, for the initial data $u_0$ constructed in the proof
of Proposition~\ref{complex-prop}, the $\omega$-limit set $\omega_0(u_0)$
is in fact equal to the right-hand side of \eqref{complex-om}. Note
that this set satisfies all the properties listed in
Proposition~\ref{omega-prop0}, including chain recurrence.
\end{rem}

\section{Discussion}\label{s:7}

We presented results on long-time behavior in scalar conservation laws on
the real line, both in the case of a general convex flux and in the special case
of Burgers' equation with quadratic flux. Our main results include a general
definition and characterization of $\omega$-limit sets, the convergence to
single shocks for monotone data, and the construction of initial data for which
the $\omega$-limit set does not consist of a constant state nor of the translates
of a single shock. The latter result was established in the context of
Burgers' equation, where a somewhat explicit representation of all entire
solutions in terms of compactly supported probability measures is
available. Since all elements in the $\omega$-limit set are entire solutions,
this characterization provides a "list" of candidates for elements in the
$\omega$-limit set.

We mentioned throughout several open problems and conjectured some answers. We
revisit some of those here in a broader context and point to some other
potentially interesting questions.

As a first step towards a complete characterization, one can ask what functions
may be found within an $\omega$-limit set. Candidates are entire solutions,
which, in the case of Burgers' equation, are described somewhat explicitly
through a one-to-one correspondence with compactly supported probability
measures on the real line. Even beyond the goal of describing the long-time
behavior of general solutions, it would be quite interesting to characterize
entire solutions of scalar conservation laws with convex, not necessarily
quadratic flux. In the absence of a direct connection with the heat equation,
we think that a description of bounded solutions in terms of their
ancient limits $t\to-\infty$ in suitably rescaled variables might provide
an avenue for progress in this direction. Note that solutions representing the
superposition of two viscous shocks as $t \to -\infty$ can be constructed
under generic assumptions on the flux function, see \cite{Se3}. 

On the other hand, it would be interesting to extend the analysis of ancient
solutions and possibly the characterization of $\omega$-limit sets to the
complex-valued Burgers' equation, where the Cole-Hopf transformation is still
at hand, but $L^\infty$ upper bounds, Oleinik's inequality, and the positivity
that is essential in the characterization of ancient solutions are not available;
see \cite{PoSv} for results on blowup in this context.

Given the characterization of entire solutions in Burgers' equation, we
conjectured in Section~\ref{s:6} that any entire solution of that equation can
be found in the $\omega$-limit set for appropriate initial data. Beyond Burgers'
equation, one may find it plausible that the existence of shock mergers in
specific $\omega$-limit sets can be established, by controlling the interaction of
shocks and rarefaction waves without the conjugation to a linear heat equation
and the associated superposition principle.

A more ambitious result would characterize the entire $\omega$-limit set. We
showed that for monotone initial data, one only finds a single shock (together
with its translates) or a family of constant states. We conjectured that in the
example considered in Proposition~\ref{complex-prop}, the $\omega$-limit set
actually consists of the shock merger itself and of the family of steady shocks
with smaller amplitude, which together form a chain recurrent set. Given the
ancient asymptotics of entire solutions, all nonconstant elements of the
$\omega$-limit sets can be thought of as continuous or discrete superpositions
of shocks and their eventual merger into a single shock. The natural question in
this direction is whether different shock mergers can occur within a single
$\omega$-limit set. Eventually, one may hope to determine which subsets of the
set of ancient solutions may occur as $\omega$-limit sets, that is, to decide if
any additional restrictions beyond compactness, connectedness, and chain
recurrence are imposed by the dynamics.

Clearly, all of the questions above can be asked for $\omega_0(u_0)$ and for
$\omega(u_0)$, that is in a fixed frame of reference or up to translations. Our
introduction of $\omega(u_0)$, while seemingly natural, can of course be
questioned. One could ask for a more narrow characterization, limiting the
allowed translations for instance to almost Galilean shifts as suggested in
Remark~\ref{other-def}. To clarify the role of the shifts, it would be
interesting to identify $\omega$-limit sets that actually depend on the class of
allowed spatial translates. More specifically, one can ask if for all $u_0$, the
set $\omega(u_0)$ defined by \eqref{e:omega} coincides with the $\omega$-limit
set obtained by restricting the class of allowed shifts to almost Galilean ones.

Beyond the structure of $\omega$-limit sets, one can investigate the dynamics
for large but finite times. While Burgers' equation is not a gradient flow, our
results basically show that the long-term asymptotics of solutions are to a
large extent determined by equilibria, up to Galilean boosts. While we did show
that solutions other than equilibria, particularly shock mergers, can occur in
the $\omega$-limit set, it is conceivable to conjecture that those occur only
"rarely" in time, a statement that one could attempt to quantify in the spirit
of the work of S. Slijep\v{c}evi\'c and the first author \cite{GS}; see Remark
\ref{lastrem}.

Finally, a number of subtle questions arise when attempting to characterize the
set of initial data that lead to a specific $\omega$-limit set. From local
stability of viscous shocks, one can conclude that the basin of attraction is
open (in appropriate topologies). On the other hand, we showed that the basin
contains all monotone initial data with the same limits at $x=\pm\infty$. The
construction of repeating shock mergers suggests robustness of this asymptotic
behavior at least in a spatially uniform topology. We note however that such
questions on the basin of attraction of an $\omega$-limit set do not appear to
have been answered in the case of mergers between layers in the Allen-Cahn
equation \cite{ER,Po1,Po2}. Despite the apparent similarities between the
results there and our construction, it is worth noticing that in our case, all
equilibria and traveling waves are asymptotically stable, whereas the Allen-Cahn
equation accommodates a large family of unstable equilibria and traveling waves,
including the zero solution, spatially periodic equilibria, and traveling waves
connecting those equilibria to stable solutions; see for instance \cite{as}. We
are not aware of results that connect the role of these unstable equilibria to
the description of long-term dynamics through $\omega$-limit sets as attempted
here and in \cite{ER,Po1,Po2}.

\appendix

\section{Appendix}\label{secA}

\subsection{Oleinik's inequality}\label{ssecA1}

If $u(t,x)$ is a solution of \eqref{SCL} with initial data $u_0$
satisfying \eqref{alphabet}, we define
\begin{equation}\label{vauxdef}
  v(t,x) \,=\, t^2 \partial_x u(t,x) - k^{-1}t\,, \qquad
  t > 0\,, \quad x \in \R\,, 
\end{equation}
where $k > 0$ is defined in \eqref{Oleinik}. The function $v$ is smooth, and it
is clear by construction that $v(t,x) < 0$ for all $x \in \R$ whenever $t > 0$ is
sufficiently small. Indeed, since $u(t,x)$ solves equation \eqref{SCL} with bounded
initial data, we know that there exist positive constants $C$ and $t_0$
such that $|\partial_x u(t,x)| \le C t^{-1/2}$ when $0 < t < t_0$. Now a direct
calculation shows that $v$ solves the equation
\begin{equation}\label{vauxeq}
  \begin{split}
  \partial_t v + f'(u)\partial_x v - \partial_x^2 v \,&=\,  2t \partial_x u
  - k^{-1} - t^2 f''(u)\bigl(\partial_x u\bigr)^2\\ \,&\le\,  2t \partial_x u
  - k^{-1} - t^2 k \bigl(\partial_x u\bigr)^2 = -k^{-1}\bigl(1-kt\partial_xu\bigr)^2
  \,\le\, 0\,,
  \end{split}
\end{equation}
where in the second line we used the fact that $f''(u) \ge k$ for all
$u \in [\alpha,\beta]$. By the maximum principle, the differential
inequality \eqref{vauxeq} implies that $v(t,x)$ stays negative for
all times $t > 0$, which gives inequality \eqref{Oleinik}. 

\subsection{Proof of Lemma~\ref{contlem}}\label{ssecA2}

We prove here that the solution of $u(t) = \cS_t(u_0)$ of \eqref{SCL} depends
continuously on the initial data $u_0$ in the topology of $L^\infty_\loc(\R)$.  
Our starting point is the integral equation associated with \eqref{SCL}, namely
\begin{equation}\label{SCLint}
  u(t) \,=\, K(t,\cdot)*u_0 - \int_0^t \partial_x K(t-s,\cdot)*f(u(s))\dd s\,, \qquad
  t > 0\,,
\end{equation}
where $K(t,x)$ is the heat kernel \eqref{K-def} and $*$ denotes convolution
in space. Straightforward calculations show that there exists a constant $C > 0$
such that
\begin{equation}\label{Kexpbound}
  \|K(t,\cdot)*u_0\|_{\exp} + \sqrt{t}\,\|\partial_x K(t,\cdot)*u_0\|_{\exp}
  \,\le\, C e^t \|u_0\|_{\exp}\,,
\end{equation}
for all $u_0 \in L^\infty(\R)$ and all $t > 0$. To prove the desired continuity
property, we fix $R > 0$ and consider two sets of initial data $u_0, v_0$ such
that $\max(\|u_0\|_{L^\infty}, \|v_0\|_{L^\infty}) \le R$. Denoting $u(t) = \cS_t(u_0)$,
$v(t) = \cS_t(v_0)$, and using \eqref{Kexpbound}, we can estimate
\[
  \|u(t) - v(t)\|_{\exp} \,\le\, C e^t \|u_0 - v_0\|_{\exp} +
  \int_0^t \frac{Ce^{t-s}}{\sqrt{t-s}}\,L \|u(s) - v(s)\|_{\exp}\dd s\,,
\]
where $L = \sup\{|f'(u)|\,;\, |u|\le R\}$. The quantity $\delta(t) :=
e^{-t} \|u(t) - v(t)\|_{\exp}$ thus satisfies an integral inequality that
can be solved using a variant of Gr\"onwall's lemma, see \cite[Lemma~7.1.1]{He}.
This gives an estimate of the form $\|\cS_t(u_0) - \cS_t(v_0)\|_{\exp} \le C_1 e^{C_2t}
\|u_0 - v_0\|_{\exp}$, for some universal constants $C_1, C_2$, which shows that
the solution of \eqref{SCLint} depends continuously on the initial data in the
topology of $L^\infty_\loc$, uniformly in time on compact intervals. 

\subsection{Proof of the $L^1$--$L^\infty$ estimate \eqref{w-est}}\label{ssecA3}

Assume that $\chi : \R \to \R$ is a smooth convex function such that
$\chi(0) = 0$ and $\chi(s) > 0$ for all $s \neq 0$. If $w(t,x)$ is a solution
of \eqref{w-eq} with initial data $w_0 \in L^\infty(\R)$, we compute
\begin{equation}\label{west1}
  \begin{split}
  \frac{\D}{\D t}\int_\R \chi(w)\dd x \,&=\, \int_\R \chi'(w)\bigl(
  \partial_x^2 w - f'(u)\partial_x w\bigr)\dd x \\ \,&=\,
  -\int_\R \chi''(w)\bigl(\partial_x w\bigr)^2\dd x
  + \int_\R f''(u)(\partial_x u)\,\chi(w)\dd x \\ \,&\le\,
  -\int_\R \chi''(w)\bigl(\partial_x w\bigr)^2\dd x \,\le\, 0\,, \qquad t > 0\,,
  \end{split}
\end{equation}
where we used the crucial observation that $f''(u)(\partial_x u) \le 0$, because
$f$ is convex and $u$ is decreasing. As a first application, we take
$\chi_\epsilon(w) = (\epsilon^2 + w^2)^{1/2} - \epsilon$, where $\epsilon > 0$ is
a small parameter. Using \eqref{west1} we easily obtain $\int_\R \chi_\epsilon\bigl(w(t,x)
\bigr)\dd x \le \int_\R \chi_\epsilon\bigl(w_0(x)\bigr)\dd x \le \|w_0\|_{L^1(\R)}$,
for any $t > 0$. Then, invoking Lebesgue's monotone convergence theorem, we can take
the limit $\epsilon \to 0$ and arrive at
\begin{equation}\label{west2}
  \|w(t,\cdot)\|_{L^1(\R)} \,\le\, \|w_0\|_{L^1(\R)}\,, \qquad t > 0\,.
\end{equation}

In a second step, we choose $\chi(w) = w^2$ in \eqref{west1} and we use the
celebrated Nash inequality
\begin{equation}\label{west3}
  \|w(t,\cdot)\|_{L^2(\R)}^3 \,\le\, C_N \|w_x(t,\cdot)\|_{L^2(\R)}
  \|w(t,\cdot)\|_{L^1(\R)}^2\,,
\end{equation}
see \cite{CL}. Taking \eqref{west2} into account and assuming $\|w_0\|_{L^1(\R)} > 0$, 
we obtain the differential inequality
\[
  \frac{\D}{\D t}\,\|w(t,\cdot)\|_{L^2(\R)}^2 \,\le\, -2 \|w_x(t,\cdot)\|_{L^2(\R)}^2
  \,\le\, -\frac{2\|w(t,\cdot)\|_{L^2(\R)}^6}{C_N^2 \|w(t,\cdot)\|_{L^1(\R)}^4} \,\le\,
  -\frac{2\|w(t,\cdot)\|_{L^2(\R)}^6}{C_N^2 \|w_0\|_{L^1(\R)}^4}\,,
\]
which can be integrated to give the $L^1$--$L^2$ estimate
\begin{equation}\label{west4}
  \|w(t,\cdot)\|_{L^2(\R)} \,\le\, C t^{-1/4} \|w_0\|_{L^1(\R)}\,,
  \qquad t > 0\,,
\end{equation}
where $C = (C_N/2)^{1/2}$.

Finally, to estimate the $L^\infty$ norm of $w$, we can either bound 
the $L^{2p}$ norm for all integers $p \ge 2$, or use a duality argument, see
\cite{FaS}. We follow here the latter approach and consider the dual equation
\begin{equation}\label{psi-eq}
  \partial_t \psi(t,x) - \partial_x\bigl(f'(u)\psi\bigr) \,=\,
  \partial_x^2 \psi(t,x)\,,   \qquad t > 0\,,\quad x \in \R\,,
\end{equation}
which has similar properties as \eqref{w-eq}. In particular, proceeding
as in \eqref{west1}, we find
\begin{equation}\label{west5}
  \begin{split}
  \frac{\D}{\D t}\int_\R \chi(\psi)\dd x \,&=\, -\int_\R \chi''(\psi)\bigl(\partial_x
  \psi\bigr)^2\dd x + \int_\R f''(u)(\partial_x u)\,\bigl(\psi\chi'(\psi) - \chi(\psi)\bigr)
  \dd x\,, \\
  \,&\le\, -\int_\R \chi''(\psi)\bigl(\partial_x\psi\bigr)^2\dd x \,\le\, 0\,,
  \end{split}
\end{equation}
because $\chi(\psi) \le \psi\chi'(\psi)$ by convexity. We deduce that estimates
\eqref{west2}, \eqref{west4} also hold for the solutions of \eqref{psi-eq}.
Now, if $w$ solves \eqref{w-eq} with initial data $w_0 \in L^2(\R)$ and $\psi$
solves \eqref{psi-eq} with initial data $\psi_0 \in L^1(\R)$, then for any
$t > 0$ the quantity $\int_\R \psi(t-s,x)w(s,x)\dd x$ is independent of $s \in [0,t]$,
as can be easily verified by differentiation. It follows that
\begin{equation}\label{west6}
  \begin{split}
  \biggl|\int_\R \psi_0(x)w(t,x)\dd x\biggr| \,=\, \biggl|\int_\R \psi(t,x)w_0(x)\dd x
  \biggr| \,&\le\, \|\psi(\cdot,t)\|_{L^2(\R)}\|w_0\|_{L^2(\R)} \\ \,&\le\, C t^{-1/4}
  \|\psi_0\|_{L^1(\R)}\|w_0\|_{L^2(\R)}\,,
  \end{split}
\end{equation}
where the last inequality follows from \eqref{west4}. Clearly \eqref{west6}
is equivalent to the $L^2$--$L^\infty$ estimate
\begin{equation}\label{west7}
  \|w(t,\cdot)\|_{L^\infty(\R)} \,\le\, C t^{-1/4} \|w_0\|_{L^2(\R)}\,,
  \qquad t > 0\,, 
\end{equation}
and the $L^1$--$L^\infty$ bound in \eqref{w-est} follows immediately by
combining \eqref{west4}, \eqref{west7}. 

\subsection{Proof of Lemmas~\ref{Seps-lem} and \ref{seps-lem}}\label{ssecA4}

\begin{proof}[\bf Proof of Lemma~\ref{Seps-lem}.]
Since $\nu_\pm$ are positive measures supported on the interval $[0,\epsilon]$, it
is clear that the functions $J^\pm_\epsilon$ introduced in \eqref{IJaux} satisfy the
estimates
\[
  |\partial_x J^\pm_\epsilon(t,x)| \,\le\, \frac{\epsilon}{2}\, J^\pm_\epsilon(t,x)\,, \qquad
  |\partial_t J^\pm_\epsilon(t,x)| \,\le\, \Bigl(\frac{b\epsilon}{2} + \frac{\epsilon^2}{4}
  \Bigr) J^\pm_\epsilon(t,x)\,, \qquad t \in \R\,, \quad x \in \R\,,
\]
which immediately imply \eqref{Sepsder} in view of the definition \eqref{waux2}
of $S_\epsilon$. On the other hand, if $0 < \epsilon' < \epsilon$, we have the identity
\begin{equation}\label{Sepsdiff}
  S_\epsilon(t,x) - S_{\epsilon'}(t,x) \,=\, \frac{1}{b}\,\biggl\{
  \log\biggl(1 + \frac{J^+_\epsilon(t,x) - J^+_{\epsilon'}(t,x)}{J^+_{\epsilon'}(t,x)}\biggr)
  -\log\biggl(1 + \frac{J^-_\epsilon(t,x) - J^-_{\epsilon'}(t,x)}{J^-_{\epsilon'}(t,x)}\biggr)
  \biggr\}\,,
\end{equation}
and proceeding as in the proof of Proposition~\ref{mu-prop1} we easily obtain,
if $|x| \le L(t)$, 
\[
  J^\pm_\epsilon(t,x) - J^\pm_{\epsilon'}(t,x) \,=\, \int_{(\epsilon',\epsilon]}
  e^{\mp xy/2}\,e^{t(by/2 + y^2/4)}\dd\nu_\pm(y) \,\le\, C^\pm_{\epsilon,\epsilon'}
  \,e^{\epsilon L(t) + tb\epsilon'/4} J^\pm_{\epsilon'}(t,x)\,,
\]
where $C^\pm_{\epsilon,\epsilon'} = \nu_\pm([\epsilon',\epsilon])/\nu_\pm([0,\epsilon'/2])$. 
If $L(t)/|t| \to 0$ as $t \to -\infty$, it follows that
\[
  \sup_{|x| \le L(t)}\,\biggl|\frac{J^\pm_\epsilon(t,x) - J^\pm_{\epsilon'}(t,x)}{J^\pm_{\epsilon'}
  (t,x)}\biggr| \,\le\, C_{\epsilon,\epsilon'}\,e^{\epsilon L(t) + tb\epsilon'/4}
  \,\xrightarrow[t \to -\infty]{}\, 0\,,
\]
which together with \eqref{Sepsdiff} implies the desired estimate \eqref{S2epsconv}. 
\end{proof}

\begin{proof}[\bf Proof of Lemma~\ref{seps-lem}.]
Let us define $\sigma_\epsilon(t) = S_\epsilon(t,0)$. For any $T > 0$, we observe
that
\begin{equation}\label{sigmadif}
  \lim_{t \to -\infty} \bigl(\sigma_\epsilon(t+T) - \sigma_\epsilon(t)\bigr) \,=\, 0\,,
  \qquad \hbox{which implies} \qquad \lim_{t \to -\infty}\frac{\sigma_\epsilon(t)}{t} \,=\,0\,.
\end{equation}
Indeed, if $0 < \epsilon' < \epsilon$, we have
\[
  \big|\sigma_\epsilon(t+T) - \sigma_\epsilon(t)\big| \,\le\,
  \big|\sigma_\epsilon(t+T) - \sigma_{\epsilon'}(t+T)\bigr| +
  \big|\sigma_{\epsilon'}(t+T) - \sigma_{\epsilon'}(t)\bigr| +
  \big|\sigma_{\epsilon'}(t) - \sigma_\epsilon(t)\bigr|\,.
\]
The middle term in the right-hand side can be estimated with the help of \eqref{Sepsder}\:
\[
  \big|\sigma_{\epsilon'}(t+T) - \sigma_{\epsilon'}(t)\bigr| \,\le\,
  \int_t^{t+T} |\partial_t S_{\epsilon'}(\tau,0)|\dd\tau \,\le\, 
  T\Bigl(\epsilon' + \frac{(\epsilon')^2}{2b}\Bigr) \,\le\,
  2T\epsilon'\,,
\]
where we used the fact that $\epsilon' < \epsilon < b$. The other terms are
controlled by \eqref{S2epsconv}, which gives
\[
  \limsup_{t \to -\infty} \big|\sigma_\epsilon(t+T) - \sigma_\epsilon(t)\big|
  \,\le\, 2T \epsilon'\,,
\]
and taking the limit $\epsilon' \to 0$ we obtain the first part of
\eqref{sigmadif}. The second claim follows by an elementary argument.

We now return to the shift function $s_\epsilon(t)$ defined by \eqref{sepsdef}. 
We have by \eqref{Sepsder}
\[
  |s_\epsilon(t) - \sigma_\epsilon(t)| \,=\, |S_\epsilon(t,s_\epsilon(t)) -
  S_\epsilon(t,0)| \,\le\, \frac{\epsilon}{b}\,|s_\epsilon(t)| \,\le\, 
  \frac{\epsilon}{b}\,|s_\epsilon(t) - \sigma_\epsilon(t)| + \frac{\epsilon}{b}\,
  |\sigma_\epsilon(t)|\,,
\]
so that
\[
  |s_\epsilon(t) - \sigma_\epsilon(t)| \,\le\, \frac{\epsilon}{b-\epsilon}
  |\sigma_\epsilon(t)|\,, \qquad |s_\epsilon(t)| \,\le\, \frac{b}{b-\epsilon}
  \,|\sigma_\epsilon(t)|\,.
\]
In view of \eqref{sigmadif}, it follows in particular that $s_\epsilon(t)/t \to
0$ as $t \to -\infty$, which is the first claim in \eqref{sepslim}. Moreover,
if $0 < \epsilon' < \epsilon$, we have
\begin{align*}
  |s_\epsilon(t) - s_{\epsilon'}(t)| \,&\le\, \bigl|S_\epsilon(t,s_\epsilon(t)) -
  S_\epsilon(t,s_{\epsilon'}(t))\bigr| + \bigl|S_\epsilon(t,s_{\epsilon'}(t)) -
  S_{\epsilon'}(t,s_{\epsilon'}(t))\bigr| \\
  \,&\le\, \frac{\epsilon}{b}\,|s_\epsilon(t) - s_{\epsilon'}(t)| \,+
  \sup_{|x| \le s_{\epsilon'}(t)} \,\bigl|S_\epsilon(t,x) - S_{\epsilon'}(t,x))
  \bigr|\,.
\end{align*}
Since $s_{\epsilon'}(t)/t \to 0$ as $t \to -\infty$, the last term in the
right-hand side converges to zero by \eqref{S2epsconv}. We deduce that
$|s_\epsilon(t) - s_{\epsilon'}(t)| \to 0$ too, which concludes the
proof of \eqref{sepslim}. 
\end{proof}

\bigskip\noindent
{\bf Thierry Gallay}\\
Institut Fourier, Universit\'e Grenoble Alpes, 100 rue des Maths, 38610 Gi\`eres,
France\\
Email\: {\tt Thierry.Gallay@univ-grenoble-alpes.fr}

\bigskip\noindent
{\bf Arnd Scheel}\\
School of Mathematics, University of Minnesota\\
127 Vincent Hall, 206 Church St.\thinspace SE, Minneapolis, MN 55455, USA\\
Email\: {\tt scheel@math.umn.edu}

\end{document}